\documentclass[letterpaper,12pt,leqno,oneside]{article}

\usepackage{float}
\usepackage{caption}
\usepackage{enumitem}
\usepackage{color}
\usepackage{ifsym}
\usepackage{bm}
\usepackage{exscale}
\usepackage{amsmath}
\usepackage{amsfonts}
\usepackage{wasysym}
\usepackage{stmaryrd}
\usepackage{amscd}
\usepackage{graphicx}
\usepackage{pagecolor,lipsum}
\usepackage{xcolor}
\usepackage{marvosym}
\usepackage{textcomp}
\usepackage{amsxtra}
\usepackage{amssymb}
\usepackage{theorem}
\usepackage{enumitem}
\usepackage[font=footnotesize]{caption}
\usepackage[sans]{dsfont}
\usepackage{mathtools}
\usepackage{esint}
\usepackage[hidelinks]{hyperref}
\usepackage[final]{epsfig}
\usepackage{eqnarray}
\usepackage{hyperref}
\usepackage{graphicx}
\let\texdisplaystyle\displaystyle
\def\displaytotextstyle{\textstyle\let\displaystyle\texdisplaystyle}

\newenvironment{talign*}
 {\let\displaystyle\displaytotextstyle\csname align*\endcsname}
 {\endalign}

\newtheorem{proposition}{Proposition}[subsection]
\newtheorem{definition}[proposition]{Definition}

\newtheorem{lemma}[proposition]{Lemma}

{\theorembodyfont{\rmfamily}\newtheorem{remark}[proposition]{Remark}}
\newtheorem{theorem}[proposition]{Theorem}
\newtheorem{corollary}[proposition]{Corollary}

{\theorembodyfont{\rmfamily}\newtheorem{example}[proposition]{Example}}

\newcommand{\mr}{\mathrm}

\newfont{\abc}{cmtt10 scaled 1200}

\def\ve{\varepsilon}

\def\ra{\rightarrow}
\def\cs{\symbol{35}}
\def\p{\partial}
\def\qed{\hfill $\Box$ \\}
\def\qeda{\hfill $\Box$}
\def\mm{\mbox}
\def\v{= \emptyset}
\def\n{\neq \emptyset}
\def\D{\mathbf{ID}}

\def\bp{\langle A \rangle}

\def\bs{\pmb{\square}}
\def\bx{\pmb{\boxplus}}
\def\si{$\mathcal{S}$}
\def\sima{\mathcal{S}}

\def\R{\mathbb{R}}

\def\Z{\mathbb{Z}}

\def\U{\mathbb{U}}
\def\D{\mathbf{ID}}
\def\P{\mathbb{P}}

\def\P{\mathbb{P}}

\def\U{\mathbb{U}}

\def\F{\mathbf{F}}
\def\G{\mathbf{G}}

\DeclareMathOperator{\su}{supp}
\DeclareMathOperator{\cu}{cut}

\def\ve{\varepsilon}

\def\ve{\varepsilon}

\def\ra{\rightarrow}
\def\cs{\symbol{35}}
\def\p{\partial}
\def\bp{\langle A \rangle}

\newcommand\diam{\mr{diam}}
\newcommand\dist{\mr{dist}}

\DeclareMathOperator{\scal}{scal}

\def\TP{\mathbb{TP}}

\def\XXint#1#2#3{{\setbox0=\hbox{$#1{#2#3}{\int}$ }
\vcenter{\hbox{$#2#3$ }}\kern-.6\wd0}}

\begin{document}

\vspace*{-0.6cm}
\begin{center}\Large{\bf{Scalar Curvature Splittings II: Removal of Singularities}}\\
\medskip
\large{\bf{Joachim Lohkamp}}\\
\smallskip
\end{center}
\noindent Mathematisches Institut, Universit\"at M\"unster, Einsteinstra\ss e 62, Germany\\
{\emph{e-mail: j.lohkamp@uni-muenster.de}}
{\footnotesize  {\center \tableofcontents}}

\setcounter{section}{1}
\renewcommand{\thesubsection}{\thesection}
\subsection{Introduction} \label{introduction}
We resume our discussion of general inductive scalar curvature splittings, involving singular minimal hypersurfaces, from the first part \cite{L1}. We also adopt the notations, concepts and results from  \cite{L1}. In this paper we establish splitting schemes with built-in regularizations to treat problems in scalar curvature geometry and general relativity in arbitrary dimensions.

\subsubsection{Statement and Discussion of Results} \label{sor}

The splitting approach involves the use of singular minimal hypersurfaces. In \cite{L1} we have seen that any singular compact area minimizer $(H^{n},d_H)$ in a $\scal > 0$-manifold $M^{n+1}$ admits a conformal deformation to some \textbf{minimal factor geometry} $(H^{n},d_{\sima})$ that shares many properties with $(H^{n},d_H)$ like the Ahlfors regularity, the validity of Poincar\'{e} inequalities and the presence of, in this case, $\scal>0$-tangent cones equipped with their minimal factor geometry. These geometries are amenable to surgery style arguments to eliminate singular sets stepwise. The basic building blocks we get are \emph{splitting with boundary} theorems. \\

\textbf{Theorem 1}  \, {\itshape  Let $H^{n}$, $n \ge 2$, be a compact area minimizer, with singular set $\Sigma$, in a $\scal >0$-manifold $M^{n+1}$.
Then there are arbitrarily small neighborhoods $U$ of $\Sigma$, so that $H \setminus U$ is \textbf{conformal} to some $\boldsymbol{\scal>0}$\textbf{-manifold} $(X_U,g_X)$ with \textbf{minimal} boundary $\p X_U=\p U$.}\\

For $\Sigma_H \v$ we choose $U \v$, $\p  U\v$. If $\Sigma \n$, then it suffices if $M^{n+1}$ has $\scal \ge 0$.  We actually get a smooth (generally non-complete) open $\scal>0$-manifold $Y_U^n$ that extends $X_U$, i.e. $\overline{X_U} \subset Y_U$, so that $\p X_U \subset Y_U$ is \emph{locally area minimizing}.

\textbf{Inductive Removal of Singularities} \,   Theorem 1 replaces the singular set $\Sigma_H$ of $(H^{n},d_{\sima})$ by a minimal boundary. The point is that although, in general, the boundary $\p X_U=\p U$ is also singular, its singular set has a \emph{lower dimension} than $\Sigma_H$. The minimality of $\p X_U$ then allows us to iteratively shift singular problems to lower dimensions before they disappear in dimension $7$. We describe some basic implementations of this idea.\\
The straightforward way to apply Theorem 1 is to use $\overline{X_U^n}$ as a replacement for $(H^{n},d_{\sima})$, for small $U$. Following the classical splitting recipe, we consider a local area minimizer $L^{n-1} \subset \overline{X_U^n}$. Since $X_U^n$ is not complete we observe that $L^{n-1}$ satisfies an additional constraint:  $L^{n-1}$ is an \emph{area minimizer with obstacle} $\p U$ that keeps $L^{n-1}$ in $\overline{X_U^n}$.  This is a standard situation in geometric measure theory, cf.~\cite[Th.\:1.20, Rm.\:1.22]{Gi}, but it is rather uncommon in the context of scalar curvature splittings.
\begin{itemize}[leftmargin=*]
\item However, since $\p U \subset Y_U$ is an unconstrained area minimizer by itself, the strict maximum principle \cite{Si} applies and it shows, componentwise,  that either $L^{n-1} \equiv \p U$ or $L^{n-1} \cap \p U \v$. This means $L^{n-1}$ always is an \emph{unconstrained} area minimizer in $Y_U^n$.
\item   $L^{n-1}$ could be singular, but when $L^{n-1}$ is compact we can apply Theorem 1 to $L^{n-1} \subset Y_U^n$,  in dimension $n-1$, since $Y_U^n$ is a \emph{smooth} $\scal>0$-manifold. In this sense Theorem 1 allows us to stay in the smooth category. Then we are ready for the next loop in the dimensional descent. Once we reach $n = 7$ all further area minimizers are regular.
\end{itemize}

\begin{figure}[htbp]
\centering
\includegraphics[scale=0.63]{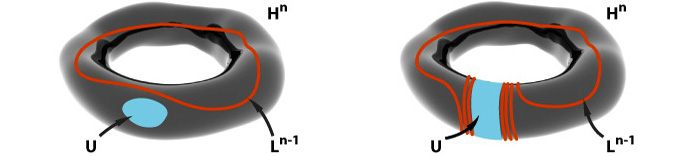}
\caption{\footnotesize $\Sigma_H$ has an empty intersection (i) with \emph{some} hypersurface in $\alpha[2] \cap \alpha[1] \cap [M]$ or not (ii).}
\label{fig:toco}
\end{figure}

\textbf{Applications }\,  We indicate some typical applications (addressed in separate accounts). We start from a compact $\scal>0$-manifold $M^{n+1}$ and classes $\alpha[1],\dots,\alpha[m] \in H^1(M^{n+1},\Z)$, $m  \le n-1$. The scheme above gives us the means to represent $\alpha[k]  \cap\cdots\cap \alpha[1] \cap [M]$, $k \le m$, by area minimizers $L^{n+1-k}$ in \emph{smooth} $\scal>0$-manifolds of dimension $n+2-k$, provided we broaden our playground and include \emph{complete spaces with periodic ends} and \emph{coarse homologies}, cf.~\cite{R} and \cite{L5} for a survey (also of underlying techniques and applications).\\
To explain the idea, we let $H^n$ represent $\alpha[1]  \cap [M]$ and apply Theorem 1 to get a conformal deformation to a $\scal>0$-geometry with minimal boundary $\p U$ for some neighborhood $U$ of $\Sigma_H$. When $L^{n-1}  =L^{n-1}(U) \equiv \p U$, the class  $\alpha[2] \cap \alpha[1] \cap [M]$ is trivial, since $L^{n-1}$ is null-homologous in $H^n \setminus U$. For non-trivial $\alpha[2] \cap \alpha[1] \cap [M]$, and small $U$, we have two possible scenarios for $L^{n-1}$, (i) and (ii), illustrated in Fig.\:\ref{fig:toco}.
\vspace{0.2cm}
\begin{enumerate}
\vspace{-0.3cm}
\item In many cases one can choose a setup so that the class contains hypersurfaces disjoint from $\Sigma$ and thus from $\overline{U}$, when $U$ is small enough. Then we find some representing \emph{compact} and unconstrained area minimizer $L^{n-1} \subset H^n \setminus \overline{U}$.  A typical situation where this happens is that where $M^{n+1}$ has \emph{geometrically large} components, like nearly flat torus summands. Examples are the \emph{positive mass/energy theorems} where we get such components from a compactification of asymptotically flat ends \cite[Ch.\:6]{L6}, \cite{L7}.
\item Otherwise the minimal factor geometry $(H^{n},d_{\sima})$ allows us to get  \emph{intrinsically complete} local area minimizers $L^{n-1}(U)$ asymptotically approaching $\p U$. For shrinking $U$, the $L^{n-1}(U)$ either approximate or represent $\alpha[2] \cap \alpha[1] \cap [M]$ in coarse homology, annihilating the infinite spinning of $L^{n-1}$ towards the compact $\p U$. We apply Theorem 1 to $\p U$ and periodically transfer the outcome to the end of $L^{n-1}$. (The deformation in Theorem 1 are assembled from local deformations and we can combine the transferred deformations  from $\p U$ with those on the compact part of $L^{n-1}$ away from these ends.)
   As a concrete application, we note that this matches the scheme in \cite[Ch.\:12]{GL1} and it shows that \emph{enlargeable manifolds} cannot admit $\scal>0$-metrics, cf.~\cite{L8}.
\end{enumerate}

Another way of applying this scheme is that of creating a smooth replacement $H_\circ^{n}$ for the singular $H^{n}$ that admits a $\scal >0$-metric. This does not change the scope of the method, but it fits into the layout of some classical arguments  for the lower dimensional cases.\\

\textbf{Variations and Extensions} \, We extend Theorem 1 in two directions. We start with the case of manifolds and area minimizers with already given boundaries. This scenario is used, for instance, in arguments in the style of classical comparison geometry, \cite{GL1}, Ch.\:12:\\

\textbf{Theorem 2} \, {\itshape  Let $H^{n}$ be a compact area minimizer with boundary $\p H\cap \Sigma \v$ in a $\scal > 0$-manifold $M^{n+1}$.  Then there are arbitrarily small neighborhoods $U \cap \p H \v$ of $\Sigma$, so that $H \setminus U$ is \textbf{conformal} to some $\boldsymbol{\scal>0}$\textbf{-manifold} $X_U$ with disjoint boundary components $\p X_0$ diffeomorphic to $\p H$ and  \textbf{minimal} $\p X_1$.}\\

The methods equally apply to larger classes of \emph{almost} minimizers, like  $\mu$-bubbles, levels sets of various geometric flows or horizons of black holes from Lorentzian geometry. The admissible classes $\cal{H}$ and $\cal{G}$ of almost minimizers are discussed and specified in \cite[Ch.\:1.2]{L1}  and \cite[Ch.\:3.1]{L4}. We have the following generalization of  Theorem 1.\\

\textbf{Theorem 3} \, {\itshape  Let $H^{n} \subset M^{n+1}$ be a compact almost minimizer with singular set $\Sigma_H$ with \si-adapted conformal Laplacian $L_H$.  Then  there are arbitrarily small neighborhoods $U$ of $\Sigma$, so that $H \setminus U$ is \textbf{conformal} to a $\boldsymbol{\scal>0}$\textbf{-manifold} $X_U$  with  \textbf{minimal}  boundary $\p X_U$.}\\

The \si-adaptedness of $L_H$ means that the principal eigenvalue of $\bp^{-2} \cdot L_H$ is positive and, thus,  $H$ can be conformally deformed into a $\scal >0$-minimal factor $(H,d_{\sima})$ by an eigenfunction of $\bp^{-2} \cdot L_H$. There is also a boundary version similar to Theorem 2.

\subsubsection{Overview of the Argument} \label{owa}

\begin{figure}[H]
\centering
\includegraphics[width=0.9\textwidth]{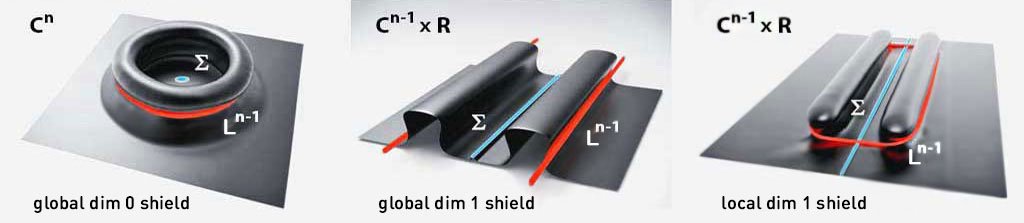}
\caption{\footnotesize In this illustration, the initial minimal factors $(C,d_{\sima})$ are entirely flat planes. The singularities are marked as bluish spots and lines. After bump deformations we get local area minimizers $L^{n-1}$ = the red rubber bands. In other words, the deformations shield (parts of) $\Sigma$ from intrusions of $L^{n-1}$.
 \normalsize}
\label{fig:a}
\end{figure}

\textbf{Basic Strategy} \, The deformations in our Theorems are similar to those in classical $\scal >0$-preserving surgeries \cite{GL2}, \cite{SY}: for a $\scal >0$-manifold $M$ and a compact submanifold $N \subset M$ of codimension $\ge 3$, there are small neighborhoods $U \subset V$ of $N$ so that  $M \setminus U$ can be conformally deformed within $V$ to a $\scal >0$-metric so that $\p U$ is a local area minimizer. This variant is taken from \cite[Ch.\:4]{L5}. For the subsequent  surgery/gluing processes in \cite{GL2}, \cite{SY}, which we do not need in our context, one appends (non-conformal) deformations of $V \setminus U$ to transform $\p U$ into a totally geodesic boundary.\\
In our case, the assumptions in our Theorems mean that the singular (almost) minimizer $H^n \subset M^{n+1}$ admits a conformal deformation to a still singular but well-behaved $\scal>0$-minimal factor geometry $(H,d_{\sima})$ from  \cite{L1}.  We take $(H,d_{\sima})$ as our initial $\scal>0$-geometry and let $\Sigma \subset (H,d_{\sima})$ play the r\^{o}le of $N \subset M$ in the classical smooth theory above. We construct small neighborhoods $U \subset V$ of $\Sigma$ and a conformal deformation,  this time of $d_{\sima}$, supported in $V \setminus \Sigma$, to another $\scal >0$-metric on $H \setminus U$ so that $\p U$ is a local area minimizer.
 Since $\Sigma$ can be rather complicated  and $H \setminus \Sigma$ degenerates towards $\Sigma$, the construction of this deformation is broken into simpler \emph{local bump} deformations of $(H,d_{\sima})$. Each of these local bumps creates a cavity oriented along some part of $\Sigma$. This gives rise to a local area minimizer $L$, stretched over the bump,  that bounds an open subset $L^+ \subset H$ containing a non-empty controllable subset $S_L \subset \Sigma$, see the right image of Fig.\:\ref{fig:a}. This is part of a bootstrapping where we use $L$ to show that the bump \emph{shields} $S_L$ against intersections not only with $L$ but with further area minimizers in $H$. The class of such minimizers is large enough to see that we can combine local bumps to deformations shielding all of $\Sigma$, proving the Theorems.

\vspace{0cm}
\vspace{-0.2cm}
\textbf{A Glance at  Technical Details} \, We work with the metric measure space $(H,d_{\sima},\mu_{\sima})$ associated to the metric completion $(H,d_{\sima})$  of $(H \setminus \Sigma,\Psi_H^{4/(n-2)} \cdot g_H)$, where $\Psi_H>0$ is a $C^{2,\alpha}$-(super)solution, for an $\alpha \in (0,1)$, with \emph{minimal growth} towards $\Sigma_H$, of
\small
\begin{equation}\label{weee}
L_{H, \lambda} (\phi):=L_H (\phi) - \lambda \cdot \bp^2 \cdot \phi := -\Delta \phi +\frac{n-2}{4 (n-1)} \cdot \scal_H \cdot \phi - \lambda \cdot \bp^2 \cdot \phi =0 \,\mm{ on } H \setminus \Sigma,
\end{equation}
\normalsize
for some \textbf{subcritical} eigenvalue $\lambda \in (0,\lambda^{\bp}_{H})$,  where $\lambda^{\bp}_{H}$ is the principal eigenvalue of ${\bp}^{-2} \cdot L_H$, to get a well-controlled geometry $(H,d_{\sima})$.  Equation (\ref{weee}) and the value of $\lambda$ remain \textbf{invariant under scalings}.  Under blow-up in any  $p \in \Sigma_H$,  $\Psi_H$ induces the unique$^\cs$ solution\footnote{\emph{unique$^\cs$}  is our abbreviation for \emph{unique up to multiplication by a positive constant}.} $\Psi_\circ>0$ of $L_{C, \lambda}  \, \phi=0$ with minimal growth towards $\Sigma_C$ on each of its tangent cones $C$  \cite[Th.\:1.5 and 1.6]{L1}. From this, the completion $(C,d_{\sima})$ of $(C, \Psi_\circ^{4/(n-2)} \cdot g_C)$  is again a cone and it is a $\scal>0$\emph{-tangent cone} of $(H,d_{\sima})$: for large $\tau \gg 1$, it locally approximates $(H, \tau \cdot d_{\sima})$ similarly to ordinary tangent cones for the original geometry $(H,\tau \cdot d_H)$. This allows us to construct local bumps on $(H,d_{\sima})$ from cone reduction arguments we indicate in the following.\\

\textbf{Dimension 7:   Isolated Singularities} \, We start in dimension $7$ since all almost minimizers are regular in lower dimensions. The singular set $\Sigma_H$ of a compact singular almost minimizer $H^7 \subset M^8$ contains at most finitely many singular points $p_1,...,p_m$ and in each $p_i \in \Sigma_H$ we have tangent cones $C_i^7$ approximating $ \tau \cdot H$ near $B_1(p_i) \subset \tau \cdot H$ for large $\tau \gg 1$. That is, for each $i$, there is a canonical $C^{2,\alpha}$-diffeomorphism $\D_i: B_1(0) \setminus B_r(0) \cap C_i^7 \ra  B_1(p_i) \setminus B_r(p_i)   \subset \tau \cdot H$, up to minor modifications near the boundary. $\D_i|_{B_1(0) \setminus B_r(0)}$ is nearly an isometry when $\tau$ is large enough, cf.~\cite[Ch.\:1.3.D.2.]{L1}, and this allows us the transfer local deformations between $C_i$ and $\tau \cdot H$.
 The $C_i^7$ are singular only in $0$ and there is some  $\lambda(H) \in (0,\lambda^{\bp}_{H})$ so that for $(\omega,r) \in  (\p B_1(0) \cap C_i^7)  \times \R^{>0} = C_i^7 \setminus \{0\}$, cf.~ \cite[Prop.\:4.6]{L4}:
\begin{equation}\label{wo}
\textstyle \Psi_\circ(\omega,r)= \psi(\omega) \cdot r^{\alpha_\circ},  \mm{ with } \textstyle{0 > \alpha_\circ \ge - (1- \sqrt{\frac 23}) \cdot \frac{7-2}{2}>  - \frac{7-2}{2}}
\end{equation}
and for the antagonist of $\Psi_\circ$, the solution $\Psi_\bullet$  with  \emph{maximal growth} towards  $\{0\}$:
\begin{equation}\label{wie}
\textstyle \Psi_\bullet(\omega,r)= \psi(\omega) \cdot r^{\alpha_\bullet},  \mm{ with } \textstyle{- \frac{7-2}{2} >- (1+ \sqrt{\frac 23}) \cdot \frac{7-2}{2} \ge \alpha_\bullet \ge - (7-2)}.
\end{equation}
For balls $B_\varrho(0) \subset C_i^7$, $\varrho >0$,   the mean curvature of $\p B_\varrho$ is negative\footnote{Our sign convention is that $\p B_1(0) \subset  (\R^n, g_{Eucl})$ has \emph{negative} mean curvature relative to $\R^n \setminus \overline{B_1(0)}$.}. The critical exponent $\beta$ for a conformal deformation $(r^\beta)^{4/7-2} \cdot g_{C_i^7}$ to flip the sign of the mean curvature is $\beta=- \frac{7-2}{2}$ where this deformation yields a $\scal>0$-cylinder making $\p B_\varrho$ totally geodesic. Thus the mean curvature of $\p B_ \varrho$ relative to $\Psi_\circ^{4/7-2} \cdot g_{C_i^7}$ is still \textbf{negative}, but relative to $\Psi_\bullet^{4/7-2} \cdot g_{C_i^7}$ it is \textbf{positive}.
We merge  these metrics to a new $\scal>0$-metric, the \textbf{bumped metric}  $g_\Omega= (\Phi_\Delta \cdot \Psi_\circ)^{4/7-2} \cdot g_{C_i^7}$, for some $C^{2,\alpha}$-regular $\Phi_\Delta \ge 1 $, cf.~the left image of Fig.\:\ref{fig:a}, so that:
\begin{itemize}[leftmargin=*]
\item there is a radius $r \in (0,1)$ so that $\su |\Phi_\Delta-1| \subset \overline{B_1(0)} \setminus B_r(0)$,
\item  there is a $\varrho \in (r,1)$ so that $L^{n-1}:=\p B_\varrho(0)$ is \textbf{locally area minimizing} relative to $g_\Omega$.
\end{itemize}
We call $\Delta_{C_i^7} := \big((\Phi_\Delta-1) \cdot \Psi_\circ\big)^{4/7-2} \cdot g_{C_i^7}$ a \textbf{local bump} pseudo metric\footnote{In this paper a pseudo metric is a symmetric \emph{positive semi-definite} bilinear form.}. Now we use the maps $\D_i$, for $\tau$ large enough, and transplant such local bumps from tangent cones to the balls $B_1(p_i)$ in $(H,\tau \cdot d_{\sima})$ adding the $\D_i$-push-forward of the pseudo metric from $C_i$ to $\tau \cdot H$ and get  Theorem 1 in dimension $7$.\\

\textbf{Dimensions $\boldsymbol \ge 8$:  Localization via Ahlfors regularity and Isoperimetry}\, As before we want use local bumps. The practical hurdle is that for a general area minimizer $H^n$, $n \ge 8$,  the singular set is no longer a discrete set.  A typical example is a product cone $C^{n-1} \times \R$ singular in $\{0\} \times \R$ with an $\R$-invariant bump pseudo metric that shields $\Sigma_{C^{n-1} \times \R}$, cf.~the middle of Fig.\:\ref{fig:a}.
These deformations have \emph{non-compact} support whereas tangent cones smoothly approximate $(H,\tau \cdot d_{\sima})$ only on subsets with \emph{compact closure} in $C^n \setminus \Sigma_C$. This obstructs the transfer of such deformations from $(C, d_{\sima})$ to $(H,\tau \cdot d_{\sima})$ as a whole.\\ We trim these pseudo metrics on $C^n \setminus \Sigma_C$ to deformations with compact support keeping $\scal>0$ to define a \textbf{local bump}. (This is indicated in the right image of Fig.\:\ref{fig:a}. The proper definition needs a bit more work illustrated in Fig.\:\ref{fig:localbump}. We give a separate overview in the beginning of Ch.\:\ref{lmc}.)
The trade-off is that the (support of) local bumps no longer topologically separate small neighborhoods of the singularities from large regular complementary parts of $H$. In view of the singular nature of the underlying space this raises the question whether such a trimming ruins the shielding effect of the bumps. \\
At this point the \textbf{Ahlfors regularity} and the \textbf{isoperimetric inequality} for $(H,d_{\sima},\mu_{\sima})$ \cite[Th.\:1.11]{L1} enter the game. They compensate this shortcoming and keep the area minimizer $L$ from entering the cavity surrounded by local bumps. This shielding property of local bumps is stable enough to survive the transfer from cones to the more general $H$.
We use a suitable covering scheme to place such bumps along $\Sigma_H$ while keeping their support disjoint.  As a direct consequence of this disjointness and the stable shielding property of local bumps the resulting global bump on $H$ shields all of $\Sigma_H$.\\

\textbf{Organization of the Paper} \, We focus on Theorem 3 and assume $\Sigma_H\n$.   Theorems 1 and 2 are special cases. \\
In Ch.\:\ref{wida}, we specify the admissible class of bumped metrics we use in this paper and we discuss results from \cite{L1} and \cite{L4} that we may apply to such metrics.\\
In Ch.\:\ref{dopc}, we construct simple bump pseudo metrics on product cones $\R^{n-k} \times C^k$ shielding the axis $\R^{n-k} \times \{0\}$. Then, in Ch.\:\ref{vcd}-\ref{localbump}, we apply localizing  modifications to this basic pseudo metric and assemble local bumps on $\R^{n-k} \times C^k$ shielding a cube $Q^{n-k}$ in $\R^{n-k} \times \{0\}$.\\
 In Ch.\:\ref{ptf}-\ref{trapro}, we prepare the transfer of local bumps to general (almost) minimizers and, in Ch.\:\ref{ifp} we place families of local bumps so that the union of the shielded subsets covers $\Sigma_H$.\\

\textbf{Acknowledgement}\, The author thanks Matthias Kemper for valuable suggestions which helped us to improve the exposition.

\subsubsection{Area Minimizers in $\Lambda$-Bumped Spaces} \label{wida}

We summarize some of the results from \cite{L2}-\cite{L4} and extend them to the class of controllably bumped metrics we work with.  Throughout the paper we exploit particularities of the conformal Laplacians $L_H$ and  ${\bp}^{-2} \cdot L_H$
for \emph{subcritical} but positive eigenvalues where we have a good control over the analysis near the singular set.\\
For the sake of consistency we oftentimes state results for all almost minimizers $H \in {\cal{G}}$. As in \cite{L1} this comes with the caveat that for regular $H \in {\cal{G}}^c$ or $H =\R^n \in {\cal{H}}^{\R}_n$ some statements become trivial/void because there is no singular set or no subcritical eigenfunction.\\

\textbf{Minimal Factors}\, The singular almost minimizer $H \in {\cal{G}}^c$ in  Theorem 3 has a positive principal eigenvalue  $\lambda^{\bp}_H>0$. This means that for any $\lambda < \lambda^{\bp}_H$ and $L_{H,\lambda} := L_H - \lambda \cdot \bp^2$ there are positive solutions of $L_{H,\lambda} \, \phi=0$, cf.~\cite{L3}. The coefficients are locally Lipschitz regular, hence any such solution is $C^{2,\alpha}$-regular for an $\alpha \in (0,1)$. For some subcritical eigenvalue $\lambda \in (0,\lambda^{\bp}_H)$, we set, following \cite{L1}:
\begin{definition} \emph{\textbf{(Minimal Factor Metrics)}}  \label{msge}  For $H \in {\cal{G}}_n$, let  $\Psi>0$ be a $C^{2,\alpha}$-supersolution  of $L_{H,\lambda} \phi =0$ on $H\setminus\Sigma_{H}$ so that in the case
\begin{itemize}
  \item  $H \in {\cal{G}}^c_n$: $\Psi$  is a solution in a neighborhood of $\Sigma$  with \textbf{minimal growth} towards $\Sigma$.
  \item $H \in {\cal{H}}^{\R}_n$:  $\Psi$ is a solution on $H\setminus \Sigma_{H}$ with \textbf{minimal growth} towards $\Sigma$.
\end{itemize}
We call the $\scal>0$-metrics $\Psi^{4/(n-2)} \cdot g_H$ on $H \setminus \Sigma$ the \textbf{minimal factor metrics}.
\end{definition}
For non-totally geodesic $H \in {\cal{H}}^{\R}_n$ there is a unique$^\cs$ such solution $\Psi_\circ$ of $L_{H,\lambda} \, \phi=0$ for subcritical $\lambda >0$ \cite[Theorem 3]{L3}. To any singular $H \in {\cal{G}}^c$ we assign one  fixed supersolution $\Psi_H$.\\

\textbf{Choice of Subcritical Eigenvalues}\, Now we select a subcritical eigenvalue $\lambda^\bullet_H$ of ${\bp}^{-2} \cdot L_H$ to gain the needed control over the growth rate of a supersolution $\Psi$, as in Def.\:\ref{msge}, towards $\Sigma_H$. To this end, we consider tangent cones $C^n$ obtained from blow-ups around a singular point of $H^n$. The blow-up process can be iterated in singular points outside $0$. After each blow-up the resulting minimal cone acquires an additional Euclidean factor. After ${n-k}$ blow-up steps we reach a singular cone $C^n= \R^{n-k} \times C^k \subset\R^{n+1}$, $k \ge 7$. When $\Sigma_{C^k} =\{0\}$ the blow-up process terminates.\\
From \cite[Theorem 3]{L4} we know that $\Psi$ induces the unique$^\cs$ solution  $\Psi_\circ=\Psi_\circ[C^n]$ of $(L_{C^n})_\lambda \phi =0$ on $C^n \setminus \Sigma_C$ with minimal growth towards $\Sigma_C$. This uniqueness implies that $\Psi_\circ$ shares the $\R^{n-k}$-translation symmetry with the underlying space $\R^{n-k} \times C^k$. This also means that $(\R^{n-k} \times C^k,d_{\sima})$ is again a cone invariant under $\R^{n-k}$-translation and this means that $\Psi_\circ|_{\{0\} \times C^k}$ satisfies the equation
$(L_{C^k,n})_\lambda \Psi_\circ(\omega,r)|_{\{0\} \times C^k}=0$ where we have set
\small
\begin{equation}\label{lck}
 L_{C^k,n}:=-\Delta - \frac{n-2}{4 (n-1)} \cdot |A|^2  \mm{ for }n \ge k \mm{ and } (L_{C^k,n})_\lambda  = L_{C^k,n}-\lambda  \cdot  \bp^2.
\end{equation}
\normalsize
$L_{C^k,n}$ is a dimensionally shifted conformal Laplacian on $C^k$. $\Delta$, $|A|^2=-\scal$ and $\bp^2$ are  defined relative to $C^k$. The dimensional shift comes from using $\frac{n-2}{4 (n-1)}$ in place of $\frac{k-2}{4 (k-1)}$.  From $\bp(\omega,\rho)= \bp(\omega,1) \cdot \rho^{-1}$ for $x=(\omega,\rho) \in S_{C^k} \setminus \Sigma_{C^k}  \times \R^{> 0}\cong C^k \setminus \{0\}$, where $S_{C^k}:= \p B_1(0) \cap C^k \setminus \Sigma_C$, we have
\small
\begin{equation}\label{pol}
(L_{C^k,n})_\lambda \, v= -\frac{\p^2 v}{\p r^2} - \frac{k-1}{r} \cdot \frac{\p v}{\p r} +\frac{1}{r^2} \cdot L_\lambda^\times \, v, \mm{ where }
\end{equation}
\vspace{-0.5cm}
\begin{equation}\label{ewg}{
L_\lambda^\times(v)(\omega):= - \Delta_{S_{C^k}} v(\omega) + \big(- \frac{k-2}{4 (k-1)} \cdot |A|^2(\omega,1) - \lambda \cdot\bp^2(\omega,1) \big)\cdot v(\omega)}.
\end{equation}
\normalsize
 We recall \cite[Prop.\:4.6]{L4}:

\begin{proposition}\emph{\textbf{(Growth Estimates)}} \label{dimshift}
Let $C^k \subset \R^{k+1}$ be a singular area minimizing cone and $n \ge k \ge 7$ where $C^k$ can also be singular outside $\{0\}$. Then there is a constant $\lambda^*_k >0$ such that for any $n \ge k$ and $\lambda \in (0,\lambda^*_k]$ there are two distinguished solutions of $(L_{C^k,n})_\lambda \phi =0$:
\begin{equation}\label{dizw}
\Psi_\circ[C^k](\omega,\rho)=\psi(\omega) \cdot \rho^{\alpha_\circ} \mm{ and }  \Psi_\bullet[C^k](\omega,\rho)=\psi(\omega) \cdot \rho^{\alpha_\bullet}, \mm{ for } (\omega,\rho) \in S_{C^k} \times \R^{>0}
\end{equation}
$\Psi_\circ$  and $\Psi_\bullet$ are the unique$^\cs$ solutions with minimal growth towards $\Sigma_{C^k}$ and $\widehat{\Sigma_{C^k}} \setminus \{0\}$ and we have some $\vartheta^*_\lambda>0$ depending only on $\lambda$ and $k$, but not on $n$, so that
\begin{equation}\label{abb}
\textstyle 0 > -\vartheta^*_\lambda >  \alpha_\circ \ge - (1- \sqrt{\frac 23}) \cdot \frac{k-2}{2}>  - \frac{k-2}{2} >- (1+ \sqrt{\frac 23}) \cdot \frac{k-2}{2} \ge \alpha_\bullet > \vartheta^*_\lambda -(k-2)
  \end{equation}
We interpret $\Psi_\circ$, $\Psi_\bullet$ as solutions $\Psi_\circ[\R^{n-k} \times C^k](x,y)=\Psi_\circ(x)$ and $\Psi_\bullet[\R^{n-k} \times C^k](x,y)=\Psi_\bullet(x)$ of $(L_{\R^{n-k} \times C^k})_\lambda \phi =0$ for $(x,y)\in\R^{n-k} \times C^k \setminus \Sigma_{\R^{n-k} \times C^k}$.
\end{proposition}
The spherical component $\psi(\omega)$ and the radial growth rates $\alpha_\circ$ of $\Psi_\circ$ and $\alpha_\bullet$ of  $\Psi_\bullet$ in  (\ref{dizw}) are related through a spherical eigenvalue equation \cite[Theorem 4.4]{L4} to is the (non-weighted) principal eigenvalue $\mu_{C^k,L^\times_\lambda}$ of the associated operator $L_\lambda^\times$ on $S_{C^k}$
for  $\psi(\omega) >0$ with  $L_\lambda^\times \, \psi =  \mu_{C^k,L^\times_\lambda} \cdot  \psi, \mm{ on } S_{C^k} \setminus \Sigma_{S_{C^k}}.$ In particular we have $\alpha_\circ + \alpha_\bullet=-(k-2)$.\\

We also recall from \cite[Lemma 3.9 and  Theorem 4.5]{L4} that for any tangent cone $C$ in some singular point of $H$, we have $\lambda^{\bp}_{H} \le \lambda^{\bp}_{C}$ and there is a constant $\widetilde{\lambda}_k^*>0$ so that $\lambda^{\bp}_{C^k} \ge \widetilde{\lambda}_k^*$ for any singular area minimizing cone $C^k$.
This suggests the following choice for a subcritical eigenvalue that respects all constraints and yields the uniform growth estimates of (\ref{abb}).
Let $H^{n} \in {\cal{G}}^c_n$ be a singular almost minimizer so that $\bp^{-2} \cdot L_H$ has principal eigenvalue $\lambda^{\bp}_{H}>0$. Then we define our \textbf{standard eigenvalue}
\begin{equation}\label{me}
\lambda^\bullet_H:= 1/2 \cdot \min\{\lambda^{\bp}_{H}, \lambda^*_7,..,\lambda^*_n, \widetilde{\lambda}_7^*,..,\widetilde{\lambda}_n^* \} >0.
\end{equation}

For the \textbf{rest of this paper}, on a given $H^n$, we choose $\lambda:=\lambda^\bullet_H>0$ and a fixed supersolution of $(L_{H^n})_\lambda \phi =0$ as in Def.\:\ref{msge}, also written as $\Psi_\circ$, and we use this $\lambda$ and the induced unique$^\cs$ solutions $\Psi_\circ$ of $(L_{C^n})_\lambda \phi =0$ on $C^n \setminus \Sigma_C$  on any of its (iterated) tangent cones $C^n$. Since we keep these choices unchanged, we generally do not mention $\lambda$ and $\Psi_\circ$  explicitly.\\

\textbf{Controlled Bumps}\, The metric completion $(H,d_{\sima})$ of $(H \setminus \Sigma,\Psi_\circ^{4/(n-2)} \cdot g_H)$ can be augmented to a metric measure space  $(H,d_{\sima},\mu_{\sima})$. As indicated in Ch.\:\ref{owa} above, our plan is to add disjointly supported local bump  deformations to $(H,d_{\sima})$. This means that all bump deformations in this paper remain \textbf{$\boldsymbol{\Lambda}$-equivalent}  to  $(H,d_{\sima})$, in the sense of $\boldsymbol{\Lambda}$\textbf{-bumped spaces} in the following result.  From this, many of the results for $(H,d_{\sima},\mu_{\sima})$, in \cite{L1}, remain valid for the bumped cases. We define, for some fixed bumping constant $\Lambda \ge 1$, the class of admissible $\Lambda$-bounded conformal deformation functions
\begin{equation}\label{lam}
{\boldsymbol{{\cal{C}}[H^n,\Lambda]}}:= \{f:H^n \ra [1,\Lambda^{(n-2)/2}]\,|\, f \in C^0(H^n)  \cap C^{2,\alpha}(H^n \setminus \Sigma_H), \mm{ for any } \alpha \in (0,1)\}.
\end{equation}

In the following we simplify the notations and only use the index $\Lambda$ for quantities actually depending on $\G \in{\cal{C}}[H,\Lambda]$ since the estimates  only depend on $\Lambda$.

\begin{theorem} \emph{\textbf{($\boldsymbol{\Lambda}$-Bumped Spaces)}}  \label{lbp}   For $\Lambda \ge 1$, $H^{n} \in {\cal{G}}_n$, and $\G \in{\cal{C}}[H,\Lambda]$ we have:
\begin{itemize}[leftmargin=*]
  \item $(\G \cdot \Psi_\circ)^{4/(n-2)} \cdot g_H$ on $H \setminus \Sigma$ can be completed to a geodesic metric space $(X^n,d_X)$. The space $(X^n,d_X)$ is  homeomorphic to $(H,d_H)$ with singular set $\Sigma_X \cong \Sigma_H$ of Hausdorff dimension $\le n-7$ relative to $(X^n,d_X)$. From this we can write $(X^n,d_X)$ as $(H,d^\Lambda_{\sima})$.
  \item The volume $(\G  \cdot \Psi_\circ)^{2 \cdot n/(n-2)}\cdot \mu_H$ on $H \setminus \Sigma_H$ extends to a Borel measure $\mu^\Lambda_{\sima}(E)$ defined by  $\int_{E \setminus \Sigma_H} (\G \cdot \Psi_\circ)^{2 \cdot n/(n-2)}\cdot d\mu_H$,  for any Borel set $E \subset H$,
and we get the metric measure space $(H,d^\Lambda_{\sima},\mu^\Lambda_{\sima})$ and we call it a $\boldsymbol{\Lambda}$\textbf{-bumped space} originating from $(H,d_{\sima},\mu_{\sima})$
  \item  $(H,d^\Lambda_{\sima},\mu^\Lambda_{\sima})$  is  \textbf{Ahlfors $\boldsymbol{n}$-regular}: there are  constants $A,B(H,\Psi_\circ,\Lambda)>0$ so that:
\begin{equation}\label{ahl}
 A \cdot r^n \le \mu^\Lambda_{\sima}(B_r(q),d_{\sima}) \le B \cdot r^n, \mm{ for any } r \in [0, \diam(H,d^\Lambda_{\sima})) \mm{ and any } q\in H.
\end{equation}
For $H \in {\cal{H}}^{\R}_n$ the constants only depend on the dimension and $\Lambda$.
\item For any $H \in {\cal{G}}$, there is constant $C_0(H,\Psi_\circ,\Lambda)>0$, depending only on $n$ and  $\Lambda$ for $H \in {\cal{H}}^{\R}_n$, so that for any ball $B\subset (H,d_\sima)$ we have a \textbf{Poincar\'{e} inequality}: for any function $u$ on $H$, integrable on bounded balls, and any upper gradient $w$ of $u$:
\begin{equation}\label{poinm0}
\textstyle \int_B |u-u_B| \,  d \mu^\Lambda_{\sima} \le C_0 \cdot \diam(B) \cdot \int_{B} w \, d \mu^\Lambda_{\sima}.
\end{equation}
\end{itemize}
\end{theorem}
The assertions readily follow from their counterparts for minimal factor geometries in \cite{L1} and the $\Lambda$-equivalence of the original $(H,d_{\sima},\mu_{\sima})$ and the bumped spaces $(H,d^\Lambda_{\sima},\mu^\Lambda_{\sima})$. From  (\ref{ahl})  and (\ref{poinm0}) we also get an isoperimetric inequality  $(H,d_{\sima},\mu_{\sima})$  and the following growth estimates. Alternatively, one can use the  $\Lambda$-equivalence of the metrics to get such applications from the estimates we have for $(H,d_{\sima},\mu_{\sima})$ from \cite[Cor.\:3.16 and Prop.\:3.17]{L1}.

\begin{corollary}\label{grr} For $(H,d^\Lambda_{\sima},\mu^\Lambda_{\sima})$, some open  subset $\Omega \subset H$ and an oriented minimal boundary $L^{n-1} \subset \Omega$ bounding an open set $L^+ \subset \Omega$, there are  $\kappa, \kappa^+(H,\Psi_\circ,\Lambda)    >0$, so that for $p \in L$:
\begin{equation}\label{est}
\kappa \cdot r^{n-1}  \le \mu^{\Lambda, n-1}_{\sima}(L \cap B_r(p)) \,\mm{ and }\, \kappa^+\cdot r^n  \le  \mu^{\Lambda}_{\sima}(L^+ \cap B_r(p)),
\end{equation}
for $r \in [0, (A/B)^{1/n} \cdot \dist(p,\p \Omega)/4)$ and where $0<A<B$ are the Ahlfors constants. For $H \in {\cal{H}}^{\R}_n$, $\kappa, \kappa^+  >0$ depend only on $n$ and $\Lambda$.
\end{corollary}

From  (\ref{ahl})  and (\ref{poinm0}) the BV-theory of Ambrosio \cite{A} and Miranda \cite{M} applies to  $(H,d^\Lambda_{\sima},\mu^\Lambda_{\sima})$. The coarea formula \cite[Prop.\:4.2]{M} and the Ahlfors $n$-regularity yield growth estimates for $\mu^{\Lambda, n-1}_{\sima}(\p B_r(q))$ in terms of a density result for the Lebesgue measure $\mu_1$ on $\R$, cf.~\cite[Ch.\:5.2]{M}:
\begin{corollary}\label{ctrl} For any $q \in H$, we define the set of radii where $\mu^{\Lambda, n-1}_{\sima}(\p B_r(q))$ exceeds $c \cdot r^{n-1}$ by  $J(q,a,b,c):= \{r \in [a,b] \,|\, \mu^{\Lambda, n-1}_{\sima}(\p B_r(q)) > c \cdot  r^{n-1}\}$, $0 \le a <b$ and $c>0$. We have
\begin{equation}\label{jes}
\mu_1(J(q,a,b,c)) \le \big(B \cdot b - A \cdot  a + (n-1) \cdot B \cdot (b-a)\big)/c.
\end{equation}
\end{corollary}
That is, for large $c>0$ it is increasingly likely that $\mu^{\Lambda, n-1}_{\sima}(\p B_r(q)) \le c \cdot r^{n-1}$ since, choosing a small $\ve>0$ and  $c=n \cdot B/\ve$, we have $\mu_1(J(q,0,b,c)) \le \ve \cdot b = \ve \cdot \mu_1([0,b])$. Using balls with radii $\notin J(q,0,b,c)$ we may exploit the low Hausdorff dimension of $\Sigma \subset (H,d^\Lambda_{\sima})$ to get:
\begin{corollary}\label{ctrl2} For $H^{n} \in {\cal{G}}_n$ and any $\ve >0$ and $R >0$ there is a locally finite cover of $\Sigma_H$ by balls $B_{r_i}(p_i)$, $p_i \in \Sigma_H$, $r_i \in (0,R)$, $i \in I$,  in $(H,d^\Lambda_{\sima},\mu^\Lambda_{\sima})$ so that for $c \ge 2 \cdot  n \cdot B$:
\begin{equation}\label{covha}
\textstyle \sum_{i \in I}\mu^{\Lambda, n-1}_{\sima}(\p B_{r_i}(p_i))   \le c \cdot \sum_{i \in I} r_i^{n-1}  < \ve.
\end{equation}
\end{corollary}
Similarly, \ref{ctrl} shows the otherwise non-obvious finiteness of the perimeter when we look for area minimizers in $(H,d^\Lambda_{\sima},\mu^\Lambda_{\sima})$.
The lower semi-continuity of perimeters \cite[Prop.\:3.6]{M} and the compactness of the BV-function space in the $L^1$-function space  \cite[Prop.\:3.7]{M} yield the existence of oriented minimal boundaries in $(H,d^\Lambda_{\sima},\mu^\Lambda_{\sima})$ as in \cite[Theorem 1.20]{Gi}.

\begin{proposition}\emph{\textbf{(Plateau Problems in $(H,d^\Lambda_{\sima},\mu^\Lambda_{\sima})$)}}\label{pla} Let $\Omega \subset H$ be a bounded open and orientable
set and let $A \subset H$ be a possibly empty \textbf{Caccioppoli set}, i.e., a set of locally finite perimeter. Then there exists a set $E\subset H$ coinciding with
$A$ outside $\Omega$ and such that
\begin{equation}\label{mininh}
\mu^{\Lambda, n-1}_{\sima}(\p E \cap \Omega) \le \mu^{\Lambda, n-1}_{\sima}(\p F \cap \Omega)
\end{equation}
for any Borel set $F  \subset H$ with $F = A$ outside $\Omega$. We say that $E$ has $\Omega$-boundary value $A$, and for open $E \cap \Omega$, we call
 $\p E$ a \textbf{Plateau solution} with $\boldsymbol\Omega$\textbf{-boundary value} $\boldsymbol{\p A}$.
\end{proposition}

\setcounter{section}{2}
\renewcommand{\thesubsection}{\thesection}
\subsection{Constructions on Cones} \label{lmc}

In this chapter we construct local bumps on product cones. We start with a more detailed description of the key steps of these deformations and their intended use.

\begin{enumerate}[label=(\arabic*),leftmargin=*]
\item \textbf{Product Bumps}\,  In Ch.\:\ref{dopc} we are starting from some singular cone $C^n=\R^{n-k} \times C^k$ with its $\scal>0$-minimal factor metric  $\Psi_\circ^{4/(n-2)} \cdot g_{C^n}$. We emphasize that $C^k$ can also be singular outside $\{0\}$.  We make an elementary deformation $(\F \cdot \Psi_\circ)^{4/(n-2)} \cdot g_{C^n}$, with $\F(x,y)=f(|y|)$, for some smooth  $f \ge 1$ with $\su |f-1| \subset [r,1]$, for some $1>r>0$.  $(\F \cdot \Psi_\circ)^{4/(n-2)} \cdot g_{C^n}$ keeps $\scal>0$ and it turns $\R^{n-k} \times \p B_\varrho(0)$, for some $\varrho \in (r,1)$, into a local area minimizer.  $(\F \cdot \Psi_\circ)^{4/(n-2)} \cdot g_{C^n}$ is the source material for all further steps where we trim and rearrange this deformation to build larger bump structures.
\item \textbf{Regular Support}\,  For the first trimming process, in Ch.\:\ref{vcd}, we define a cut-off function $\phi_\beta:C \setminus \Sigma_C \ra [0,1]$ with  $\phi_\beta =0$ on $\delta_{\bp}^{-1}(\R^{\le \beta})$ and $\phi_\beta =1$ on $\delta_{\bp}^{-1}(\R^{\ge \delta(\beta)})$, for some $\delta(\beta) \ra 0$, for $\beta \ra0$.
    We set $\F_\beta=\phi_\beta \cdot (\F-1) +1$ for usage in the new metric $(\F_\beta  \cdot \Psi_\circ)^{4/(n-2)} \cdot g_{C^n}$ and observe that $\su |\F_\beta-1| \cap \Sigma_C \v$.  We think of $\delta_{\bp}^{-1}(\R^{\le \beta})$, where we have trimmed $\F-1$ to zero, as a \emph{tunnel} that surrounds the singularities. This tunnel passes (and removes) part of the support of $(\F  \cdot \Psi_\circ)^{4/(n-2)} \cdot g_{C^n}$ reaching from $\R^{n-k} \times \p B_1(0)$ to $\R^{n-k} \times \p B_r(0)$, cf.~Fig.\:\ref{fig:tunnel}.  We choose $\delta_{\bp}$ in place of the metric distance since it commutes with the convergence of the underlying spaces as needed to transfer the bumps between different spaces and it gives us a sufficiently good control to carry this out keeping $\scal(\F_\beta  \cdot \Psi_\circ)^{4/(n-2)} \cdot g_{C^n}>0$.
\item \textbf{Tightness of Regular Bumps}\,  For $(\F_\beta  \cdot \Psi_\circ)^{4/(n-2)} \cdot g_{C^n}$ the boundary $\R^{n-k} \times \p B_\varrho(0)$ is no longer minimal, but in Ch.\:\ref{rol} we show that for small $\beta>0$, there is an area minimizer ${\cal{T}}^{\beta}$, in this metric which is Hausdorff close to $\R^{n-k} \times \p B_\varrho(0)$. In particular, ${\cal{T}}^{\beta} \cap \delta_{\bp}^{-1}(\R^{\le \beta})$
    remains close to $\delta_{\bp}^{-1}(\R^{\le \beta}) \cap \R^{n-k} \times \p B_\varrho(0)$ and we will see that this is true for any area minimizer relative to $(\F_\beta  \cdot \Psi_\circ)^{4/(n-2)} \cdot g_{C^n}$ within a certain threshold distance from $\R^{n-k} \times \p B_\varrho(0)$. In this sense, the tunnel is \emph{tight} and  $(\F_\beta  \cdot \Psi_\circ)^{4/(n-2)} \cdot g_{C^n}$  still \emph{shields} a neighborhood of $\R^{n-k} \times \{0\}$ against intrusions of these area minimizers.
\item \textbf{Compact Support}\, In Ch.\:\ref{localbump}  we apply an additional cut-off, now along the $\R^{n-k}$-axis, to $(\F_\beta  \cdot \Psi_\circ)^{4/(n-2)} \cdot g_{C^n}$, to make the support compact, cf.~the right image of Fig.\:\ref{fig:a}. We call the resulting deformation of $\Psi_\circ^{4/(n-2)} \cdot g_{C^n}$ an \emph{(n-k)-bump element} and they are the elementary particles of our bump construction. In this step we can preserve $\scal>0$, but it is difficult to prevent minimizers from shrinking along the $\R^{n-k}$-axis.
\item \textbf{Tight Configurations}\,  To resolve this issue, we inductively place ($m$)-bump elements along the ($m$)-dimensional faces of the unit cube $Q_1^{n-k}  \times \{0\}\subset \R^{n-k} \times C^k$. (We get an ($m$)-bump element writing $C^n=\R^m \times (\R^{n-k-m} \times C^k)$.) In this process, the cut-off ends of ($m$)-elements belong to subsets shielded  from ($q$)-elements for $q <m$, cf.~Fig.\:\ref{fig:localbump}. The tightness (3) ensures that these bump configurations shield $Q_1^{n-k}$ against intrusion of area minimizers. We call them \emph{local bumps}. They deform $\Psi_\circ^{4/(n-2)} \cdot g_{C^n}$ to $\scal>0$-metrics containing open sets $U \supset Q_1^{n-k} \times \{0\}$ with locally minimal boundary $\p U$.
\item \textbf{Stable Shielding}\,  The bump elements in local bumps (and in all larger configurations we build until we reach global bumps) have disjoint supports, but the generated minimizing boundaries intersect the support of other bump elements and they change whenever we add further bump elements. In general, any such addition produces several new locally minimizing boundaries. Some of them jump far away and they become useless for our purposes, but the central tightness result in Ch.\:\ref{rol} also shows that there is always at least one minimizing boundary only slightly perturbed compared to the initial one.
\end{enumerate}

\subsubsection{Deformations on Product Cones} \label{dopc}

The main building blocks for our deformations are blends of the two extremal metrics associated to  $\Psi_\circ$ and $\Psi_\bullet$ on $\R^{n-k} \times C^k$:
\begin{equation}\label{extr}
g^\circ_{C^n}:=\Psi^{4/(n-2)}_\circ  \cdot g_{C^n} \, \mm{ and } \, g^{\bullet}_{C^n}:=\Psi^{4/(n-2)}_\bullet  \cdot g_{C^n}.
\end{equation}
We compute the mean curvature of distance tubes of $\R^{n-k} \times \{0\}$ relative to these two metrics.
\begin{lemma} \emph{\textbf{(Bending Effects)}} \label{etl} \, We consider $u= f(\omega) \cdot r^\beta$, for  $\omega \in S_{C^k}:= \p B_1(0)
\cap C^k$,  $r(x):=\dist(x, \R^{n-k} \times \{0\})$, for $x \in \R^{n-k} \times C^k$ and some positive $f \in C^{2,\alpha}(S_{C^k} \setminus \Sigma_{S_{C^k}},\R)$. The \textbf{mean curvature} $tr\, A_{\mathcal{T}_\varrho}$ in regular points of the boundary $\mathcal{T}_\varrho$ of distance tubes $\R^{n-k} \times B_\varrho(0)
\cap C^k \subset C^n$ of $\R^{n-k} \times \{0\}$,  of radius $\varrho >0$ relative to $g_{C^n}$,  is given by
\begin{equation}\label{te}
u^{4/(n-2)} \cdot tr\, A_{\mathcal{T}_\varrho}(u^{4/(n-2)} \cdot g_{C^n}) =  - \left( (k-1) +  2 \cdot (n-1) \cdot  \beta/(n-2)\right) \cdot \varrho^{-1}.
\end{equation}
In particular,  for the two $\scal >0$-metrics $g^\circ_{C^n}$ and  $g^{\bullet}_{C^n}$ on $C^n \setminus \Sigma_{C^n}$, we get for any $\varrho>0$:
$tr\, A_{\mathcal{T}_\varrho}(g^\circ_{C^n})  <0$ and $tr\, A_{\mathcal{T}_\varrho}(g^{\bullet}_{C^n})  >0$.
For the superposition, we get some $\varrho_0(\alpha_\circ ,n,k)>0$ which continuously depends on $\alpha_\circ$ so that
\begin{equation}\label{m4}
 tr\, A_{\mathcal{T}_\varrho}\big((\Psi_\circ + \Psi_\bullet)^{4/(n-2)} \cdot g_{C^n}\big) >0 \mm{ for } \varrho < \varrho_0, =0 \mm{ for } \varrho = \varrho_0 \mm{ and } <0  \mm{ for } \varrho > \varrho_0.
\end{equation}
$\mathcal{T}_{\varrho_0}$ is an area minimizer in the metric measure space associated to the completion of $(C^n \setminus \Sigma_{C^n},(\Psi_\circ + \Psi_\bullet)^{4/(n-2)} \cdot g_{C^n})$. For any bounded open set $\Omega \subset C^n$, $\mathcal{T}_{\varrho_0} \cap \Omega$ is the \textbf{unique} oriented minimal boundary with $\Omega$-boundary value $\mathcal{T}_{\varrho_0}$ bounding an open set containing $\R^{n-k} \times \{0\}$.
\end{lemma}

\textbf{Proof} \,  The second fundamental form $A_L(g)$ of a submanifold $L$ with respect to a metric $g$ transforms under conformal deformations $g \ra u^{4/(n-2)} \cdot g$ as follows \cite[1.163, p.\:60]{Be}:
\begin{equation}\label{aaw}
\textstyle A_L(u^{4/(n-2)} \cdot g)(v,w) = A_L(g)(v,w) - \frac{2}{n-2}\cdot {\cal{N}} (\nabla u / u) \cdot g(v,w),
\end{equation}
where ${\cal{N}} (\nabla u / u)$ is the normal component of $\nabla u / u$.   $C^k \subset \R^{k+1}$ is a cone and $\R^{n-k} \times \{z\} \subset \mathcal{T}_\varrho$, $z \in C^k$, is totally geodesic. Thus we have $tr\,A_{\mathcal{T}_\varrho}( g ) = -(k-1)/\varrho$.  The trace of the second summand is $2/(n-2)  \cdot \beta/\varrho$ multiplied  by $n-1$. This yields:
{\small \[\Psi^{4/(n-2)}_\circ  \cdot tr\, A_{\mathcal{T}_\varrho}(\Psi^{4/(n-2)}_\circ \cdot g_{C^n}) =
- \left( (k-1) +   2 \cdot (n-1) \cdot  \alpha_\circ/(n-2)\right) \cdot \varrho^{-1}< -  (k-2)/2 \cdot \varrho^{-1}<0,\]
\[\Psi^{4/(n-2)}_{\bullet} \cdot tr\, A_{\mathcal{T}_\varrho}(\Psi^{4/(n-2)}_{\bullet} \cdot g_{C^n}) =
- \left( (k-1) +   2 \cdot (n-1) \cdot  \alpha_\bullet/(n-2)\right) \cdot \varrho^{-1}> (k-2)/2 \cdot \varrho^{-1}>0,\]}
since $k \ge 7$ and, in particular, $k\ge 3$. For $(\Psi_\circ + \Psi_\bullet)^{4/(n-2)} \cdot g_{C^n}$ we get
\begin{equation}\label{sk}
(\Psi_\circ + \Psi_\bullet)^{4/(n-2)} \cdot tr\, A_{\mathcal{T}_\varrho}((\Psi_\circ + \Psi_\bullet)^{4/(n-2)} \cdot g_{C^n})=
\end{equation}
\vspace{-0.8cm}
\begin{equation}\label{sl}
 -(k-1)/\varrho +2(n-1)/(n-2) \cdot (\alpha_\circ \cdot \varrho^{\alpha_\circ-1} +  \alpha_\bullet  \cdot \varrho^{\alpha_\bullet-1})/( \varrho^{\alpha_\circ} +  \varrho^{\alpha_\bullet}).
\end{equation}
The mean curvature is constant on each $\mathcal{T}_\varrho$. To understand its sign, we recall that
\begin{itemize}
  \item $\alpha_\circ + \alpha_\bullet=-(k-2)$ and $2 \cdot \alpha_\circ+1+(k-2)>0$,
  \item $|2(n-1)/(n-2) \cdot \alpha_\circ|<k-1$ and  $|2(n-1)/(n-2) \cdot \alpha_\bullet|>k-1$.
\end{itemize}
We multiply (\ref{sl}) by $( \varrho^{\alpha_\circ} +  \varrho^{\alpha_\bullet})$ and $\varrho^{\alpha_\circ+1+(k-2)}$ and get:
\[ -((k-1) +2(n-1)/(n-2) \cdot \alpha_\circ) \cdot \varrho^{2 \cdot \alpha_\circ+1+(k-2)}-((k-1) +2(n-1)/(n-2) \cdot \alpha_\bullet).\]
From this we see that there is a unique $\varrho_0(\alpha_\circ ,n,k)>0$ with
\[tr\, A_{\mathcal{T}_\varrho}\big((\Psi_\circ + \Psi_\bullet)^{4/(n-2)} \cdot g_{C^n}\big) >0  \mm{ for $\varrho < \varrho_0$,  $=0$  for  $\varrho = \varrho_0$ and $<0$ for $\varrho > \varrho_0$}.\]

Hence, the radial projection in the $C^k$-factor $\pi_{s,\varrho}: \mathcal{T}_s \ra \mathcal{T}_{\varrho_0}$ strictly decreases the hypersurface area element for any $\varrho \neq \varrho_0$. From this we see that $\mathcal{T}_{\varrho_0}$ is the unique oriented minimal boundary with boundary value $\mathcal{T}_{\varrho_0}$ bounding an open set containing $\R^{n-k} \times \{0\}$.\\
This carries over to the metric completion of $(C^n \setminus \Sigma_{C^n},(\Psi_\circ + \Psi_\bullet)^{4/(n-2)} \cdot g_{C^n})$ since the Hausdorff dimension of $\Sigma_C \setminus \R^{n-k} \times \{0\}$ is still $\le n-7$ from Theorem \ref{lbp}.   \qed

We use these estimates in our basic bump deformation of the minimal factor geometry $(C^n,d_{\sima})$. The construction uses two particularities of $(C^n,d_{\sima})$ while we approach $\R^{n-k} \times \{0\}$.
\begin{itemize}[leftmargin=*]
\item A first transition from $\Psi_\circ$ to $\Psi_\bullet$ uses, in the guise of $tr\, A_{\mathcal{T}_\varrho}(g^{\bullet}_{C^n})  >0$, that the codimension of the singular set is $\ge 3$ matching the condition in the $\scal >0$-preserving surgeries of \cite{GL2}, \cite{SY}. This deformation yields a local area minimizer $\mathcal{T}_\varrho$ with $\mathcal{T}_\varrho \cap \R^{n-k}  \times \{0\}$.
\item Then we append a second transition back from $\Psi_\bullet$ to $\Psi_\circ$ and return to the geometry $(C,d_{\sima})$ near $\R^{n-k} \times \{0\} \subset \Sigma_C$. This step has no classical counterpart in the mentioned $\scal >0$-preserving surgeries. It uses the scaling invariance of the $\scal>0$-cone $(C,d_{\sima})$.
\end{itemize}

\begin{proposition} \emph{\textbf{(Basic Deformations)}}\label{el}\, For any singular cone $C^n= \R^{n-k} \times C^k \in {\cal{SC}}_n$ there are a smooth function $\mathbf{f}[C]: \R^{>0} \ra \R^{\ge 1}$ and a \textbf{core radius} $r[n] \in (0,1/100)$,  such that for $z=(x,y) \in C^n \setminus \Sigma_{C^n}= \R^{n-k} \times (C^k \setminus \Sigma_{C^k})$ and $\F(x,y):= \mathbf{f}(|y|)$, the metric
\[g^{n,n-k}_{C^n} := (\F \cdot \Psi_\circ)^{4/(n-2)} \cdot g_{C^n}=(\mathbf{f}(|y|) \cdot \Psi_\circ(x,y))^{4/(n-2)} \cdot g_{C^n} \, \mm{ has the following properties: }\]
\begin{enumerate}
\item $\F=\F[C]$ depends continuously on $C \in \mathcal{SC}_{n}$:  for a converging sequence $C_i \ra C$ in $\cal{SC}$ we have a compact   $C^{2,\alpha}$-convergence $\F[C_i] \circ \D \ra \F[C]$ on $C \setminus \Sigma_C$.
\item  $\su |\mathbf{f}-1| \subset [r,1]$,  $\mathbf{f} \le \Lambda^{(n-2)/2 }$, for some $\Lambda[n] \ge 1$, i.e. $\F \in{\cal{C}}[H,\Lambda[n]]$,
\item   $L_{C^n} (\F \cdot \Psi_\circ) \ge \lambda/4 \cdot \bp^2 \cdot \F \cdot \Psi_\circ$
and, in particular, $\scal(g^{n,n-k}_{C^n}) >0$,
\item there are constants $\vartheta_i[C]\in (r,1)$, $i=1,2$, with $\vartheta_1 <   \vartheta_2$  and some $\varrho \in (\vartheta_1,\vartheta_2)$, all depending continuously on $C$, so that   $\mathcal{T}_\varrho = \R^{n-k} \times \p B_\varrho(0)$ is a \textbf{local area minimizer}.
   For any bounded open set $\Omega \subset C^n$, $\mathcal{T}_{\varrho} \cap \Omega$  is the \textbf{unique} oriented minimal boundary of some Caccioppoli set in $\R^{n-k} \times \overline{B_{\vartheta_4}(0)} \setminus B_{\vartheta_3}(0)$ with $\Omega$-boundary value $\mathcal{T}_{\varrho}$.
\end{enumerate}
We call $\Delta^{n, n-k}:=\big((\F-1) \cdot \Psi_\circ\big)^{4/(n-2)} \cdot g_{C^n}$  the associated \textbf{bump pseudo metric} and $\R^{n-k} \times B^k_{r[n]}(0) \subset \R^{n-k} \times C^k$ its \textbf{core}.
\end{proposition}

\textbf{Proof} \, We deform $\Psi_\circ$ on $C$ consecutively in 3 steps while we gradually approach the axis $\R^{n-k} \times \{0\}$. To simplify the definition of the deformations we work also with large radii, but since the setup is scaling invariant we can finally rescale the steps to compose them seamlessly.\\

\textbf{Step 1} \textbf{($\Psi_\circ  \rightarrowtail   \Psi_\bullet$)}\, We first pass from $\Psi_\circ$ to $\Psi_\circ+\Psi_\bullet$.
We choose a cut-off $\phi \in C^\infty(\R,[0,1])$ with $\phi \equiv 1$ on $\R^{\le 0}$, $\phi \equiv 0$  on $\R^{\ge 1}$ and $\phi' \le 0$. For large $D >0$, we set $h_D(z):=\phi_D(\rho):=\phi(\rho-D)$, where $\rho=\dist(z,\R^{n-k} \times \{0\})$ and get (usually omitting to write the coordinate $z$)
\begin{equation}\label{ers}
(\Psi_\circ(z) +h_D(z)\cdot \Psi_{\bullet}(z))^{(n+2)/(n-2)} \cdot \scal\left( (\Psi_\circ(z) + h_D(z)\cdot \Psi_{\bullet}(z))^{4/(n-2)} \cdot g_{C^n}\right)=
\end{equation}
\[\textstyle -(\Delta\Psi_\circ +h_D \cdot \Delta\Psi_{\bullet}) + \frac{(n-2)}{4(n-1)} \cdot \scal(g_{C^n}) \cdot (\Psi_\circ + h_D\cdot \Psi_{\bullet}) - (\Delta h_D  \cdot \Psi_{\bullet} + 2 \cdot \langle \nabla h_D, \nabla \Psi_{\bullet} \rangle) \]
\vspace{-0.4cm}
\begin{equation}\label{dd1}
 = \lambda \cdot \bp^2 \cdot (\Psi_\circ + h_D\cdot \Psi_{\bullet})  - (\Delta h_D  \cdot \Psi_{\bullet} + 2 \cdot \langle \nabla h_D, \nabla \Psi_{\bullet} \rangle)
\end{equation}
\vspace{-0.8cm}
\begin{equation}\label{uab1}
\ge \psi_C(\omega) \cdot \Big(\lambda \cdot   \rho^{-2} \cdot (\rho^{\alpha_\circ} + \phi_D\cdot \rho^{\alpha_\bullet}) -  (\phi''_D + (k-1)/\rho \cdot \phi'_D) \cdot \rho^{\alpha_\bullet} - 2 \cdot  \alpha_\bullet  \cdot \phi'_D \cdot \rho^{\alpha_\bullet-1}  \Big),
\end{equation}
 from $\bp(z) \ge 1/\dist(z,\Sigma)$ and $|\phi'_D|, |\phi''_D| \le c(\phi)$, for some $c(\phi) >0$, independent of $D$. For large $D(n,k,\alpha_\circ,\alpha_\bullet)>0$, independent of the chosen cone the expression (\ref{uab1}) is \emph{positive}, since $\lambda>0$, $\alpha_\circ-\alpha_\bullet \ge  2 \cdot \sqrt{2/3}  \cdot \frac{k-2}{2}\ge 3$ and, hence, $\rho^{\alpha_\circ-2}/\rho^{\alpha_\bullet} \ra \infty$, for $\rho \ra \infty$.\\

Now we pass from $\Psi_\circ+\Psi_\bullet$ to $\Psi_{\bullet}$.
For $E >0$,  $h_E(z):=\phi_E(\rho):=1-\phi(\rho/E -1)$ we have
\[(h_E(z)\cdot \Psi_\circ(z) + \Psi_{\bullet}(z))^{(n+2)/(n-2)} \cdot \scal\left( (h_E(z)\cdot \Psi_\circ(z) + \Psi_{\bullet}(z))^{4/(n-2)} \cdot g_{C^n}\right) \ge \]
\vspace{-0.9cm}
\begin{equation}\label{uab21}
\psi_C(\omega) \cdot \left(\lambda \cdot  \rho^{-2} \cdot (\phi_E\cdot \rho^{\alpha_\circ} + \rho^{\alpha_\bullet}) -  (\phi''_E + (k-1)/\rho \cdot \phi'_E) \cdot \rho^{\alpha_\circ} - 2 \cdot  \alpha_\circ  \cdot \phi'_E \cdot \rho^{\alpha_\circ-1}  \right).
\end{equation}
For small $E(n,k,\alpha_\circ,\alpha_\bullet)>0$ this is \emph{positive}. Namely, $|\phi'_E| \cdot E, |\phi''_E| \cdot E^2\le c^*(\phi)$, for some $c^*(\phi) >0$
and  $\alpha_\circ-\alpha_\bullet >0$ and, hence, for $\rho \le E$: $\rho^{\alpha_\bullet-2}/(\rho^{\alpha_\circ} \cdot E^{-2}) \ge \rho^{\alpha_\bullet-2}/\rho^{\alpha_\circ-2} \ra \infty$, for $\rho \ra 0$.
Summarizing we have some $D>10$, an $E \in (0,1/10)$ and a smooth function $\Psi>0$ so that $\scal(\Psi^{4/(n-2)} \cdot g_{C^n}) >0$ on $C^n \setminus \Sigma,$ with $\Psi \equiv \Psi_\circ$ for $\rho \ge D+1$ and $\Psi \equiv \Psi_\bullet$ for $\rho \le E.$ We observe that $D$ and $E$ can be chosen depending continuously on $\alpha_\circ$, $\alpha_\bullet$ and, hence, on $C$.\\

\textbf{Step 2} \textbf{($\Psi_\bullet  \rightarrowtail  c \cdot \Psi_\circ$)}\, Now we pass from $\Psi_\bullet$ back to $c \cdot \Psi_\circ$, for some $c>0$. Since
\begin{itemize}
  \item $\psi_C(\omega) \cdot \rho^{\alpha_\bullet}$ and   $c \cdot \psi_C(\omega) \cdot \rho^{\alpha_\circ}$ solve
$L_{C^k,n} \phi = \lambda \cdot \bp^2 \cdot \phi$
  \item $(E/8)^{\alpha_\bullet-\alpha_\circ} \cdot \rho^{\alpha_\circ}|_{\rho=E} = \rho^{\alpha_\bullet}|_{\rho=E}$ and $(E/8)^{\alpha_\bullet-\alpha_\circ} \cdot  (\rho^{\alpha_\circ})'|_{\rho=E/8} >   (\rho^{\alpha_\bullet})'|_{\rho=E/8}$,
\end{itemize}
we can find a smooth interpolation $F>0$, continuously depending on $\alpha_\circ$, $\alpha_\bullet$, such that  $F\le \min \{(E/8)^{\alpha_\bullet-\alpha_\circ} \cdot \rho^{\alpha_\circ},\rho^{\alpha_\bullet}\}$,  $F(\rho)=E^{\alpha_\bullet-\alpha_\circ} \cdot \rho^{\alpha_\circ}$ for $\rho <99/100 \cdot (E/8)$, $F(\rho)= \rho^{\alpha_\bullet}$ for $\rho>101/100 \cdot (E/8)$, with
 \[L_{C^k,n}(\psi_C(\omega) \cdot F(r)) = -\psi_C(\omega) \cdot \left(F''  + (k-1)/\rho \cdot F'\right) + L_{\lambda}^\times \psi(\omega) \cdot F(r)= \]
 \vspace{-0.9cm}
 \begin{equation}\label{sui}
  \psi_C(\omega) \cdot \left(- (F''  + (k-1)/\rho \cdot F') + \mu_{C^k,L_{\lambda}^\times} \cdot F\right) \ge \lambda/2 \cdot \bp^2 \cdot \psi_C(\omega) \cdot F(r).
 \end{equation}

\textbf{Step 3} \textbf{($c \cdot \Psi_\circ  \rightarrowtail  \Psi_\circ$)}\, In Step 2 we returned to $E^{\alpha_\bullet-\alpha_\circ} \cdot \Psi_\circ$. To reach $\Psi_\circ$, we now show how to pass from $c \cdot \Psi_\circ$, for $c :=E^{\alpha_\bullet-\alpha_\circ}>1$, to $\Psi_\circ$ while keeping $\scal>0$. For $d \in \Z$, we  set  $h^*_d(z):=\phi^*_d(\rho):= \phi(2^d \cdot (\rho-2^{-d}) )$. For $\zeta \in (-1/2,0)$, we have:
\[((1+\zeta \cdot h^*_d)\cdot\Psi_\circ)^{(n+2)/(n-2)} \cdot \scal\left(((1+\zeta \cdot h^*_d)\cdot\Psi_\circ)^{4/(n-2)} \cdot g_{C^n}\right)\]
\vspace{-0.9cm}
\begin{equation}\label{hd3}
=   \lambda  \cdot \bp^2   \cdot (1+\zeta \cdot h^*_d) \cdot \Psi_\circ  -  \zeta \cdot \big(\Delta h^*_d \cdot \Psi_\circ + 2 \cdot \langle \nabla h^*_d, \nabla \Psi_\circ \rangle\big)
\end{equation}
Now we write $\Psi_\circ(\omega,\rho)=\psi(\omega) \cdot \rho^{\alpha_\circ}$ and recall that $\bp(\omega,\rho)= \bp(\omega,1) \cdot \rho^{-1}$ and, since $ |\delta_{\bp}(p)- \delta_{\bp}(q)|   \le  d_H(p,q)$, we have $\bp(\omega,1) \ge 1$. Thus (\ref{hd3}) is lower bounded by
\begin{equation}\label{let}\,
\psi(\omega) \cdot \Big(\lambda \cdot \rho^{-2} \cdot (1+\zeta \cdot \phi^*_d)  -  |\zeta| \cdot |(\phi^*_d)'' + (k-1)/\rho \cdot (\phi^*_d)'| - 2 \cdot |\zeta \cdot \alpha_\circ  \cdot (\phi^*_d)'| \cdot \rho^{-1}  \Big) \cdot \rho^{\alpha_\circ}.
\end{equation}
From (\ref{abb}) we have negative upper and lower bounds on $\alpha_\circ$, independent of $C$, and, hence, when $|\zeta| \le \zeta_0(d)$, for some small $\zeta_0(d)>0$,    (\ref{hd3}) is \emph{positive}.  $\zeta_0(d)$ can actually be chosen independently of $d$, i.e.  $\zeta_0=\zeta_0(n,k)$.
Namely, since $\Psi_\circ$ is unique$^\cs$, scaling by $2^d$ transforms
\begin{equation}\label{resca}
(1+\zeta \cdot h^*_d)^{4/(n-2)}\cdot\Psi_\circ^{4/(n-2)} \cdot g_{C^n} \mm{ into } (1+\zeta \cdot h^*_0)^{4/(n-2)}\cdot\Psi_\circ^{4/(n-2)} \cdot g_{C^n},
\end{equation}
up to multiplication by a positive constant. (More explicitly, the scaling effects $(2^{-d} \cdot \rho)^{-2}=2^{2 \cdot d} \cdot \rho^{-2}$, $(\phi^*_d)'(2^{-d} \cdot z) \cdot (2^{-d} \cdot \rho)^{-1}=2^{2 \cdot d} \cdot (\phi^*_1)'(z)\cdot \rho^{-1}$ and $(\phi^*_d)''(2^{-d} \cdot z)=2^{2 \cdot d} \cdot (\phi^*_1)''(z)$ are all the same.)
The interiors of $\su |\nabla h^*_{d_1}|$ and $\su |\nabla h^*_{d_2}|$ are disjoint when $d_1 \neq d_2$ and, hence, we can iteratively multiply by $1+\zeta \cdot h^*_d$ keeping $\scal>0$ when $|\zeta| \le \zeta_0(n,k)$:
\begin{equation}\label{sal15}
\textstyle \scal((\prod^j_{d=1}(1+\zeta \cdot h^*_d)\cdot\Psi_\circ)^{4/(n-2)} \cdot g_{C^n}) >0, \mm{ for any } j \ge 1.
\end{equation}
In turn, we observe that for any fixed $\zeta \in (-1/2,0)$ and $\eta>0$ there is some $j >0$ so that $\textstyle \prod^j_{d=1}(1+\zeta \cdot h^*_d) \le \eta$ on $\R^{\le 2^{-j}}$. \\
Therefore, there are some $\zeta < 0$ and $j \ge 1$  so that (\ref{sal15}) holds, $\prod^j_{d=1}(1+\zeta \cdot h^*_d) = 1$, for $\rho \ge 1$, and  $= 1/c$, for $\rho \le 2^{-j}$.  Since $c =E^{\alpha_\bullet-\alpha_\circ}>1$ is upper bounded, there is a large $j (i,n,k)$,  independent of $C$, so that we can choose $\zeta < 0$ depending continuously on $c$ and thus on $C$.\\

\textbf{Conclusion}\, We scale the constructions in steps 1-3 individually. In step 1 and 2 we rescale $C^n$ by  $(D+1)^{-1}$ and multiply $\Psi_\circ$ by  $(D+1)^{-\alpha_\circ}$. Then the deformation starts from $\mathcal{T}_1$. To append step 3, we rescale  $C^n$ by $E/10$ and
multiply $\Psi_\circ$ by  $(E/10)^{-\alpha_\circ}$. This defines $\F(C) \cdot \Psi_\circ(C)$ as (i). The various radii and $\Lambda$ depend on $n$ and $k$, but since there are only finitely many $k \le n$ we can drop $k$ choosing common estimates.
From the compactness of $\mathcal{SC}_{n}$  we find some lower bound $r(n)>0$ so that $\su |\mathbf{f}(C)-1| \subset [r,1]$ and some upper bound $f \le \Lambda^{(n-2)/2 }$ for any $C \in \mathcal{SC}_{n}$ as claimed in (ii).  We observe that each step still works, with adjusted parameters and function  $\F$, when we replace the dominating term  $\lambda  \cdot \bp^2  \cdot \phi$ for $\lambda/2  \cdot \bp^2  \cdot \phi$ in (\ref{dd1}), (\ref{uab21}) and (\ref{hd3}). Then we finally add $\lambda/4 \cdot \bp^2 \cdot \F \cdot \Psi_\circ$ and get $L_{C^n} (\F \cdot \Psi_\circ) \ge \lambda/4 \cdot \bp^2 \cdot \F \cdot \Psi_\circ$ and (iii).
(iv)   follows from Step 1 and Lemma \ref{etl}.
 \qeda

\subsubsection{Removal of Intersections with $\Sigma_C$} \label{vcd}

\begin{figure}[htbp]
\centering
\includegraphics[width=0.5\textwidth]{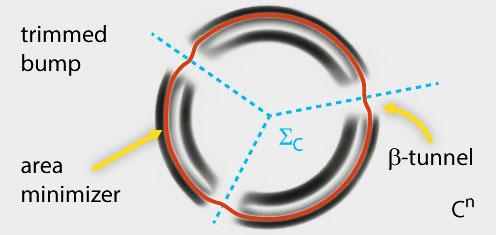}
\caption{\footnotesize Illustration of a $\boldsymbol{\beta}$\textbf{-trimming} near $\Sigma_C$, here for $n-k=0$, to get bumps supported in $C \setminus \Sigma_C$. In Ch.~\ref{rol} we show that for small $\beta>0$ the resulting $\boldsymbol{\beta}$\textbf{-tunnels} are still \textbf{tight} enough to shield the core. \normalsize }
\label{fig:tunnel}
\end{figure}

The \emph{metric distance} $\dist(\cdot,\Sigma)$ has allowed us to find the area minimizer $\mathcal{T}_\varrho$ in Prop.\:\ref{el} and we have $\su g^{n,n-k}_\Delta \cap \R^{n-k} \times \{0\} \v$. In general, $C^k$ is singular also outside $\{0\}$ and $\su g^{n,n-k}_\Delta \cap \R^{n-k} \times (\Sigma_{C^k} \cap \overline{B_{1}(0)} \setminus B_r(0)) \n$. An important case is that of $\Delta^{n,m}$, for $0 \le  m < n-k$, on $\R^{n-k}  \times C^k$ considered as $\R^{m} \times  (\R^{n-k-m}  \times C^k)$, cf.~Fig.\:\ref{fig:tunnel}.\\
We use the \emph{\si-distance} $\delta_{\bp}=1/\bp$ to trim such bump pseudo metrics so that their support no longer intersects $\Sigma_C$ while keeping $\scal>0$. This is a preparation to be able to combine and transfer local bumps to general hypersurfaces.  Since $\delta_{\bp}(z)$ is merely Lipschitz regular, we also use its  Whitney smoothing $\delta_{\bp^*}$.  $\delta_{\bp^*}$ has the same coarse geometric properties as $\delta_{\bp}$, cf.~\cite[Appendix B]{L1}. In \emph{normal coordinates} $\p/ \p x_i$, $i=1,...,n$, in $z \in C \setminus \Sigma_C$, we have:
\begin{equation}\label{smot}
1/c_* \cdot \delta_{\bp}(z) \le \delta_{\bp^*}(z)  \le c_* \cdot \delta_{\bp}(z)\quad\mm{and}\quad |\p^\beta \delta_{\bp^*}  / \p x^\beta |(z) \le c_{**}(\beta) \cdot \delta_{\bp}^{1-|\beta|}(z)
\end{equation}
for constants $c_*(n)\ge 1$, $c_{**}(n,\beta)>0$. $\beta$ is a multi-index for derivatives. In particular,
\begin{equation}\label{nael}
|\nabla \delta_{\bp^*}|(z) \le c_\nabla \mm{ and } |\Delta \delta_{\bp^*}|(z) \le c_\Delta \cdot \delta_{\bp}^{-1}(z) = c_\Delta \cdot \bp(z), \mm{ for  } c_\nabla(n), c_\Delta(n)\ge 1.
\end{equation}
Now we introduce the secondary or $\beta$-trimming deformation.
\begin{proposition} \emph{\textbf{($\boldsymbol{\beta}$-Trimming)}}\label{elt}   There is a $\beta_0(n) \in (0,1)$ so that for $\beta \in [0,\beta_0)$ there is a smooth $\F_\beta(x,y):= \mathbf{f}_\beta(y)$ depending continuously on $C$ in $C^{2,\alpha}$-norm, with $\F\ge \F_\beta \ge 1$ and $\F_{\beta} \in {\cal{C}}[C^n, \Lambda[n]]$ and some $\delta(\beta,n) > \beta$, with $\delta(\beta) \ra 0$ for $\beta \ra 0$,  so that:
\begin{enumerate}
\item  $g^{n, n-k \,|\,\beta}_{C^n}= (\F_\beta \cdot \Psi_\circ)^{4/(n-2)} \cdot g_{C^n}$ equals $g^{[n, n-k]}_{C^n} \mm{ for }  \delta_{\bp^*}(z) \ge \delta$,  $g^\circ_{C^n} \mm{ for } \delta_{\bp^*}(z) \le \beta$, i.e.  the $\boldsymbol{\beta}$\textbf{-trimmed bump} $ \Delta^{n, n-k \,|\,\beta}:=\big((\F_\beta-1) \cdot \Psi_\circ(x,y)\big)^{4/(n-2)} \cdot g_{C^n}$  is supported on a subset of $\delta^{-1}_{\bp^*}(\R^{\ge \beta})$ and  $g^{n, n-k \,|\,0}_{C^n}=g^{n, n-k}_{C^n}$,
     \item  $L_{C^n} (\F_\beta \cdot \Psi_\circ) \ge \lambda/8 \cdot \bp^2 \cdot \F_\beta \cdot \Psi_\circ$, in particular, $g^{n, n-k \,|\,\beta}_{C^n}$ has $\scal>0$.
\end{enumerate}
We call \emph{$\delta_{\bp^*}^{-1}(\R^{\le \beta}) \cap \su g^{n,n- k}_{C^n}$} a $\boldsymbol{\beta}$\textbf{-tunnel}.
\end{proposition}

\textbf{Proof of \ref{elt}} \, The idea of the $\beta$-trimming process is the same as in  Step 3 of \ref{el} where we have deformed $c \cdot \Psi_\circ$, for some $c >1$, to
$\Psi_\circ$ near  $\R^{n-k} \times \{0\}$. That process can be interpreted as a trimming of the support of $(c-1) \cdot \Psi_\circ$. Here we trim $\Delta^{n, n-k}$ towards the remainder of $\Sigma_C$, i.e., we gradually deform $\F-1$ to zero keeping $\scal>0$. We exploit again that we are using the unique$^\cs$ solution $\Psi_\circ>0$ of the equation $-\Delta \Psi_\circ - \frac{n-2}{4 (n-1)} \cdot |A|^2 \cdot \Psi_\circ = \lambda  \cdot  \bp^2 \cdot \Psi_\circ$.
\vspace{-0.1cm}
For  a cut-off function $\phi \in C^\infty(\R,[0,1])$ with $\phi \equiv 1$ on $\R^{\le 0}$, $\phi \equiv 0$  on $\R^{\ge 1}$, we define
 $h^\diamondsuit_d(z):=\phi^\diamondsuit_d(\delta_{\bp^*}(z)):= \phi(2^d \cdot (\delta_{\bp^*}(z)-2^{-d}) )$, for $d \in \Z$ and $z \in C \setminus \Sigma_C$. We collect some basic properties $h^\diamondsuit_d$ where we denote the interior of a set $A$ by $A^\circ $:
\begin{itemize}[leftmargin=*]
\item $\su |\nabla h^\diamondsuit_{d}| \subset \delta_{\bp^*}^{-1}([2^{-d},2^{1-d}])$ and $\su |\nabla h^\diamondsuit_{d_1}|^\circ \cap \su |\nabla h^\diamondsuit_{d_2}|^\circ \v$, for $d_1 \neq d_2$,
  \item  from (\ref{nael}) and $|\phi'|, |\phi''| \le c_\phi$, for some $c_\phi \ge 1$, we get  for $z \in C \setminus \Sigma_C$:
\begin{equation}\label{diaest}
|\nabla h^\diamondsuit_d|(z) \le n \cdot  2^d \cdot  c_\nabla \cdot c_\phi \, , \,  |\Delta h^\diamondsuit_d|(z) \le n \cdot c_\phi \cdot  \big( 2^{2 \cdot d} \cdot  c^2 _\nabla + 2^d \cdot   c_\Delta  \cdot \bp(z)\big).
\end{equation}
\end{itemize}
For $\mathbf{P}_i^j(\zeta):=\prod^j_{d=i}(1+\zeta \cdot h^\diamondsuit_d)$, where $\zeta<0$ with norm $|\zeta| < 1/2$, we have
\begin{equation}\label{die}
\textstyle|\nabla \mathbf{P}_i^j(\zeta)|(z) \le |\zeta| \cdot n\cdot  2^a \cdot  c_\nabla \cdot c_\phi, \mm{ on } \su |\nabla h^\diamondsuit_a|,
\end{equation}
\vspace{-0.8cm}
\begin{equation}\label{die2}
\textstyle|\Delta\mathbf{P}_i^j(\zeta))|(z) \le |\zeta| \cdot n \cdot c_\phi \cdot  \big( 2^{2 \cdot a} \cdot  c^2 _\nabla + 2^a \cdot   c_\Delta  \cdot \bp(z)\big), \mm{ on } \su |\nabla h^\diamondsuit_a|,
\end{equation}
and for any integers $i  \le  a$ and any $\eta >0$ there is an $j(i,\eta,\zeta) > a$ with
\begin{equation}\label{n}
\textstyle 0 <\mathbf{P}_i^j(\zeta)(z)   \le \eta, \mm{ for }  z \in C \setminus \Sigma_C  \mm{ with }  \delta_{\bp^*}(z) \le 2^{-j}.
\end{equation}

\textbf{Step 1 (Cut-Off adapted to $\boldsymbol{\delta_{\bp^*}}$)} \,  For $\zeta \in (-1/2,0)$ and integers $ i  <  j$,  we define
\begin{equation}\label{a}
\textstyle g^{n, n-k}_{C^n}(i,j):= \big((\mathbf{P}_i^j(\zeta)\cdot (\F-1) +1) \cdot \Psi_\circ\big)^{4/(n-2)} \cdot g_{C^n}.
\end{equation}
We note that $g^{n, n-k}_{C^n}(i,j)=(\F \cdot \Psi_\circ)^{4/(n-2)}(z) \cdot g_{C^n}$ for $\delta_{\bp^*}(z) \ge 2^{1-i}$ and that $g^{n, n-k}_{C^n}(i,j)$ asymptotically approaches $\Psi_\circ^{4/(n-2)} \cdot g_{C^n}$ towards $\Sigma$, from (\ref{n}). We claim that there is an  $\zeta_0(n,k) \in (0,1/2)$ so that, for any $\zeta<0$ with norm $|\zeta| \le \zeta_0(n,k)$ and any pair $i < j$, we have  $\scal(g^{n, n-k}_{C^n}(i,j))>0$ on $C \setminus \Sigma_C$. As before, we start from the $\scal$-transformation law:
\[\textstyle  \big((\mathbf{P}_i^j(\zeta)\cdot (\F-1) +1) \cdot \Psi_\circ\big)^{\frac{n+2}{n-2}} \cdot \scal(g^{n, n-k}_{C^n}(i,j))=\]
\[\textstyle -\Delta \big((\mathbf{P}_i^j(\zeta)\cdot (\F-1) +1) \cdot \Psi_\circ\big) - \frac{(n-2)}{4(n-1)} \cdot |A|^2 \cdot \big((\mathbf{P}_i^j(\zeta)\cdot (\F-1) +1) \cdot \Psi_\circ\big)\]
Using the brackets $[$, $\llbracket$ and $[[$ to group the terms the second line of this equation becomes
\[\textstyle [-\Delta \mathbf{P}_i^j(\zeta)\cdot \F \cdot \Psi_\circ]  + [-2 \cdot \langle\nabla(\mathbf{P}_i^j(\zeta),\nabla (\F\cdot \Psi_\circ)\rangle] + [[- \mathbf{P}_i^j(\zeta)\cdot \Delta (\F \cdot \Psi_\circ)]] +...\]
\[[-\Delta (1 -\mathbf{P}_i^j(\zeta) \cdot \Psi_\circ] + [-2 \cdot  \langle\nabla(1 - \mathbf{P}_i^j(\zeta) , \nabla \Psi_\circ \rangle] +\llbracket-(1 - \mathbf{P}_i^j(\zeta)) \cdot \Delta \Psi_\circ\rrbracket + ...\]
\[\textstyle [[- \frac{(n-2)}{4(n-1)} \cdot |A|^2 \cdot \mathbf{P}_i^j(\zeta)\cdot \F \cdot \Psi_\circ]] +\llbracket - \frac{(n-2)}{4(n-1)} \cdot |A|^2 \cdot \big(1-\mathbf{P}_i^j(\zeta)\big)\cdot \Psi_\circ\rrbracket\]

\begin{itemize}[leftmargin=*]
\item For the sum of the $[$-terms we use $\nabla (\F\cdot \Psi_\circ) = \F \cdot \nabla  \Psi_\circ + \Psi_\circ \cdot \nabla \F$ and (\ref{diaest})-(\ref{die2}).
The definition of $\F$ in Prop.\:\ref{el} shows that $\F, |\nabla \F|, |\Delta \F| \le c_\F$ on $C \setminus \Sigma_C$, for some $c_\F(n,k)>0$. This gives us the following lower estimate on $\delta_{\bp^*}^{-1}([2^{-a},2^{1-a}])$:
\begin{equation}\label{hd6a}
  - |\zeta| \cdot n \cdot c_\phi \cdot  \big( 2^{2 \cdot a} \cdot  c^2 _\nabla + 2^a \cdot   c_\Delta  \cdot \bp(z)\big)\cdot c_\F \cdot \Psi_\circ - 2  \cdot n\cdot |\zeta| \cdot 2^a \cdot  c_\nabla \cdot c_\phi \cdot c_\F \cdot \Psi_\circ -...
\end{equation}
\[2  \cdot n\cdot |\zeta| \cdot 2^a \cdot  c_\nabla \cdot c_\phi \cdot c_\F \cdot  |\nabla \Psi_\circ | - |\zeta| \cdot n \cdot c_\phi \cdot  \big( 2^{2 \cdot a} \cdot  c^2 _\nabla + 2^a \cdot   c_\Delta  \cdot \bp(z)\big)\cdot \Psi_\circ -2  \cdot n\cdot |\zeta| \cdot 2^a \cdot  c_\nabla \cdot c_\phi \cdot  |\nabla \Psi_\circ | .\]
\item The sum of the $[[$ and $\llbracket$-terms  can be lower estimated from \ref{el}(iii) and $\F \ge 1$:
\begin{equation}\label{bb}
\lambda \cdot \bp^2 \cdot \big(1-\mathbf{P}_i^j(\zeta)\big)\cdot \Psi_\circ + \lambda/4 \cdot \bp^2 \cdot \mathbf{P}_i^j(\zeta) \cdot  \F \cdot \Psi_\circ \ge \lambda/4 \cdot \bp^2 \cdot   \Psi_\circ.
\end{equation}
\end{itemize}

\textbf{Step 2 (Elliptic Estimates relative to $\boldsymbol{\delta_{\bp^*}}$)} \, To compare (\ref{hd6a}) with (\ref{bb}), we estimate  $|\nabla \Psi_\circ |$ in terms of $|\Psi_\circ |$.  $\su |\nabla h^\diamondsuit_a| \subset \delta_{\bp^*}^{-1}([2^{-a},2^{1-a}])$  is geometrically well-controlled: for $a=1$ we have $\delta_{\bp^*}^{-1}([1/2,1]) \subset \delta_{\bp}^{-1}([c^{-1}_*/2, c_*])$ and thus
\begin{equation}\label{bgeo}
|A|(z) \le 2 \cdot c_* , \mm{ for any } z \in \delta_{\bp^*}^{-1}([1/2,1]) \mm{ from } \bp \ge |A|.
\end{equation}
From Gau\ss\ equations this yields a uniform control over the geometry
of $\delta_{\bp^*}^{-1}([1/2,1])$ even on the neighborhood $V:=\bigcup \{ B_{1/(2 \cdot \bp(p)}(p)\,|\, p \in \delta_{\bp}^{-1}([c_a^{-1}/2,c_a])\}$ from  \cite[Lemma B.2]{L1}:
\begin{equation}\label{aum}
|\bp(p)/\bp(q)-1| \le \bp(p) \cdot d_H(p,q),\mm{ for }q \in B_{1/(2 \cdot \bp(p)}(p).
\end{equation}
We get
for any $\epsilon>0$ an $r(n,\epsilon)>0$ so that for $x \in \delta_{\bp^*}^{-1}([1/2,1])$ the exponential map $\exp_x:B_{r(n,\epsilon)}(0) \ra C \setminus \Sigma_C$ is a local diffeomorphism with $|\exp^*_x(g_C)-g_{T_xC}|_{C^5(B_{r(n,\eta)}(0))}\le \epsilon$. $g_{T_xC}$ is the flat metric on the tangent space $T_xC$ in $x$. \\
The coefficients of $-\Delta \Psi_\circ - \frac{n-2}{4 (n-1)} \cdot |A|^2 \cdot \Psi_\circ = \lambda  \cdot  \bp^2 \cdot \Psi_\circ$ are uniformly bounded in $C^\alpha$-norm on any $B_{r(n,\eta)}(x)$, $\alpha \in (0,1)$, hence we get uniform constants for the Harnack inequality and elliptic estimates for  $\Psi_\circ$ on $B_{r(n,\eta)}(x)$ independent of $C$, i.e., there is some $c_\Box(n,k)>0$ with
\begin{equation}\label{nao}
|\nabla \Psi_\circ|   \le c_\Box \cdot \Psi_\circ \mm{ on }  V  \mm{ and, thus, } |\nabla \Psi_\circ|   \le c_\Box \cdot 2^{a} \cdot \Psi_\circ \mm{  on } 2^{-a} \cdot V,
\end{equation}
where we have also used that $\Psi_\circ>0$ is the unique$^\cs$ solution on $C \setminus \Sigma_C$. Since $\delta_{\bp_{c \cdot H}}=c \cdot \delta_{\bp}$, for any $c>0$, and thus  $c \cdot c^{-1}_* \cdot \delta_{\bp^*} \le \delta_{\bp^*_{c \cdot H}} \le c \cdot c_* \cdot \delta_{\bp^*}$,  we have
\begin{equation}\label{suu}
\delta_{\bp^*}^{-1}([2^{-a},2^{1-a}]) \subset \delta_{\bp}^{-1}([c^{-1}_* \cdot 2^{-a}, c_* \cdot 2^{1-a}])= 2^{-a} \cdot V.
\end{equation}
From  (\ref{nao})  we find some $\zeta_0(1) \in (0,1/2)$ so that for $\zeta<0$ with  $|\zeta| \le \zeta_0$,  the norm of  (\ref{hd6a}) is smaller that $\lambda/16 \cdot \bp^2 \cdot   \Psi_\circ$ on $V$. As in \ref{el}, we observe that, due to the uniqueness of $\Psi_\circ$, we have, under scaling by $2^a$, that, up to a constant multiple
\begin{equation}\label{say}
\Psi_\circ \mm{ and } h^\diamondsuit_a \mm{ on } 2^{-a} \cdot V \mm{  transform into } \Psi_\circ \mm{ and } h^\diamondsuit_1 \mm{ on }  V.
\end{equation}
From (\ref{suu}), (\ref{nao}) and $\bp_{2^{-a} \cdot H}= 2^{a} \cdot \bp_{H}$ we therefore see that, for $a >0$ and $|\zeta| \le \zeta_0$, the norm of  (\ref{hd6a}) is smaller than $\lambda/100 \cdot \bp^2 \cdot   \Psi_\circ$
 on $2^{-a} \cdot V$, i.e. $\zeta_0$ can be chosen to depend on $n$ and $k$ and not on $a$. This shows that for $\zeta<0$ with $|\zeta| \le \zeta_0$ and $0 < i \le j$ on $C \setminus \Sigma_C $:
\begin{equation}\label{ijj}
\textstyle L_{C^n} (\mathbf{P}_i^j(\zeta)\cdot (\F-1) +1) \cdot \Psi_\circ \ge \lambda/6 \cdot (\mathbf{P}_i^j(\zeta)\cdot (\F-1) +1) \cdot \Psi_\circ.
\end{equation}

\textbf{Step 3 (Iterated Cut-Offs)} \, We use this $\zeta_0(n,k)$ and some $\zeta<0$ with small enough $|\zeta| \le \zeta_0$ and, from (\ref{n}), some sufficiently large $j$  so that  the trimmed pseudo metric
\begin{equation}\label{diam}
\textstyle \big(((1-h^\diamondsuit_j) \cdot \mathbf{P}_i^j(\zeta)\cdot (\F-1) +1) \cdot \Psi_\circ\big)^{4/(n-2)} \cdot g_{C^n}
\end{equation}
has $\scal>0$ on $C \setminus \Sigma_C$. Namely, we consider the transformation under scaling by $2^j$:
\begin{equation}\label{sayy}
\textstyle \mathbf{P}_i^j(\zeta) \mm{ on } 2^{-j} \cdot V  \mm{  transform into }  \mathbf{P}_{i-(j-1)}^1(\zeta)(1+\zeta \cdot h^\diamondsuit_d) \mm{ on }  V
\end{equation}
and this uniformly $C^3$-converges to zero on $V$, when $j \ra \infty$. From this a suitably large $j(i,n,k)$ can be chosen  so that the metric (\ref{diam}) has $\scal>0$  on $C \setminus \Sigma_C$.\\
To choose the parameter $\delta(\beta)>0$, we start with $i=1$ and $j(1)$ and set $\beta_0(n,k) :=2^{-j(1)}$. For any $\beta \in (0,\beta_0(n,k))$, we therefore have some $\delta(\beta) \le 1$ as claimed. For $2^{-j(i+1)} \le \beta < 2^{-j(i)}$, we choose $\delta(\beta):=2^{-i}$ and we see that
$\delta(\beta) \ra 0$ for $\beta \ra 0$. In turn, we associate these $i$ and $j(i+1)$ to any $\beta$ with $2^{-j(i+1)} \le \beta < 2^{-j(i)}$ and we set
\begin{equation}\label{fff}
\textstyle \F_\beta:=((1-h^\diamondsuit_j) \cdot \mathbf{P}_i^j(\zeta) \cdot (\F-1) +1).
\end{equation}
As before, we can drop the dependency on $k$ to simplify our notations. From the discussion above we see that for $\F_\beta$ we have $L_{C^n} (\F_\beta \cdot \Psi_\circ) \ge \lambda/8 \cdot \bp^2 \cdot \F_\beta \cdot \Psi_\circ$ and it satisfies the other asserted properties directly from the construction.
\qeda

\subsubsection{Tightness of Cavities} \label{rol}

We show that there are area minimizers $\mathcal{T}^\beta$ relative to the trimmed bump metric $g^{n, n-k \,|\,\beta}_{C^n}= (\F_{\beta} \cdot \Psi_\circ)^{4/(n-2)} \cdot g_{C^n}$ and, as a part of the argument, we also prove that $\mathcal{T}^\beta$ Hausdorff converges to $\mathcal{T}_\varrho$, for $\beta \ra 0$, where $\mathcal{T}_\varrho$ is the minimizer relative to  $g^{[n, n-k]}_{C^n}$ of Prop.\:\ref{el}.\\
In our assembly of global bumps, we superpose disjointly supported families of trimmed bumps to $(C,d_{\sima})$ and later to  $(H,d_{\sima})$, for $H \in {\cal{G}}^c$, cf.~Fig.\:\ref{fig:localbump}. Therefore we consider a more general $\Lambda$-bumped situation and include the following type of additional deformations:
\begin{equation}\label{gve}
(\G_{\beta} \cdot \F_{\beta} \cdot \Psi_\circ)^{4/(n-2)} \cdot g_{C^n} \mm{ with } \su |\G_{\beta} -1| \cap \su |\F_{\beta} -1|  \v \mm{ and } \G_\beta \in {\cal{C}}[C^n, \Lambda[n]].
\end{equation}
The $\G_\beta$ are the placeholders for further local bumps. They do \emph{not} need to be translation $\R^{n-k}$-invariant and we note from $ \su |\G_{\beta} -1| \cap \su |\F_{\beta} -1|  \v$ that $\G_{\beta} \cdot \F_{\beta} \in {\cal{C}}[C^n, \Lambda[n]]$.
We denote the associated hypersurface and volume measures by $\mu^{\Lambda, n-1}_{\sima}[\beta,\G_{\beta}]$ and $\mu^{\Lambda}_{\sima}[\beta,\G_{\beta}]$, where $[\beta,\G_{\beta}]= [0,1]$ means the case of the untrimmed metric $(\F \cdot \Psi_\circ)^{4/(n-2)} \cdot g_{C^n}$.\\

\textbf{Local versus Global Minimizers} \, Now we consider local area minimizers in the metric measure space $(C,d^\Lambda_{\sima},\mu^\Lambda_{\sima})$, of  Theorem \ref{lbp},  associated to $(\G_{\beta} \cdot \F_{\beta} \cdot \Psi_\circ)^{4/(n-2)} \cdot g_{C^n}$. The attribute \emph{local} means that these minimizers are area minimizing under, not quantitatively specified, small perturbations. The hypersurfaces $T_\varrho$ in Prop.\:\ref{el} or the minimal boundaries in our main Theorems belong to this class. In some cases we can specify an open test set $\Omega \subset (H,d^\Lambda_{\sima},\mu^\Lambda_{\sima})$ so that the area minimizer becomes a \emph{global} area minimizer when compared to all other hypersurfaces within $\Omega$ and with the same (Plateau) boundary data along $\p \Omega$. In the BV-approach, $\p \Omega$ also serves as an \emph{obstacle} that keeps interior points of these minimizers to stay in $\overline{\Omega}$, cf.~\cite[Rm.\:1.22]{Gi}. We choose the auxiliary
\begin{equation}\label{ome}
\Omega :=  Q_\ell^{n-k} \times C^k, Q_\ell^{n-k}:=(-\ell/2,\ell/2)^{n-k}, \ell \ge 1 \mm{ and } O_\rho:=\R^{n-k} \times B_\rho(0), \rho>0,
\end{equation}
to get existence results and, in a second step, we check that the minimizers actually avoid $\p \Omega$ away from their Plateau boundary.
Concretely, for small $\beta$  we get the existence of global minimizers in a suitably chosen $\Omega \subset (C,d^\Lambda_{\sima},\mu^\Lambda_{\sima})$ close to the unique $T_\varrho$ we have for $\beta=0$. For $\beta >0$, and even in $\Omega$, such global minimizers need not to be unique, but  we observe a \textbf{stability property}: all competing global area minimizers, relative to $\Omega$, stay close to $T_\varrho$.

\begin{proposition} \emph{\textbf{(Tightness of  $\beta$-Tunnels)}}\label{lue2}
There is some $\beta^*[n,\ell] \in (0,1)$, so that for $\beta \in [0,\beta^*]$, $\F_\beta$ from Prop.\:\ref{elt}, and any $\G_\beta \in {\cal{C}}[C^n, \Lambda[n]]$ with \emph{$\su |\G_{\beta} -1| \cap \su |\F_{\beta} -1|  \v$}, we have for $\vartheta_i[n,C] \in (r,1)$, $\varrho \in (\vartheta_1,\vartheta_2)$ from Prop.\:\ref{el} and  $\zeta=1/4 \cdot \min \{|\varrho - \vartheta_1|,|\varrho - \vartheta_2|\}$:
\begin{enumerate}
\item there is an open Caccioppoli set $U[n-k,\ell,\G_\beta]\subset C^n$ with $ \overline{O_{\varrho-\zeta}} \subset U \subset O_{\varrho+\zeta}$ so that $\mathcal{T}[n-k,\ell,\G_\beta] :=\p U$ and $\mathcal{T}  \cap \Omega$ is area minimizing with $\Omega$-boundary value $\mathcal{T}_{\varrho}$, relative to $(C,d^\Lambda_{\sima},\mu^\Lambda_{\sima})$ associated to $(\G_{\beta} \cdot \F_{\beta} \cdot \Psi_\circ)^{4/(n-2)} \cdot g_{C^n}$,
\item for any such oriented minimal boundary $\mathcal{T}_\bullet :=\p U_\bullet \subset \overline{O_{\varrho+ 2 \cdot \zeta}} \setminus O_{\varrho- 2 \cdot \zeta}$ with the obstacle $\p (O_{\varrho+  2 \cdot \zeta} \setminus O_{\varrho-2 \cdot \zeta})$, we already have $\mathcal{T}_\bullet\subset O_{\varrho+\zeta} \setminus \overline{O_{\varrho-\zeta}}$. We write  $\mathcal{P}[n-k,\ell,\G_\beta]$ for the class of all such Plateau solutions $\mathcal{T}_\bullet$.
\end{enumerate}
In the case $k=n$, $\ell$ is a void parameter, but we keep writing it for a consistent notation. In turn, for known dimension of the Euclidean factor we also drop writing $n-k$.
\end{proposition}

\textbf{Proof} \, We consider the auxiliary \emph{obstacle problem} for open Caccioppoli sets $U$
\begin{equation}\label{sus}
\R^{n-k} \times B_{\vartheta_1}(0) \subset  \R^{n-k} \times B_{\varrho-2 \cdot \zeta}(0) \subset U \subset \R^{n-k} \times B_{\varrho+2 \cdot \zeta}(0) \subset \R^{n-k} \times B_{\vartheta_2}(0).
\end{equation}
To be able to use convergence arguments, we ensure compactness properties for the obstacles. For $n=k$ any such $\p U \subset \overline{B_{\vartheta_2}(0)} \setminus B_{\vartheta_1}(0)$ is compact.
For $n>k$ we choose the compact set $\overline{Q_\ell^{n-k}} \times \overline{B_{\vartheta_2}(0)} \setminus B_{\vartheta_1}(0) \subset \Omega$. We let $\mathcal{T} \subset \Omega$ be an oriented minimal boundary $\mathcal{T}=\p U$ with $\Omega$-boundary value $\mathcal{T}_{\varrho}$ solving the obstacle problem (\ref{sus}). We observe that $Q_\ell^{n-k}$ is an intersection of Euclidean halfspaces and thus $\p \Omega$ is locally outer minimizing, that is, its area increases under local outward deformations, cf. the argument of \ref{oma}. That is, a \emph{free} Plateau solution with $\Omega$-boundary value $\mathcal{T}_{\varrho}$ does not leave $\overline{\Omega}$. We claim that $\mathcal{T}$ is also a free solution in the radial direction if $\beta>0$ is small enough and  we even show that $\mathcal{T}[\ell,\G_\beta]\subset O_{\varrho+\zeta} \setminus \overline{O_{\varrho-\zeta}}$.\\

\textbf{Step 1}\, We define the $\Omega$\emph{-flat norm} as the $\mu^\Lambda_{\sima}[0,1]$-volume of the difference set $(U \, \Delta \, O_\varrho) \cap  \Omega$. In this step we show that for the oriented minimal boundaries $\mathcal{T}=\p U$,
\begin{equation}\label{fl}
\mathcal{T}[\ell,\G_\beta]  \ra  \mathcal{T}_\varrho \mm{ in }\Omega\mm{-flat norm, for }\beta \ra 0,
\end{equation}
uniformly for all singular area minimizing cones $C^n= \R^{n-k} \times C^k$ and for all admissible $\mathbf{G}_\beta$.
From the naturality of $\bp_H$ on $\cal{G}$ and its properness on $H \setminus \Sigma_H$, for any non-totally geodesic $H \in \cal{G}$, we get some constant $\nu(\beta,\Omega)>0$, independent of $\R^{n-k} \times C^k$, so that
\begin{equation}\label{c22}
 \mu^{\Lambda, n-1}_{\sima}[0,1]\big(\delta^{-1}_{\bp^*}(\R^{\le \delta(\beta)}) \cap \mathcal{T}_\varrho \cap \Omega\big) \le \nu \mm{ and } \nu(\beta,\Omega)\ra 0 \mm{ for } \beta \ra 0.
\end{equation}
The minimality of the $\mathcal{T}$ and of $\mathcal{T}_\varrho$ and  $\G_{\beta}\ge 1$ show that
\begin{equation}\label{c23}
\mu^{\Lambda, n-1}_{\sima}[0,1]\big(\mathcal{T}_\varrho \cap \Omega\big) \le \mu^{\Lambda, n-1}_{\sima}[\beta,\G_{\beta}]\big(\mathcal{T} \cap   \Omega\big)
\le \mu^{\Lambda, n-1}_{\sima}[0,1]\big(\delta^{-1}_{\bp^*}(\R^{\ge \delta}) \cap \mathcal{T}_\varrho \cap \Omega\big) + \Lambda^{n-1} \cdot \nu.
\end{equation}
Since $\mu^{\Lambda, n-1}_{\sima}[\beta,\G_{\beta}](\mathcal{T} \cap \Omega) \ge \mu^{\Lambda, n-1}_{\sima}[0,1](\mathcal{T} \cap \Omega)$, we infer that for $\beta \ra 0$,
\begin{equation}\label{c1}
\mu^{\Lambda, n-1}_{\sima}[\beta,\G_{\beta}](\mathcal{T}[\ell,\G_{\beta}]  \cap \Omega) \mm{ and } \mu^{\Lambda, n-1}_{\sima}[0,1](\mathcal{T}[\ell,\G_{\beta}] \cap \Omega) \ra \mu^{\Lambda, n-1}_{\sima}[0,1](\mathcal{T}_\varrho \cap \Omega).
\end{equation}

From this, the BV-compactness relative to  $\mu^{\Lambda}_{\sima}[0,1]$, \cite[Prop.\:3.7]{M}, shows that for any sequence $\beta_i \ra 0$, for $i \ra \infty$, there is a $\Omega$-flat norm converging subsequence of the $\mathcal{T}[\ell,\G_{\beta_i}] \cap \Omega$. From  the lower semicontinuity of the perimeter \cite[Prop.\:3.6]{M} and the uniqueness of $\mathcal{T}_\varrho$ we infer that the limit is $\mathcal{T}_\varrho$. Another loop of the argument shows the convergence of the entire sequence $\mathcal{T}[\ell,\G_{\beta_i}]  \cap \Omega$,
\begin{equation}\label{delta}
\mu^{\Lambda}_{\sima}[0,1]\big((U[\ell,\G_{\beta_i}]  \, \Delta \, O_\varrho) \cap \Omega\big) \ra 0 \mm{ for } i \ra \infty.
\end{equation}
From compactness results for minimal factor geometries on cones and from the continuous dependence of $\F_\beta$ and $\varrho$ on $C$, this convergence is uniform for all $C^n= \R^{n-k} \times C^k$ and  $\G_\beta$.\\

\textbf{Step 2}\, We apply the $\mu^{\Lambda}_{\sima}$-growth estimate of  Corollary \ref{grr} (\ref{est}) to the difference set to upgrade the \emph{flat norm} convergence from step 1 to \emph{Hausdorff} convergence.\\
 If $\overline{U[\ell,\G_{\beta_i}]} \cap \R^{n-k} \times \p B_{\varrho + \zeta}(0) \n$ for all $i$, then there is a $p_i \in (U[\ell,\G_{\beta_i}] \, \Delta \, O_\varrho) \cap \R^{n-k} \times \p B_{\varrho + \zeta/2}(0)$. We have $B_{\zeta/4}(p_i) \cap \R^{n-k} \times \p B_{\varrho}(0) \v$, $p_i \in \overline{U[\ell,\G_{\beta_i}]}$ and the  $\mathcal{T}[\ell,\G_{\beta_i}]$ are (global) area minimizers relative to the test set $B_{\zeta/4}(p_i)\subset\Omega$.  Thus, Cor.\:\ref{grr} (\ref{est}) shows
 \setlength{\jot}{12pt}
\begin{equation}\label{mki}
 \kappa^+ \cdot (\zeta/8)^n  \le  \mu^{\Lambda}_{\sima}[\beta,\G_{\beta}](U[\ell,\G_{\beta_i}] \cap B_{\zeta/8}(p_i)) \le \Lambda[n]^{n-1} \cdot  \mu^{\Lambda}_{\sima}[0,1]\big((U[\ell,\G_{\beta_i}] \, \Delta \, O_\varrho) \cap \Omega\big).
\end{equation}

This contradicts (\ref{delta}). For $\overline{U[\ell,\G_{\beta_i}]}\cap \R^{n-k} \times \p B_{\varrho - \zeta}(0) \v$ we argue similarly. Since only $1 \le \G_{\beta} \le \Lambda$ for $\G_{\beta}$ were used, compactness arguments
for the space of all singular cones $C^n= \R^{n-k} \times C^k$ with their minimal factor geometry show that there is some $\beta^*[n,\ell]>0$, independent of $C^n$ and $\G_{\beta}$, so that for any $\beta \in [0,\beta^*]$:  $\mathcal{T}[\ell,\G_{\beta}]  \subset O_{\varrho+\zeta} \setminus \overline{O_{\varrho-\zeta}}$ proving claim (i). The argument also implies (ii) since it applies to all such area minimizers.\qeda

\subsubsection{Local Bumps and Shields} \label{localbump}
The $\beta$-trimmed bump pseudo metrics $\Delta^{n, n-k \,|\,\beta}$ are supported away from $\Sigma_C$. This is one of the requirements
to transfer such bumps from $C$ to other hypersurfaces. Since  tangent cone approximations generally only work on bounded subsets of the cone, we additionally need to trim the bump pseudo metrics so that their support becomes \emph{compact} in $C \setminus \Sigma_C$.\\
For this we inductively use the $\Delta^{n, q \,|\,\beta}$, for all $q$ starting from $q=0$ up to $q= n-k$, on  $C^n=\R^q  \times  (\R^{n-k-q}  \times C^k)= \R^{n-k} \times C^k$. All of the $\Delta^{n, q \,|\,\beta}$ are defined on the same cone $\R^{n-k} \times C^k$, but for $q < n-k$, the definition of  $\Delta^{n, q \,|\,\beta}$ ignores the $\R^{n-k-q}$-symmetry and treats $\R^{n-k-q}  \times C^k$ as the cone factor. $\Delta^{n, 0 \,|\,\beta}$ already has compact support.

\begin{figure}[htbp]
\centering
\includegraphics[width=0.5\textwidth]{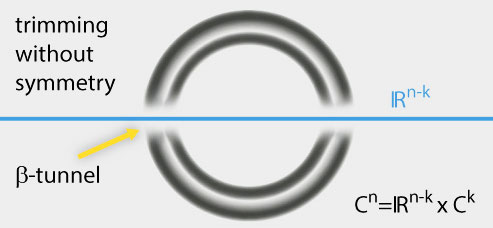}
\caption{\footnotesize The illustration shows $g^{n, 0 \,|\,\beta}_{C^n}$ for $\Sigma_{C^n}=\R^{n-k} \times \{0\}$ where the $\beta$-trimming ignores the $\R^{n-k}$-translation symmetry of $C^n$. By contrast, $g^{n, n-k \,|\,\beta}_{C^n}$ looks like the middle  image of Fig.\:\ref{fig:a} without $\beta$-tunnels.\normalsize }
\label{fig:regtunnel}
\end{figure}
\vspace{-0.5cm}
\begin{definition}\emph{\textbf{(Bump Elements)}}\label{bumel} For $q=0,...n-k$, $\ell \ge 1$
and $\phi \in C^\infty(\R,[0,1])$ with $\phi \equiv 1$ on $\R^{\le 0}$, $\phi \equiv 0$ on $\R^{\ge 1}$, we set $\F^{\phi,Q_\ell^q}_\beta :=\phi \left(\dist\left(z, Q_\ell^q \times  \R^{n-k-q}  \times C^k \right)\right) \cdot (\F_\beta-1)+1$:
\begin{itemize}
\item $\Delta^{n, q  \,|\, \beta, \ell}:=\big(\F^{\phi,Q_\ell^q}_\beta-1\big)^{4/(n-2)} \cdot g_{C^n}$ is the \textbf{(q)-bump element} aligned to $Q_\ell^q$,
\item  $\big(\F^{\phi,Q_\ell^q}_\beta \cdot \Psi_\circ \big)^{4/(n-2)} \cdot g_{C^n}$ is the associated \textbf{bumped metric},
\item  \emph{$\cu \Delta^{n, q  \,|\, \beta, \ell} :=\su \Delta^{n, q  \,|\, \beta, \ell} \cap \{z \in  C^n \,|\, 0 < \dist\left(z, Q_\ell^q ...\right) <1\}$} is the \textbf{cut-off region}.
\end{itemize}
\end{definition}
For our purposes we can think of $\dist$ as a \emph{smooth} function. As in Ch.\:\ref{vcd}, we may choose Whitney smoothings of $\dist$ whenever needed.

\begin{remark}\textbf{($\boldsymbol{\scal}$ on  $\boldsymbol{\cu \Delta}$)}\label{assu}  On $C \setminus (\cu \Delta \cup \Sigma_C)$ the bumped metric  coincides with $g^{n, n-k \,|\,\beta}_{C^n}$ respectively with $\Psi_\circ^{4/(n-2)} \cdot g_{C^n}$. In both case we have $L_{C^n} (\F^{\phi,Q_\ell^q}_\beta \cdot \Psi_\circ) \ge \lambda/8 \cdot \bp^2 \cdot \F^{\phi,Q_\ell^m}_\beta \cdot \Psi_\circ$. In turn, we anticipate that $\cu \Delta$ will be inductively shielded by other bump elements, i.e. it belongs to the finally \emph{deleted} neighborhood $U$ of $\Sigma$ with minimal boundary $\p U$ of $H \setminus U$ we described in our theorems. This makes knowing $\scal$ on $\cu \Delta$ dispensable. Even so, to avoid exceptional subsets in our further discussion we indicate how to ensure that on $\cu \Delta$:
\begin{equation}\label{assum}
L_{C^n} (\F^{\phi,Q_\ell^q}_\beta \cdot \Psi_\circ) \ge \lambda/10 \cdot \bp^2 \cdot \F^{\phi,Q_\ell^m}_\beta \cdot \Psi_\circ \mm{ and, hence, } \scal(\big(\F^{\phi,Q_\ell^q}_\beta \cdot \Psi_\circ \big)^{4/(n-2)} \cdot g_{C^n})>0.
\end{equation}
We observe that  $\F_\beta \cdot \Psi_\circ$ and $\Psi_\circ$ satisfy the relation $L_{C^n} F \ge \lambda/8 \cdot \bp^2 \cdot F$. When we replace $\phi$ by $\phi_\tau$, $\tau >0$, with $\phi_\tau(t):=\phi(\tau \cdot  t)$, we have $|\phi'_\tau|, |\phi''_\tau| \ra 0$, for $\tau \ra 0$. Thus the terms involving $|\phi'_\tau|$ and $|\phi''_\tau| $ are  majorized by $\phi_\tau \cdot L_{C^n} (\F_\beta \cdot \Psi_\circ)+(1-\phi_\tau) \cdot L_{C^n}  \Psi_\circ$ and we get, for small $\tau>0$, that $L_{C^n} (\F^{\phi,Q_\ell^q}_\beta \cdot \Psi_\circ) \ge \lambda/10 \cdot \bp^2 \cdot \F^{\phi,Q_\ell^q}_\beta \cdot \Psi_\circ$. For fixed $\phi$ this $\tau$ depends only on $n$ and $\beta$. To avoid non-essential parameters we henceforth assume that $\tau=1$. \qed
\end{remark}

\textbf{Local Bumps}\, For $0 \le  q \le n-k$, $\cu \Delta^{n, q  \,|\, \beta, \ell}$ is compact. For $q=0$, this leaves $\Delta^{n, 0  \,|\, \beta}$ unchanged, but for the sake of consistency we keep the parameter $\ell$. For (0)-bump elements, Prop.\:\ref{lue2} shows that for small $\beta>0$ there are minimizing boundaries ${\cal{T}}^\beta=\p U$ for an open $U$ with $O_r \subset U \subset O_1$.
For $m>0$ we inductively and disjointly place \textbf{scaled bump elements} along all lower dimensional faces of $Q^{n-k}_1$, so that the cut-off regions of the ($q$)-bump elements belong to the union of the ($\phi$-trimmings of the) \textbf{cores} $\R^q \times B^{n-k-q}_{r[n]}(0) \subset \R^q  \times  (\R^{n-k-q}  \times C^k)$, from \ref{el}, of the already assigned  ($q$)-bump elements, for $q <m$, cf.~Fig.\:\ref{fig:localbump}. For each ($q$)-dimensional face $Q^{q,f}_1 \subset Q^{n-k}_1$, $f \in F_q$ for some index set $F_q$, we  choose an affine transformation
\begin{equation}\label{trtr}
T[\eta,f]:\R^{n-k} \times C^k \ra \R^{n-k} \times C^k \mm{ with } (x,y) \mapsto \eta\cdot  (x, A_f(y)+ v_f),
\end{equation}
where $\eta \in (0,\ell^{-1}], A_f \in O(n-k),  v_f \in \R^{n-k},$ that maps the model ($q$)-cube $Q_\ell^q$ to $Q^{q,f}_1$ so that $T[\ell^{-1},f](Q_\ell^q)=Q^{q,f}_1$. From Prop.\ref{dimshift} we have $\Psi_\circ(\omega,\rho)=\psi(\omega) \cdot \rho^{\alpha_\circ}$  and hence
\begin{equation}\label{trtf}
T[\eta,f]_*(\Psi_0)= \Psi_0 \circ T[\eta,f]^{-1} = \eta^{-\alpha_\circ} \cdot \Psi_0.
\end{equation}
The parameters $\beta$, $\ell$ and $\eta$ are determined inductively in the following construction. The choices anticipate additions of ($q$)-bump elements for $q\ge 1$, as shown in Fig.\:\ref{fig:localbump}, and of local bumps when we form larger configurations finally reaching global bumps. To this end, we recall from \ref{lue2} that  any area minimizer $\mathcal{T}_\bullet \in \mathcal{P}[q,\ell,\G_\beta]$ with $\Omega$-boundary value $\mathcal{T}_{\varrho}$, for $\Omega_q = Q_\ell^q \times \R^{n-k-q} \times C^k$ and $O_\rho=\R^{n-k} \times B_\rho(0), \rho>0$,  satisfies
\begin{equation}\label{tor}
\mathcal{T}_\bullet \cap \Omega \subset (O_{\varrho+\zeta} \setminus \overline{O_{\varrho-\zeta}}) \cap  \Omega \subset (O_1 \setminus \overline{O_r}) \cap  \Omega,
\end{equation}
where $r=r[n] \in (0,1)$ is the core radius of $\big(\F  \cdot \Psi_\circ \big)^{4/(n-2)} \cdot g_{C^n}$ in Prop.\:\ref{el}.

\begin{description}[leftmargin=*]
\item[0-faces] For $\eta_0=1/10$, we place $\eta_0 $-scaled  (0)-bump elements $\Delta^{n, 0  \,|\, \beta, \ell}$ in each 0-dimensional face (=vertex) $v_f:=p_f$ of $Q^{n-k}_1$.
For $\beta^*[n] \in (0,1)$ from \ref{lue2}, $\beta_0 \in (0,\beta^*]$, $\xi \in [r[n],1]$, we set
\[\textstyle  \pmb{\square}^{\beta_{0},\eta_{0}, \ell_0} :=  \sum_{f\in F_0} T[\eta_0 ,f]_*\big(\Delta^{n, 0 \,|\,\beta_0, \ell_0}\big)  \mm{, }\,  \mathbb{O}^{\eta_0,\ell_0}_{0,\xi} :=\bigcup_{f\in F_0}   O^{\eta_0,\ell_0}_{f,0,\xi}:= \bigcup_{f\in F_0} T[\eta_0 ,f]\big(O_\xi\big)\]
for the void parameter $\ell_0=1$. We observe that \[\su T[\eta_0 ,f_1]_*\big(\Delta^{n, 0 \,|\,\beta_0, \ell_0}\big) \cap \su T[\eta_0 ,f_2]_*\big(\Delta^{n, 0 \,|\,\beta_0, \ell_0}\big) \v \mm{ for } f_1 \neq f_2 \in F_0.\]
\item[q-faces] We continue inductively for  $q$ with $1 \le q \le n-k$. We choose a large $\ell_q \ge 1$, some $\eta_q \in (0, \min\{\eta_0 \cdot \beta_0/10,...,\eta_{q-1} \cdot \beta_{q-1}/10,\ell_q^{-1}\} )$ and an $\beta_q \in [0,\beta^*[n,\ell_q]]$, so that after placing the $\eta_q$-scaled ($q$)-bump elements $\Delta^{n, q  \,|\, \beta_q, \ell_q}$ along  the ($q$)-faces of the cube $Q^{n-k}_1$, the support of any two
($q$)-bump elements and of any  ($q$)-bump element and  the already placed  ($m$)-bump element, $0 \le  m <q$ are \emph{disjoint}, cf.~(i) and (ii) below. The cut-off region of each  ($q$)-bump element  belongs to the union of the cores of the  ($m$)-bump elements for $0 \le  m <q$, cf.~(iii).
\begin{enumerate}
\item $\su T[\eta_q,f_1]_*(\Delta^{n, q  \,|\, \beta_q, \ell_q}) \cap \su T[\eta_q,f_2]_*(\Delta^{n, q \,|\, \beta_q, \ell_q}) \v$, for $f_1\neq f_2  \in F_q$,
\item $\su T[\eta_q,f]_*(\Delta^{n, q \,|\, \beta_q, \ell_q}) \cap \bigcup_{k =0,..q-1} \su \bs^{\beta_{k},\eta_{k},\ell_{k}} \v$, for $f \in F_q$,
\item $\bigcup_{f =F_q}  T[\eta_q,f](\cu \Delta^{n, q\,|\, \beta_q, \ell_q})\subset  \bigcup^{q-1}_{m=0} \mathbb{O}^{\eta_m,\ell_m}_{m,r[n]}$,
\end{enumerate}
and for the induction step we set for $q$ the following unions of bump elements and cylinders:
\begin{align}\label{oo}
\bs^{\beta_{q},\eta_{q}, \ell_q} := & \textstyle  \sum_{f =F_q} T[\eta_q, f]_*\big(\Delta^{n, q \,|\, \beta_q, \ell_q}\big)\mm{ and }\\
 \mathbb{O}^{\eta_q, \ell_q}_{q,\xi} := & \textstyle \bigcup_{f\in F_q}   O^{\eta_q,\ell_q}_{f,q,\xi}:=
\bigcup_{f\in F_q} T[\eta_q ,f]\big(O_\xi \cap \Omega_q\big).  \nonumber
\end{align}

\end{description}
$\bigcup^{n-k}_{q=0} \mathbb{O}^{\eta_q,\ell_q}_{q,\xi}$ is an open neighborhood  of $Q^{n-k}_1 \times \{0\},$ for any $\xi \in [r[n],1]$, and from Prop.\:\ref{el} and Prop.\:\ref{lue2}, we notice that the inductively chosen parameters
\begin{equation}\label{indpar}
\ell_q \ge 1, \eta_q>0 \mm{ and } \beta_q>0, \mm{ for } 0 < q \le n-k,
\end{equation}
are independent of $C^n=\R^{n-k} \times C^k \in {\cal{SC}}_n$. We make one \textbf{fixed choice} of these parameters and keep it for the rest of this paper.\\

{\begin{figure}[htbp]
\centering
\includegraphics[width=0.7\textwidth]{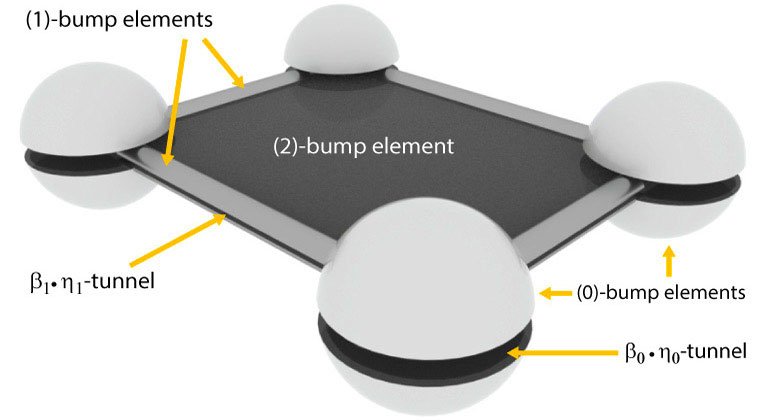}
\caption{\footnotesize An illustration of the \textbf{support} of a \textbf{local bump} shielding $Q^{n-k}_1 \times \{0\}$. \normalsize}
\label{fig:localbump}
\end{figure}}

\begin{definition}\emph{\textbf{(Local Bumps)}}\label{locbu}  We call the pseudo metric
\begin{equation}\label{deflb}
\textstyle \bs(Q^{n-k}_1):= \sum^{n-k}_{q=0} \bs^{\beta_q,\eta_q, \ell_q} \mm{ an  \textbf{(n-k)-local bump} on } C^n=\R^{n-k} \times C^k.
\end{equation}
We define the associated \textbf{(n-k)-local bump metric}
\begin{equation}\label{defbb}
\textstyle (\Phi_\square \cdot \Psi_\circ )^{4/(n-2)} \cdot g_{C^n}:= \left(\big(\prod^{n-k}_{q=0} \prod_{f =F_q} (\F^{\phi,Q_\ell^m}_{\beta_q} \circ T[\eta_q,f]^{-1})\big) \cdot \Psi_\circ \right)^{4/(n-2)} \cdot g_{C^n}.
\end{equation}
and we call $\mathbb{O}^C_{n-k}:=\bigcup^{n-k}_{q=0} \mathbb{O}^{\eta_q,\ell_q}_{q,r[n]}$ the \textbf{core} of the local bump.
\end{definition}
$\Phi_\square$ is the multiplicative representation of a local bump with $\Phi_\square\in {\cal{C}}[C^n, \Lambda[n]]$ and we have $\su |\Phi_\square - 1| = \su \bs(Q^{n-k}_1)$
\begin{equation}\label{jo}
L_{C^n} (\Phi_\square  \cdot \Psi_\circ) \ge \lambda/10 \cdot \bp^2 \cdot \Phi_\square  \cdot \Psi_\circ.
\end{equation}
This follows from \ref{assu} (\ref{assum}) and (\ref{trtf}) since the $\su |\F^{\phi,Q_\ell^m}_{\beta_q} \circ T[\eta_q, f]^{-1}-1|$ for any two different ($q$)-bump elements are disjoint and the relation (\ref{jo}) remains invariant under translations and scalings.  We use Prop.\:\ref{lue2} to understand the effect of local bumps.
\begin{proposition}\emph{\textbf{(Shielding)}}\label{oma}
For $\G \in {\cal{C}}[C^n, \Lambda[n]]$  with \emph{$\su |\G -1| \cap \su |\Phi_\square -1|  \v$} and  the $\ell_q, \eta_q,\beta_q$ from (\ref{indpar}), we choose for any ($q$)-face $f\in F_q$ some
\begin{equation}\label{upen}
\mathcal{T}_f  = \p U_f \in \mathcal{P}[q,\ell_q,\G_{\beta_q}] \mm{ for an open } O_{r[n]} \cap \Omega_q  \subset  U_f \cap \Omega_q \subset O_1 \cap \Omega_q
\end{equation}
where $\G_{\beta_q}:=\G \circ T[\eta_q ,f]$ on $O_1 \cap \Omega$. Combining these $U_f$ we get a \textbf{local shield} for $Q^{n-k}_1$:
\begin{equation}\label{open}
\textstyle \mathbb{U}(Q^{n-k}_1):=\bigcup^{n-k}_{q=0} \mathbb{U}^{\beta_q,\eta_q,\ell_q}_{q}:= \bigcup^{n-k}_{q=0} \bigcup_{f\in F_q} T[\eta_q ,f]\big(U_f \cap \Omega_q\big) \supset \mathbb{O}^C_{n-k} \supset Q^{n-k}_1 \times \{0\},
\end{equation}
\vspace{-0.8cm}
 \begin{equation}\label{ibm}
 \p \mathbb{U} \mm{  is a local \textbf{inner area minimizer} relative to } \mathbb{U} \mm{ and }(\G \cdot \Phi_\square \cdot \Psi_\circ )^{4/(n-2)} \cdot g_{C^n}.
 \end{equation}
That is, $\mu^{\Lambda, n-1}_{\sima}(\p\mathbb{U}) \le \mu^{\Lambda, n-1}_{\sima}(\p\big(\mathbb{U}\setminus A\big))$ for sufficiently small open sets $A \subset C$.
\end{proposition}

This means we take the standard (q)-bump element aligned to $Q_\ell^q$, defined in \ref{bumel},  and use $\G \circ T[\eta_q ,f]$
for the standard (q)-bump element to apply Prop.\ref{lue2}.  Then the resulting $\mathcal{T}_f$ and $U_f$ are shifted by $T[\eta_q ,f]$ to their associated face $f$ of $Q_1^{n-k}$.\\

\textbf{Proof}\, From Prop.\:\ref{lue2} we see that $\mathbb{U}$ is an open neighborhood of $Q^{n-k}_1 \times \{0\}$.
For (\ref{ibm}) we note that the union of two locally inner area minimizing $U$ and $V$ is again locally inner area minimizing:  for an open $A$ we first use this inner minimality of $U$, for $A  \setminus V$, and then that of $V$, for $A  \cap V$, and write $(U \cup V) \setminus A= \big((U \cup V) \setminus (A \setminus V)\big) \setminus A \cap V.$\\
This applies to $\p \mathbb{U}$ since, dropping $\Omega_q$ for $q=0$, $\p \mathbb{U} \cap \bigcup^{n-k}_{q=0} \bigcup_{f\in F_q} T[\eta_q ,f]\big(U_f \cap \Omega_q\big) \v$.
We only need to consider  $A \cap \bigcup^{n-k}_{q=0} \bigcup_{f\in F_q} T[\eta_q ,f]\big(U_f \cap \Omega_q\big) \v$ and we argue by induction over all $T[\eta_q ,f]\big(U_f \cap \Omega_q\big)$ in the union $\mathbb{U}$.
\qed

\textbf{Groupings}\, We also need bumps shielding arbitrarily large cubes $Q^{n-k}=Q^{n-k}_{a_1,...,a_{n-k}}=[-a_1,a_1] \times \cdots\times[-a_{n-k},a_{n-k}] \times \{0\} \subset \R^{n-k} \times C^k$, for some $a_i>0$. To control the interplay with other bumps we must keep the size of both, the (support of) local bumps in the $C^k$-directions and of the tunnel parameters $\beta_q$, unchanged. This means, we cannot just scale the bump construction. Alternatively, when we try to adjust the parameters in (\ref{deflb}) to shield a larger cube $Q$,  we find from \ref{lue2} that, starting from an initial $\beta_0$-tunnel for the $(0)$-faces and $(0)$-bump elements, the subsequent $\beta_i$ have to be chosen the smaller the larger $Q$ becomes.\\
The way out are \emph{groupings} of local bumps, each shielding a translated copy of the unit cube $Q^{n-k}_1$.
Concretely,  we consider disjoint unions $\bigcup_{v \in S} (Q^{n-k}_1+ v) \subset Q^{n-k}_{L,...,L} \subset  \R^{n-k}$, for  finite, otherwise arbitrary subsets $S \subset \Z^{n-k}$ and some $L=L(S)>0$ large enough. \\
We repeat the inductive process we have used to define local bumps and place ($m$)-bump elements, $0 \le  m \le n-k$, along the ($m$)-faces of the cubes in this union. We count all faces in this union with \textbf{multiplicity one} even when they belong to more than one cube $Q^{n-k}_1+ v$. Applying \ref{locbu} and \ref{oma} to this periodic version of the local bump construction yields:

\begin{corollary}\label{perg}\emph{\textbf{(Periodic Groupings)}} We get \textbf{periodic bumps} and \textbf{shields} from the obvious extensions of the definitions (\ref{deflb}) and (\ref{open}):
\begin{equation}
 \textstyle \bs(\bigcup_{v \in S} (Q^{n-k}_1+ v)), \mm{ the bumped metric }   (\Phi_\square \cdot \Psi_\circ )^{4/(n-2)} \cdot g_{C^n} \mm{ and } \mathbb{U}(\bigcup_{v \in S} (Q^{n-k}_1+ v)),
 \end{equation}
 and we define the \textbf{core} of the periodic bump by
 \begin{equation}\label{coo}
 \textstyle \mathbb{O}^C_{n-k}(\bigcup_{v \in S} (Q^{n-k}_1+ v)):=\bigcup_{v \in S} \big(\bigcup^{n-k}_{q=0} \mathbb{O}^{\eta_q,\ell_q}_{q,r[n]} + v\big).
\end{equation}
For $\G \in {\cal{C}}[C^n, \Lambda[n]]$  with \emph{$\su |\G -1| \cap \su |\Phi_\square -1|  \v$}, and the $r$-distance neighborhood $U_r$ in $C^n$ of the family $S$ of cubes, we have:
 \begin{enumerate}
\item $\mathbb{U}(\bigcup_{v \in S} (Q^{n-k}_1+ v)) \supset \mathbb{O}^C_{n-k}(\bigcup_{v \in S} (Q^{n-k}_1+ v)) \supset \bigcup_{v \in S} (Q^{n-k}_1+ v)$,
\item  $\p \mathbb{U}$ is locally inner area minimizing relative to $\mathbb{U}$, i.e. $\mathbb{U}$ shields $\mathbb{O}^C_{n-k}(\bigcup_{v \in S} (Q^{n-k}_1+ v))$,
\item   $L_{C^n} (\Phi_\square  \cdot \Psi_\circ) \ge \lambda/10 \cdot \bp^2 \cdot \Phi_\square  \cdot \Psi_\circ$.
 \end{enumerate}
These results extend to $\bigcup_{v \in S} (Q^{n-k}_1+ v) \subset Q^{n-k}_{L,...,L} \subset  \R^{n-k}$ after scaling by $2^a$, $a \in \Z^{\ge 1}$, where we interpret $2^a  \cdot Q^{n-k}_1$ as a union of $2^{a \cdot (n-k)}$ appropriately translated copies of $Q^{n-k}_1$. Rescaling the resulting bump by $2^{-a}$ yields the $2^a$\textbf{-rescaled periodic bumps} $\bs[a]$ and \textbf{shields} with
\begin{equation}\label{perg2a}
\textstyle \mathbb{U}[a](\bigcup_{v \in S} (Q^{n-k}_1+ v)) \supset  \mathbb{O}^C_{n-k}[a](\bigcup_{v \in S} (Q^{n-k}_1+ v)) \supset \bigcup_{v \in S} (Q^{n-k}_1+ v),
\end{equation}
with $\mathbb{U}[a]:=2^{-a} \cdot \mathbb{U}(2^a  \cdot \bigcup_{v \in S} (Q^{n-k}_1+ v))$ and $\mathbb{O}^C_{n-k}[a]:=2^{-a} \cdot \mathbb{O}^C_{n-k}(2^a \cdot \bigcup_{v \in S} (Q^{n-k}_1+ v))$.
\end{corollary}

\textbf{Proof}\,  We note that, unless to faces of different cubes coincide, all bump elements have again disjoint support. From this we can repeat the same bump element-wise argument as in the proof of Prop.\:\ref{oma} for the same parameters with the same conclusions. \qeda

\begin{remark}\label{perg2}\textbf{(Estimates for $\su \bs$)}\, We set $r^\circ[n]:=  \min\{\beta_0 \cdot \eta_0,...,\beta_{n-7} \cdot \eta_{n-7}\} \cdot r[n]/4$
as a lower bound for all radii, independent of the dimension of the factor $\R^{n-k}$, and note that
\begin{equation}\label{uaa}
\textstyle U_{2^{-a}}(\bigcup_{v \in S} (Q^{n-k}_1+ v)) \setminus U_{2^{-a} \cdot r^\circ[n]}(\R^{n-k} \times \{0\}) \supset \su \bs[a](\bigcup_{v \in S} (Q^{n-k}_1+ v)).
\end{equation}
Estimates in terms of $\delta_{\bp}$ are compatible with transfers between $H$ and tangent cones.
From (\ref{elt})(i) and $\delta_{\bp}(z) \le dist(z,\Sigma_C)$, we get for any $p \in \Sigma$
\begin{align}\label{bus}
\textstyle \overline{B_{2^{-a} \cdot r^\circ[n]}(p)} \cap  \su \bs[a](\bigcup_{v \in S} (Q^{n-k}_1+ v)) & \subset ... \\
\textstyle  \{ z \in C \setminus \Sigma_C \,|\, \delta_{\bp}(z) \le 2 ^{-a}\cdot r^\circ[n] \}  \cap  \su \bs[a](\bigcup_{v \in S} (Q^{n-k}_1+ v)) & \v. \nonumber
\end{align}
and from (\ref{uaa}) complementary to (\ref{bus}):
\begin{equation}\label{bus2}
\textstyle   \{ z \in C \setminus \Sigma_C \,|\, \delta_{\bp}(z) \ge 2 ^{-a} \} \cap \su \bs[a](\bigcup_{v \in S} (Q^{n-k}_1+ v)) \v.
\end{equation}  \qeda
\end{remark}
\begin{example}\label{perg3}\textbf{(Shielding Balls)} The rescaling of periodic bumps shrinks the distance of the support to $\bigcup_{v \in S} (Q^{n-k}_1+ v)$ but it also
improves the control over the shielded subset. For $B^{n-k}_1(0) \subset \R^{n-k} \times \{0\}$ we interpret the rescaling as a partition of the cubes into $2^{a \cdot (n-k)}$-subcubes of side length $2^{-a}$, $a \in \Z^{\ge 1}$, and consider the \emph{minimal set} $S_{a,n-k}\subset 2^{-a} \cdot \Z^{n-k}$ of such cubes with  $B^{n-k}_1(0) \subset \bigcup_{v \in S_{a,n-k}} (Q^{n-k}_{2^{-a}}+ v)$. We have from $\diam (Q^{n-k}_1) = \sqrt{n-k}$ :
\begin{equation}\label{rtr}
\textstyle \overline{B^{n-k}_{1+{2^{-a}} \cdot \sqrt{n-k}}(0)} \supset \bigcup_{v \in S_{a,n-k}} (Q^{n-k}_{2^{-a}}+ v) \supset \overline{B^{n-k}_1(0)}.
\end{equation}
From this and (\ref{uaa}) we have for some $a_n \in \Z^{\ge 1}$ large enough and any $m \in \Z^{\ge 1}$
\begin{itemize}[leftmargin=*]
  \item  $B^n_{2}(0) \supset B^n_{1+{2 \cdot 2^{-a_n}} \cdot \sqrt{n-k}}(0) \supset \su \bs[\kappa_n \cdot (a_n + m)](\bigcup_{v \in  S_{a_n,n-k}} (Q^{n-k}_{2^{-a}}+ v))$,
  \item $\su \bs[\kappa_n \cdot a_n](\bigcup_{v \in S_{a_n,n-k}} (Q^{n-k}_{2^{-a}}+ v))  \cap \su \bs[\kappa_n \cdot (a_n + m)](\bigcup_{v \in  S_{a_n,n-k}} (Q^{n-k}_{2^{-a}}+ v)) \v,$
\end{itemize}
where $\kappa_n:=[\,\mm{smallest integer} \ge  |\ln (r^\circ[n])|\,]$, for the radius $r^\circ[n]>0$  in \ref{perg2}. \qeda
\end{example}

\setcounter{section}{3}
\renewcommand{\thesubsection}{\thesection}
\subsection{From Local to Global Deformations} \label{gl}

We turn to a general hypersurface $H \in {\cal{G}}_n$ and gradually scale $H$ around some $p \in \Sigma_H$, by large $\tau\gg1$. Eventually there is some not necessarily unique area minimizing cone $C$ that Hausdorff approximates $\tau \cdot H$ on the unit ball. On compact parts of $C \setminus \Sigma_C$  this can be improved to a smooth norm convergence, cf.~\cite[Appendix A]{L2}. This smooth approximation comes with a uniquely determined smooth map, which we call an $\D$-map, locally parameterizing $\tau \cdot H$ as a section of the normal bundle over $C$, we can use to transfer local bump configurations on $C$ to $H \setminus \Sigma_H$.
The challenge to build global bumps along $\Sigma_H$, from local transfers, is to compensate the differences between the singular sets of $H$ and such cones.
\begin{enumerate}[label=(\arabic*),leftmargin=*]
\item \textbf{Basic Idea}\, We would like to have ball covers of $\Sigma_H$ with (i) common ball radii and controlled intersection numbers $c(n)$, depending only on $n$ and so that (ii) each of the balls admits a good cone approximation to place cone bumps on $H \setminus \Sigma_H$ and to estimate their effect.
Yet, we cannot accomplish both goals (i) and (ii) at the same time. For any $p \in \Sigma_H$ and $\ve>0$ there is a critical radius, the $\ve$\emph{-accuracy radius} $\mathcal{R}_\ve(p)$. Balls of radius $r \in (0,\mathcal{R}_\ve(p))$, suitably rescaled,  admit fine approximations by cones of a (yet to define) accuracy $\ve$. One hurdle is that $\mathcal{R}_\ve(p)$ does not depend continuously on $p \in \Sigma$.
\item \textbf{Self-Similar Covers}\,  To resolve this discontinuity issue, we define a hierarchy on $\Sigma$, the \emph{accuracy decomposition} $\Sigma = \dot{\bigcup}_{i\ge 0}\Sigma_i$, with
$\Sigma_0 = \{ p \in \Sigma \,|\, \mathcal{R}_\ve(p) > 1/\ve \}$ and $\Sigma_i = \{ p \in \Sigma \,|\,  s^{i-1}/\ve \ge \mathcal{R}_\ve(p) >  s^{i}/\ve\}$, for $i \ge 1$ and some $s \in (0,1/100)$, in Ch.\:\ref{ptf}.   For each $\Sigma_i $  we can find covers by balls of radius $s^i$ with the two conditions (i) and (ii) in (1). Scaling by $s$ shifts the family index by $1$ and makes the cover self-similar under scaling by $s$ in the sense that properties we derive for one $\Sigma_i$ carry over to the other families. For small $s>0$, the bumps placed on ball families for different $i$ do not interfere.
\item \textbf{Auto-Alignment}\, For the placement of bumps on families of equally sized balls we observe that
singular and highly curved directions in any two neighboring balls align like compass needles in a common field. This controls not only the relative position of the singular set in these balls but also the positioning of the balls themselves.
\end{enumerate}

\subsubsection{$\D$-Maps and Ball Covers} \label{ptf}

Around any $q \in \Sigma_H$ in some $H \in {\cal{G}}_n$, we have tangent cones of $H$:  $H_i= \tau_i \cdot H$, for some sequence $\tau_i \ra \infty$, $i \ra \infty$, flat norm subconverges to some (generally non-unique) area minimizing tangent cone $C^n \subset \R^{n+1}$. Allard theory shows that this convergence can be upgraded.\\

$\D$\textbf{-maps}\, For $B_R(q) \cap C \setminus \Sigma_C$, $R>0$, large $i$ and $B_R(q_i) \cap H_i$ for suitable $q_i \in H \setminus \Sigma$, there is a local $C^5$-section $\Gamma_i :B_R(q) \cap C \ra  B_R(q_i) \cap H_i\subset\nu$ of  the normal bundle $\nu \mm{ of } B_R(q)\cap C$
 up to minor adjustments near $\p B_R(q)$. For $i \ra \infty$, $\Gamma_i$ converges in $C^k$-norm to the zero section, which we identify with $B_R(q) \cap C$. We call  $\D := \Gamma_i$ the asymptotic identification map, briefly
the $\D$\textbf{-map}. Since they are assigned canonically, we generally omit writing indices. We also recall the definition of the \textbf{\si-pencil $\P(p,\omega)$} pointing to $p\in \Sigma$, for some $\omega>0$: $\P(p,\omega):= \{x \in H \setminus \Sigma_H \,|\, \delta_{\bp_H} (x) > \omega \cdot d_H(x,p)$.
We consider the $\D$-maps on bounded subsets of $\P(p,\omega)$, the \textbf{truncated \si-pencils} $\TP$.
\begin{equation}\label{trpen}
\TP_H(p,\omega,R,r):=B_{R}(p) \setminus B_{r} (p) \cap  \P(p,\omega)\subset H.
\end{equation}
For increasingly large $\tau >0$, $\tau \cdot  \TP_H(p,\omega,R/\tau,r/\tau)$ is gradually better $C^5$-approximated by the truncated \si-pencil in \emph{some} tangent cone, i.e., the twisting of $\P(p,\omega)$ slows down as $\tau\to\infty$, it \emph{asymptotically freezes}. More precisely, we have \cite[Prop.\:3.10]{L4}:
\begin{proposition}[\si-Freezing]\label{freezfu}
Let $H\in {\cal{G}}_n$ and $p \in \Sigma_H$. We pick an $\ve > 0$, some $R > 1 > r >0$ and an $\omega \in (0,1)$. Then there is a \textbf{smallest} $\tau^* = \tau^*(\ve,\omega, R , r, p) \ge 1$ so that for \textbf{any} $\tau> \tau^*$ there is some tangent cone $C^\tau_p$ so that  for the $\D$-map to $\tau \cdot H$:
\begin{itemize}
  \item $|\D - id_{C^\tau_p}|_{C^5(\TP_{C^\tau_p}(0,\omega,R,r))} \le \ve $, where we view the zero section of $\nu$ as the identity $id$,
  \item $\Psi_H$ induces the uniquely$^\cs$ determined solution $\Psi_\circ>0(C^\tau_p)$ on $C^\tau_p$ from Prop.\:\ref{dimshift} and we have, normalizing $\Psi_\circ$ appropriately,  $|\Psi_H \circ \D / \Psi_\circ-1|_{C^{2,\alpha}(\TP_{C^\tau_p}(0,\omega,R,r) )} \le \ve$.
\end{itemize}
We call $\ve$ the \textbf{accuracy}  of the $\D$-map and $\mathcal{R}_\ve(H,p):=  1 / \tau^*(\ve, \ve,1/\ve, \ve, p)$ the \textbf{$\boldsymbol\ve$-accuracy radius} of $p \in \Sigma_{H}$. That is, any ball $B_\rho(p) \subset H$, $0 < \rho \le \mathcal{R}_\ve$, rescaled by  $(\ve \cdot \rho)^{-1}$ admits an $\ve$-accurate $\D$-map over $\TP_C(0,\ve,1/\ve, \ve) \subset B_{1/\ve}(0) \setminus B_{\ve} (0) \cap C$,
  for some tangent cone $C$.
\end{proposition}

\begin{remark} For Euclidean area minimizers, the flat norm $d^\flat$-convergence implies Hausdorff $d_{\cal{H}}$-convergence from volume growth estimates \cite[Prop.\:5.14]{Gi}. We already used this argument in Step 2 of \ref{lue2} for area minimizers in $(C,d^\Lambda_{\sima},\mu^\Lambda_{\sima})$ based on  Corollary \ref{grr} (\ref{est}).\\ The spaces $\mathcal{C}_{n}$ and ${\cal{H}}^{\R}_n$ are compact  under compact convergence in the flat norm topology. $\mathcal{K}_{n}$ is compact in flat norm topology.  $\mathcal{SC}_n \subset\mathcal{C}_n$ and ${\cal{SH}}^{\R}_n \subset {\cal{H}}^{\R}_n$ are closed  and, hence, compact. Consequently $\overline {\cal T}_H \subset \mathcal{SC}_n$, also we have $\overline {\cal T}_p =  {\cal T}_p$.\qeda
\end{remark}

We also recall a variant for singular cones, \cite[Cor.\:3.11]{L4},  where we use the Hausdorff distance on ${\cal{SC}}_n$ measured in terms of the distance of the unit spheres $S_C=\p B_1(0) \cap C$.
\begin{corollary}[Cone Transitions]\label{freezfu2}
Let $\alpha\in(0,1)$, $\ve$, $\omega >0$ and $R > 1 > r >0$. Then there is a $\zeta(\ve,\omega , R , r)>0$ so that for any $C$ and $C' \in {\cal{SC}}_n$ with $d_{\cal{H}}(S_C, S_{C'}) \le \zeta$:
\begin{itemize}
  \item $|\D - id_{C}|_{C^5(\TP_{C}(0,\omega,R,r))} \le \ve$,
  \item $\Psi_\circ(C')$ on $C'$ induces $\Psi_\circ(C)$ on $C$ with $|\Psi_\circ(C) /\Psi_\circ(C') \circ \D -1|_{C^{2,\alpha}(\TP_{C}(0,\omega,R,r))}  \le \ve.$
\end{itemize}
\end{corollary}

\textbf{Self-Similar Ball Covers}\,  We write $\Sigma$ as a disjoint union of subsets $\Sigma_i$, $i \in \Z^{\ge 1}$, sorted according to their $\ve$-accuracy radii. Then we cover these subsets by balls with radii $s^i$, for some small $s>0$. This makes the scenario self-similar under scaling by powers of $1/s$.

\begin{definition} \emph{\textbf{(Accuracy Decomposition)}}\label{add}\, For $H \in {\cal{G}}_n$, some $\ve>0$ and $s \in (0,s_0)$ we define the $\boldsymbol\ve$-\textbf{accuracy decomposition} of $\Sigma = \dot{\bigcup}_{i\ge 0}\Sigma_i:$
 \begin{equation}\label{acde}
 \Sigma_0 = \{ p \in \Sigma \,|\, \mathcal{R}_\ve(p) > 1/\ve \} \mm{ and } \Sigma_i = \{ p \in \Sigma \,|\,  s^{i-1}/\ve \ge \mathcal{R}_\ve(p) >  s^{i}/\ve\}.
 \end{equation}
\end{definition}
 In general, the $\Sigma_i$ will be uncountable, but each of them contains a \textbf{countable} and \textbf{dense}
subset $\Sigma^{*}_i \subset \Sigma_i$ and we set $\Sigma^{*}:=\bigcup_{i \ge 0}\Sigma^{*}_i \subset \Sigma$. We write the elements of $\Sigma^*_i$ accordingly as
$a_{i,1}, a_{i,2},...$.  We assign a ball $\overline{B_{s^{i}}(p)}$ to any of the points $p \in \Sigma^*_i$ and set
\begin{equation}\label{acde2}
\textstyle {\cal{B}}[i]:= \{\overline{B_{s^{i}}(p)}\,|\, p \in \Sigma^*_i\}\mm{ and }{\cal{B}}:= \bigcup_{i\ge 0}{\cal{B}}[i].
\end{equation}
We note that ${\cal{B}}[i]$ is a cover of $\Sigma_i$ and, hence,  ${\cal{B}}$ is a cover of $\Sigma$ since for any $x  \in \Sigma_i$  there is some $p_x  \in \Sigma^*_i$ with $d_H(x,p_x) < s^{i}/2$ and thus $x \in \overline{B_{s^{i}}(p_x)}$.

\begin{proposition} \emph{\textbf{(Self-Similar Covers)}}\label{brop}\, For $H \in {\cal{G}}_n$, with $\Sigma_H \n$, some accuracy  $\ve \in (0,10^{-5})$, $R \in [10^2,\ve^{-1}/10^2]$ and any
self-similarity factor $s \in (0,s_0)$, the following holds:
\begin{enumerate}
\item For each $i \ge 0$, there is a locally finite family ${\cal{A}}[i]\subset {\cal{B}}[i]$, from (\ref{acde2}), of  closed balls
\[{\cal{A}}[i] = \{\overline{B_{s^i}(p)} \subset H\,|\, p \in A[i]\}, \mm{ for some discrete set } A[i] \subset \Sigma^*_i.\]
\item For any $k \ge 0$ we have: $\bigcup_{i \le k} {\cal{A}}[i]$ covers $\bigcup_{i \le k} \Sigma_i$.
 \item  $q \notin \overline{B_{s^i}(p)}$, $p \notin \overline{B_{s^k}(q)}$ for any two different $\overline{B_{s^i}(p)} \in {\cal{A}}[i]$, $\overline{B_{s^k}(q)} \in {\cal{A}}[k]$.
\end{enumerate}
For $i \ge 0$, there is a subset $A^\bullet[i] \subset  A[i]$
and a family ${\cal{A}}^\bullet[i] = \{\overline{B_{R  \cdot s^i}(p)} \subset H\,|\, p \in A^\bullet[i]\}$ with
\begin{equation}\label{22}
\textstyle \bigcup_{x \in A[i]} \overline{B_{s^i}(x)} \subset \bigcup_{z \in A^\bullet[i]} \overline{B_{R  \cdot s^i}(z)} \mm{ and } z^* \notin \overline{B_{(R-1)  \cdot s^i}(z)}, \mm{ for } z \neq z^*, z,z^* \in A[i]
\end{equation}
and so that, for some constant $c(H)$,  ${\cal{A}}^\bullet[i]$ splits into disjoint families ${\cal{A}}^\bullet[i,1],...,{\cal{A}}^\bullet[i,c]$ with
 \begin{equation}\label{10}
 B_{10 \cdot R \cdot s^i}(p)
     \cap B_{10 \cdot R \cdot s^i}(q) \v, \mm{ for $\overline{B_{R \cdot s^i}(p)}, \overline{B_{R \cdot s^i}(q)}$ in the same family }{\cal{A}}^\bullet[i,j].
 \end{equation}
We call the  ${\cal{A}}^\bullet[i]$ the \textbf{layers} of the ball cover ${\cal{A}}^\bullet:=\bigcup_{i \ge 1}{\cal{A}}^\bullet[i]$ of $\Sigma$ around the points of $A^\bullet:=\bigcup_{i \ge 1}A^\bullet[i]=\bigcup_{i \ge 1, c \ge j \ge 1}A^\bullet[i,j]$.  The ${\cal{A}}^\bullet[i,j]$ are the \textbf{sublayers} of the ${\cal{A}}^\bullet[i]$.
\end{proposition}

\begin{remark}\label{bror} 1.~To explain the meaning of \ref{brop} we anticipate that in the assembly of global bumps in Ch.\ref{ifp}  we transfer a periodic bump to $\overline{B_{R  \cdot s^i}(p)} \in {\cal{A}}^\bullet[i]$ that shields all balls $\overline{B_{s^i}(x)}\subset \overline{B_{R  \cdot s^i}(p)}$ with $\overline{B_{s^i}(x)} \in {\cal{A}}[i]$.  The balls in ${\cal{A}}^\bullet[i]$ may intersect but we can ensure that the support of the bumps in the at most $c(n,R)$ different sublayers ${\cal{A}}^\bullet[k,i]$ are disjoint. For balls in ${\cal{A}}^\bullet[i]$, ${\cal{A}}^\bullet[j]$, with $i \neq j$, the bumps are disjoint when $s>0$ is sufficiently small.\\
2.~The covering number $c(H)$ depends only on $n$ if $\ve \cdot R \le r_{H,\Phi}$, for  $r_{H,\Phi}>0$ in \cite[Th.~3.4]{L1}, or when $H \in {\cal{H}}^{\R}_n$, since, cf.~(\ref{mual}) below, in both cases the Ahlfors constants depend only on $n$. From (\ref{10}) we see that for any $t \in (0,10)$ and any point $z \in H$:
\begin{equation}\label{cnu}
\cs \{ p \in A^\bullet[i]\, | \, z \in B_{t \cdot R \cdot s^i}(p)\} \le c.
\end{equation}
3.~After scaling the balls in ${\cal{A}}^\bullet[i]$ by $(\ve \cdot s^i/\ve)^{-1}=s^{-i}$, $B_{s^i}(p)$ transforms
to $B_1(p)$ and $B_{R \cdot s^i}(p)$ becomes $B_R(p)$ with a tangent cone approximating with accuracy $\ve$ over $\TP_C(0,\ve,1, \ve) $ respectively over $\TP_C(0,\ve,R, \ve)$, since $R \le \ve^{-1}/10$. We observe that $\mathcal{R}_\ve(p) <\diam(H,d_{\sima})$, for any $p \in H$. Thus, the Ahlfors  regularity estimates (\ref{ahl})  apply to the balls in these covers.\qeda
\end{remark}

 \textbf{Proof}\, We inductively define two different selection maps $\kappa$ on $\Sigma^{*}$. In a first step we select $A \subset \Sigma^{*}$ for the small balls $B_{s^i}$ and then $A^\bullet \subset A$ for the large balls $B_{R \cdot s^i}$.\\

\textbf{Step 1 (Small Cover ${\cal{A}}$)}\,  We define $\kappa_i: \Sigma^{*}_i \ra \Z_2.$ When $\kappa_i(p)=0$, then we delete
$\overline{B_{s^{i}}(p)}$, otherwise, we keep it.  To start the (double) induction over $i$, we define $\kappa_i$ on the \emph{first non-empty} $\Sigma^*_i$:\\

\textbf{Start on $\Sigma^*_i$}:\,For $a_{i,1} \in \Sigma^*_i$, we set $\kappa_i(a_{i,1}):=1$.\\

 \textbf{Step on $\Sigma^*_i$}:\,We assume $\kappa_i$ has been defined for $a_{i,j}$, $j \le m$. Then we set
\[\kappa_i(a_{i,m+1}):=0,\, \mm{ if }\, a_{i,m+1} \in \textstyle \bigcup_{j \le m} \{\overline{B_{s^i}(a_{i,j})} \,|\,
a_{i,j} \in \Sigma^*_i,\kappa_i(a_{i,j}) = 1\}, \mm{ otherwise } \kappa_i(a_{i,m+1}):=1.\]

\textbf{Start on $\Sigma^*_{i+1}$}:\,Assuming $\kappa_a$ on $\Sigma^*_a$ has been defined for any $a \le i$ we continue with the definition of $\kappa_{i+1}$ on $\Sigma^*_{i+1}$. We start with $a_{i+1,1} \in \Sigma^*_{i+1}$  and set $\kappa_{i+1}(a_{i+1,1}):=0$ when
\[a_{i+1,1} \in  \textstyle \bigcup_{l \le i} \{\overline{B_{s^l}(a_{l,j})} \,|\, a_{l,j} \in \Sigma^*_l,\kappa_l(a_{l,j}) =1 \}\,\mm{ and } \kappa_{i+1}(a_{i+1,1}):=1  \mm{ otherwise.}\]

\textbf{Step on $\Sigma^*_{i+1}$}:\,We assume $\kappa_{i+1}$ has been defined for $a_{i+1,j}$, $j \le m$. When
\begin{talign*}
 \textstyle a_{i+1,m+1} \in & \bigcup_{j \le m} \{\overline{B_{s^{i+1}}(a_{i+1,j})} \,|\, a_{i+1,j} \in \Sigma^*_{i+1},\kappa_{i+1}(a_{i+1,j}) =1 \} \, \cup\\
  &\bigcup_{l \le i} \{\overline{B_{s^l}(a_{l,j})} \,|\,
a_{l,j} \in \Sigma^*_l,\kappa_l(a_{l,j}) = 1\}, \mm{ we set } \kappa_{i+1}(a_{i+1,m+1}):=0
\end{talign*}
and otherwise, $\kappa_{i+1}(a_{i+1,m+1}):=1$.
 Now we define  for $i \ge 0$ the following sets of balls and of center points of these balls

 \begin{itemize}
   \item $\textstyle  {\cal{A}}[i] := \{\overline{B_{s^i}(p)}\,|\, p \in \Sigma^*_i, \kappa_i(p)=1\},\, {\cal{A}} := \bigcup_{i\ge 0} {\cal{A}}[i],$
   \item $\textstyle A[i]:= \{p \in \Sigma^*_i\,|\, \kappa_i(p)=1\},\, A := \bigcup_{i\ge  0} A[i].$
 \end{itemize}
We observe that $\overline{B_{s^i/3}(x)} \cap \overline{B_{s^j/3}(y)} \v$, for $x \neq y, x \in A[i], y \in A[j]$. From this the Ahlfors $n$-regularity, \cite[Remark 3.7]{L1}, of $(H,d_H,\mu_H)$ on $H \in {\cal{G}}_n$ shows that in a given ball $B \subset H$ there are at most \emph{finitely} many balls belonging to ${\cal{A}}[i]$, for any given $i \ge 0$. We infer that
\begin{equation}\label{abi}
\textstyle \bigcup_{i \le k} {\cal{A}}[i]\mm{ covers }\bigcup_{i \le k} \Sigma_i \mm{ for any } k \ge  0.
\end{equation}
Otherwise, there were a $q \in \Sigma_l \setminus   \bigcup_{p \in \bigcup_{i \le k} A[i]} \overline{B_{s^i}(p)}$, for some $l \le k$. The local finiteness of each ${\cal{A}}[i]$ shows that $\bigcup_{p \in \bigcup_{i \le k} A[i]} \overline{B_{s^i}(p)}$ is a closed subset of $H$ and, thus, $q$ belongs to the \emph{open} complement.
Since $\Sigma^*_l \subset \Sigma_l$ is dense, there is some $q'= a_{l,m}\in  \Sigma^*_l \setminus \bigcup_{p \in \bigcup_{i \le k} A[i]} \overline{B_{s^i}(p)}$. For $\kappa[l]_1(a_{l,m})=0$, $a_{l,m}\in \bigcup_{p \in  A[l]} \overline{B_{s^i}(p)} \subset \bigcup_{p \in  \bigcup_{i \le k} A[i]} \overline{B_{s^i}(p)}$. Alternatively, for $\kappa[l]_1(a_{l,m})=1$, we have $a_{l,m}\in \overline{B_{s^i}(a_{l,m})} \subset  \bigcup_{p \in \bigcup_{i \le k} A[i]} \overline{B_{s^i}(p)}$. That is, both cases lead to contradictions.\\

\textbf{Step 2 (Large Cover ${\cal{A}}^\bullet$)}\,  Now we select subfamilies ${\cal{A}}^\bullet[i]$ of $\{\overline{B_{R \cdot s^i}(p)} \subset H\,|\, p \in A[i]\}$ with (\ref{22}) and  (\ref{10}). We define another selection map $\kappa_i^\bullet: A[i] \ra \Z^{\ge 0}$, writing the elements of $A[i]$ as $p_{i,1}, p_{i,2},..$.   When $\kappa_i^\bullet(p)=0$, then we delete $\overline{B_{R \cdot s^{i}}(p)}$, otherwise, we keep it.\\

\textbf{Start}:\,For $p_{i,1} \in A[i]$, we set $\kappa_i^\bullet(p_{i,1}):=1$.\\

 \textbf{Step}:\,We assume $\kappa_i^\bullet$ has been defined for $p_{i,j}$, $j \le m$.  We set
\[\kappa_i^\bullet(p_{i,j+1}):=0, \mm{ when } \overline{B_{s^i}(p_{i,j+1})} \subset \textstyle \bigcup_{k \le j} \{\overline{B_{R \cdot s^i}(p_{i,k})} \,|\, \kappa_i^\bullet(p_{i,k}) \ge 1 \}, \mm{ and otherwise }\]
 \[\kappa_i^\bullet(p_{i,j+1}):= \min \big(\{k \le j\,|\, \kappa_i^\bullet(p_{i,k})\ge 1 \mm{ with } d(p_{i,j+1},p_{i,k}) >  20 \cdot R \cdot s^i\} \cup \{j+1\}\big).\]
For $j \ge 1$ we define  the following sets of balls and of center points of these balls
\begin{itemize}
  \item $\textstyle  {\cal{A}}^\bullet[i,j] := \{\overline{B_{R \cdot s^i}(p)}\,|\, p \in A[i], \kappa_i^\bullet(p)=j\},\, {\cal{A}}^\bullet := \bigcup_{i\ge 0} {\cal{A}}^\bullet[i]:= \bigcup_{i\ge 0,j\ge 1} {\cal{A}}^\bullet[i,j],$
  \item $\textstyle A^\bullet[i,j]:= \{ p \in A[i]\,|\, \kappa_i^\bullet(p)=j\},\, A^\bullet := \bigcup_{i\ge 0} A^\bullet[i]:= \bigcup_{i\ge 0,j\ge 1} A^\bullet[i,j]$.
\end{itemize}
From the induction step we notice that $\textstyle \bigcup_{x \in A[i]} \overline{B_{s^i}(x)} \subset \bigcup_{z \in A^\bullet[i]} \overline{B_{R  \cdot s^i}(z)}$ and $z^* \notin \overline{B_{(R-1)  \cdot s^i}(z)}$, for $z \neq z^*, z,z^* \in A^\bullet[i]$. From Ahlfors $n$-regularity (\ref{ahl}) we have for any $z \in H$:
\begin{equation}\label{mual}
\mu_H(B_{R/4  \cdot s^i}(z))\ge  A \cdot (R/4  \cdot s^i)^n\mm{ and } \mu_H(B_{20 \cdot R \cdot s^i}(z)) \le B \cdot  (20 \cdot R \cdot s^i)^n,
\end{equation}
Now we set $c:=  [\mm{\emph{the smallest integer}}  \ge 100^n \cdot B/A]$ and we claim  ${\cal{A}}^\bullet[i,j] \v \mm{ for } j > c.$\\

Otherwise, we had some $\overline{B_{R  \cdot s^i}(p)} \in {\cal{A}}^\bullet[i,c+1]$. Then there are at least $c$ different $p_m \in A^\bullet[i,j_m]$ with $B_{10 \cdot R  \cdot s^i}(p) \cap B_{10 \cdot R  \cdot s^i}(p_m) \n$. Since the $B_{R/4  \cdot s^i}(x_m)$ are pairwise disjoint we get $\mu_H(B_{20 \cdot R \cdot s^i}(p)) \ge  100^n \cdot B/A \cdot A \cdot (R/4  \cdot s^i)^n$.
But this contradicts the upper estimate for $\mu_H(B_{20 \cdot R \cdot s^i}(p))$, and, hence, we have ${\cal{A}}^\bullet[i,c+1] \v$. \qed

\subsubsection{Transfer of Local Bumps} \label{trapro}

We cannot transfer shields on cones $C$ via $\D$-maps to a general singular $H \in{\cal{G}}$ since they intersect $\Sigma_C$. Instead we transfer periodic bumps from cones and use them to reproduce shields on $H$. (But we use the shields on $C$ to derive estimates for shields on $H$ in \ref{luem} below.)
We consider $\bigcup_{v \in S} (Q^{n-k}_1+ v) \subset Q^{n-k}_{L,...,L} \subset  \R^{n-k}$, for  some $L >0$, and the periodic bump $\bs(\bigcup_{v \in S} (Q^{n-k}_1+ v))$ on  $\R^{n-k} \times C^k \in {\cal{SC}}_n$ for $\omega_0(n,L)$, $R_0(n,L)$, $r_0(n,L)$ with
$\omega \in (0,\omega_0)$, $R > R_0$, $r \in (0, r_0)$ and  some $\beta_\bullet(\omega,r,R) \in (0,\beta_{n-k} \cdot \eta_{n-k})$:
 \begin{equation}\label{aoop}
\textstyle \su \bs(\bigcup_{v \in S} (Q^{n-k}_1+ v)) \subset \delta^{-1}_{\bp^*}(\R^{\ge  \beta_\bullet}) \cap B_R(0) \subset\TP_{\R^{n-k} \times C^k}(0,\omega,R,r)
 \end{equation}

\begin{corollary}\emph{\textbf{(Periodic Bumps on $H$)}} \label{tran} For $H \in{\cal{G}}$, some $p \in \Sigma_H$ and a tangent cone $(\R^{n-k} \times C^k,0)$ in $p$ so that the $\D$-map over $\TP_{\R^{n-k} \times C^k}(0,\omega,R,r) \subset C^n$, for $\omega \in (0,\omega_0)$, $R > R_0$, $r \in (0, r_0)$, with $\Psi_\circ$ appropriately normalized,  satisfies, for some  $\ve >0$:
\begin{equation}\label{lal}
|\D - id_{C}|_{C^5(\TP_{C}(0,\omega,R,r))} \le \ve \mm{ and } |\Psi_H \circ \D / \Psi_\circ-1|_{C^{2,\alpha}(\TP_{C}(0,\omega,R,r) )} \le \ve.
\end{equation}
For some small  $\ve_0(n) \in (0,1/10)$,  $\ve \in [0,\ve_0]$, (\ref{aoop}) shows that $\Phi^H_\square:=\Phi_\square \circ \D^{-1}$ is well-defined and extends to $H \setminus \Sigma_H$ with
\begin{equation}\label{ll}
L_{H^n}(\Phi^H_\square\cdot  \Psi_H) \ge \lambda/16 \cdot \bp^2 \cdot \Phi^H_\square \cdot  \Psi_H \mm{ on  }H \setminus \Sigma_H.
\end{equation}
We call $\big(\Phi^H_\square \cdot \Psi_H\big)^{4/(n-2)} \cdot g_H$ a \textbf{periodic bump metric} and the pseudo metric $\big((\Phi^H_\square -1) \cdot \Psi_H\big)^{4/(n-2)} \cdot g_H$ a \textbf{periodic bump}.
\end{corollary}
\begin{remark} We have $\Phi^H_\square \in {\cal{C}}(H^n, \Lambda[n])$, $\su |\Phi^H_\square- 1| = \su (\bs \circ \D^{-1})$. When we add $\bs_H:=\bs \circ \D^{-1}$  to $\Psi_H^{4/(n-2)} \cdot g_H$ this may change the conformal class since $\D$ need not to be conformal.
Using the multiplicative local bump $\Phi^H_\square$ resolves this issue. The two resulting deformations coincide in the limit of $\ve \ra 0$ and we continue speaking of added local bumps but technically we use the multiplicative version.\qeda
\end{remark}

\textbf{Tightness of Bumps on $H$} \,
We use the Ahlfors $n$-regularity of $(H,d^\Lambda_{\sima},\mu^\Lambda_{\sima})$ to extend the tightness result \ref{lue2} to the $\D$-images of bump elements in $H \in {\cal{G}}_n$ by some perturbation argument when $H$ is close to a product tangent cone in the sense that (\ref{lal}) of Lemma~\ref{tran} holds for $\ve \in [0,\ve_1]$, for some $\ve_1[n,\ell] \in (0,\ve_0[n,\ell])$. In general $\Sigma_H$ has no local product structure. For a given product tangent cone
$\R^{n-k} \times C^k$  we can, however, use induced structures when $H$ is (locally) closely approximated by $\R^{n-k} \times C^k$ in a common ambient space $\R^{n+1}$.
\begin{equation}\label{ome2}
\Omega_\ell := Q_\ell^{n-k} \times \R^{k+1}, \mm{ matching (\ref{ome}) from }  Q_\ell^{n-k} \times C^k= \R^{n-k} \times C^k \cap Q_\ell^{n-k} \times \R^{k+1},
\end{equation}
\vspace{-0.7cm}
\begin{equation}\label{oar}
O^C_{\rho}:= \R^{n-k} \times (C^k \cap B^{k+1}_\rho(0))\mm{  and  } O^H_{\rho} :=H \cap \R^{n-k} \times  B^{k+1}_\rho(0),
\end{equation}
 keeping in mind that $O^H_{\rho}$ depends on the chosen approximation by $\R^{n-k} \times C^k$ in a common ambient space $\R^{n+1}$.
For $\ve>0$  small enough so that for $\vartheta_i[n,C] \in (0,1)$, $\varrho \in (\vartheta_1,\vartheta_2)$ and  $\zeta=1/4 \cdot \min \{|\varrho - \vartheta_1|,|\varrho - \vartheta_2|\}$
of \ref{lue2}, we have for $\Psi_\circ$ appropriately normalized:
\begin{equation}\label{dell}
\, |\D - id_{C}|_{C^5(\delta^{-1}_{\bp^*}(\R^{\ge \beta_\bullet}) \cap O^C_{\varrho+3 \cdot \zeta} \setminus \overline{O^C_{\varrho-3 \cdot \zeta}} \cap \Omega_{\ell+2})}  \le \ve,
\end{equation}
\vspace{-0.7cm}
\begin{equation}\label{eq}
\textstyle |\Psi_H \circ \D/\Psi_\circ-1|_{C^{2,\alpha}(\delta^{-1}_{\bp^*}(\R^{\ge \beta_\bullet}) \cap O^C_{\varrho+3 \cdot \zeta} \setminus \overline{O^C_{\varrho-3 \cdot \zeta}} \cap \Omega_{\ell+2})}  \le \ve,
\end{equation}
and after small adjustments of the $\D$-map near the boundary in $C^ {2,\alpha}$-norm  for small $\ve > 0$:
\begin{equation}\label{hp}
\D\big(\delta^{-1}_{\bp^*}(\R^{\ge \beta_\bullet}) \cap O^C_{\varrho+2 \cdot \zeta} \setminus \overline{O^C_{\varrho-2 \cdot \zeta}} \cap \Omega_{\ell+1}\big) = \delta^{-1}_{\bp^*}(\R^{\ge \beta_\bullet}) \cap O^H_{\varrho+2 \cdot \zeta} \setminus \overline{O^H_{\varrho-2 \cdot \zeta}} \cap \Omega_{\ell+1},
\end{equation}
As a counterpart of $\mathcal{T}^C_{\varrho}$ in \ref{lue2} we set $\mathcal{T}^H_{\varrho}:=\Omega_{\ell+1} \cap \p O^H_{\rho}$.  For almost every $\rho>0$,  $\Omega_{\ell+1} \cap O^H_{\rho}$ is a Caccioppoli set.
Since $\mathcal{T}^C_{\varrho} \subset \R^{n-k} \times C^k$ is the level set of the radial distance function $d$ with $|\nabla d|=1$,  $\mathcal{T}^H_{\varrho} \cap \delta^{-1}_{\bp^*}(\R^{\ge \beta_\bullet})$ is smooth for small $\ve >0$ and we have a smooth convergence
\begin{equation}\label{idtq}
\D^{-1}(\mathcal{T}^H_{\varrho} \cap \delta^{-1}_{\bp^*}(\R^{\ge \beta_\bullet})) \cap \Omega_{\ell+1} \ra \mathcal{T}^C_{\varrho}  \cap \delta^{-1}_{\bp^*}(\R^{\ge \beta_\bullet}) \cap \Omega_{\ell+1} \subset C^n,
\mm{ for } \ve \ra 0.
\end{equation}
In general, $\mathcal{T}^H_{\varrho}$ is not a minimal hypersurface. We use $\mathcal{T}^H_{\varrho}$ to prescribe an appropriate Plateau boundary data in the following extension of  \ref{lue2}  for $\zeta_\bullet:=102/100 \cdot \zeta$.

\begin{proposition} \emph{\textbf{(Tightness and $\D$-Maps)}}\label{luem}
There is an $\ve_1[n,\ell] \in (0,\ve_0[n.\ell])$ so that for $\ve \in [0,\ve_1]$ in \ref{tran}, $\beta \in [0,\beta^*]$, $\G^H_\beta \in {\cal{C}}[H^n, \Lambda]$ with
\begin{equation}\label{ask}
\su |\G^H_{\beta} -1| \cap \su |\F_\beta \circ \D^{-1} -1| \cap O^H_{\varrho+3 \cdot \zeta}\setminus \overline{O^H_{\varrho-3 \cdot \zeta}} \cap \Omega_{\ell+2} \v
\end{equation}
\begin{enumerate}
\item there is a Caccioppoli set $U^\beta[\ell,n,k,H^n] \subset H \cap \Omega$ with $\overline{O^H_{\varrho-\zeta_\bullet}} \cap \Omega_{\ell+1}\subset U^\beta\subset O^H_{\varrho+\zeta_\bullet} \cap \Omega_{\ell+1}$ so that $\mathcal{T}_H^{\beta}[\ell,n,k,H^n] :=\p U_H^\beta$  and $\mathcal{T}_H^{\beta}  \cap \Omega_\ell$ is area minimizing with $\Omega_\ell$-boundary value $\mathcal{T}^H_{\varrho}$, relative to $(H,d^\Lambda_{\sima},\mu^\Lambda_{\sima})$ associated to $(\G^H_{\beta} \cdot \F_\beta \circ \D^{-1} \cdot \Psi_H)^{4/(n-2)} \cdot g_{H^n}$,
\item for any such oriented minimal boundary $\mathcal{T}_\bullet :=\p U_\bullet \subset \overline{O^H_{\varrho+ 2 \cdot \zeta}} \setminus O^H_{\varrho- 2 \cdot \zeta}$, allowing intersections with the obstacle $\p (O^H_{\varrho+  2 \cdot \zeta} \setminus O^H_{\varrho-2 \cdot \zeta})$, we already have $\mathcal{T}_\bullet\subset O^H_{\varrho+\zeta_\bullet} \setminus \overline{O^H_{\varrho-\zeta_\bullet}}$. We write $\mathcal{P}_H[n-k,\ell,\G^H_\beta]$ for the class of all such Plateau solutions $\mathcal{T}_\bullet$.
\end{enumerate}
\end{proposition}
\textbf{Proof}\, We assume there is no such approximation threshold $\ve_1[n,\ell] >0$, then we have sequences $H_i \in {\cal{G}}_n$ and of singular minimal cones $\R^{n-k} \times C_i^k$, so that the pairs $H_i, \R^{n-k} \times C_i^k$  satisfy (\ref{dell}) and (\ref{eq}) for some $\ve_i \ra 0$, when $i \ra \infty$, violating (i). (As in  \ref{lue2}, claim (ii) follows from the argument for (i).) Since the cones subconverge to another such cone $\R^{n-k} \times C^k$, we may assume that $\R^{n-k} \times C_i^k=\R^{n-k} \times C^k=C^n$.
We consider an area minimizer $\mathcal{T}_{H_i}^{\beta}$ solving the obstacle problem for open Caccioppoli sets $U_i$ with $O^{H_i}_{\varrho-2 \cdot \zeta} \cap \Omega_{\ell+1}\subset U_i \cap \Omega_{\ell+1} \subset O^{H_i}_{\varrho+ 2 \cdot \zeta} \cap \Omega_{\ell+1}$ so that $\p U_i$ has $\Omega_{\ell}$-boundary value $\mathcal{T}^{H_i}_{\varrho}$.
As in the cone case \ref{lue2} we can assume that $\p \Omega_\ell \cap O^{H_i}_{\varrho+ 3 \cdot \zeta}$ is locally outer minimizing. Namely, for $i$ large enough, we can replace the intersection of the hyperplanes (in the ambient space) with $O^{H_i}_{\varrho+3 \cdot \zeta} \subset H_i$ by Hausdorff approximating area minimizing hypersurfaces in $H_i$, where we use the hyperplanes to define the boundary value. This follows from Cor.\ref{grr} (\ref{est}) using a tightness argument as in \ref{lue2}. In the following we set $\Omega=\Omega_\ell$.\\

\textbf{Step 1 (Using $C$ for Estimates in $H_i$)}\, The first step is to use $\mathcal{T}_C^{\beta} \subset \R^{n-k} \times C^k$ to upper estimate the area of
$\mathcal{T}_{H_i}^{\beta}$. For such a comparison, we use $\D$-maps to embed $\mathcal{T}_C^{\beta}$ and $\mathcal{T}_{H_i}^{\beta}$ into the same space $H_i$.
These maps are, however, controllably defined only away from the singular set, typically on $\delta^{-1}_{\bp^*}(\R^{\ge \alpha})$, for some $\alpha \in (0,\beta)$. We therefore seek for a way to extend
$\D(\mathcal{T}_C^{\beta} \cap \delta^{-1}_{\bp^*}(\R^{\ge \alpha}) \cap \Omega)$ as a boundary of an open set in $H_i$ so that the area of the extension tends to zero for  $\alpha \ra 0$. One candidate would be a suitable subset of $\p(\delta^{-1}_{\bp^*}(\R^{\le \alpha}) \cap \Omega \cap O^{H_i}_{\varrho+ 2 \cdot \zeta} \setminus \overline{O^{H_i}_{\varrho-2 \cdot \zeta}})$, but it is difficult to estimate the area of this set.\\
Our workaround is to first choose finite ball covers  $B_{r_k}(p_k) \subset C$, $p_k \in \Sigma_C \cap \Omega \cap O^C_{\varrho+ 2 \cdot \zeta} \setminus \overline{O^C_{\varrho-2 \cdot \zeta}}$, $r_k \in (0,\zeta/4)$, $k \in I$, of $\Sigma_C \cap \Omega \cap O^C_{\varrho+ 2 \cdot \zeta} \setminus \overline{O^C_{\varrho-2 \cdot \zeta}}$ as in \ref{ctrl2} so that
\begin{equation}\label{covha2}
\textstyle  \sum_{k \in I}\mu_{\sima}^{\Lambda, n-1}(\p B_{r_k}(p_k))   \le c \cdot \sum_{k \in I} r_k^{n-1}  < \xi.
\end{equation}
For any ball $B^C_r(q) \subset B_1(0) \cap C$ of radius $r \in (0,1)$ in $(C,d_{\sima})$, there is a (non-unique) \emph{corresponding} ball $B^{H_i}_r(q_i)\subset (H_i,d_{\sima})$ so that
$B^{H_i}_r(q_i) \ra B^C_r(q)$ in flat and Hausdorff norm for $i \ra \infty$, cf.~\cite[Lemma 2.25]{L1}.\\ From \ref{ctrl2} we may assume (after generic changes of the radii in $[r(q_i), 2 \cdot r(q_i)]$ that (\ref{covha2})  also  holds for the corresponding balls, up to a factor $\le 2^{n-1}$. We note that the Hausdorff convergence shows that for $i$ large enough the $B^{H_i}_r(q_i)$ cover $\Sigma_{H_i} \cap \Omega \cap O^{H_i}_{\varrho+ 2 \cdot \zeta} \setminus \overline{O^{H_i}_{\varrho-2 \cdot \zeta}}$ and we even have, for $\alpha(\xi)>0$ small enough, in $C$ and any $H_i$ for $i$ large enough,
\begin{equation}\label{inb}
\textstyle \delta^{-1}_{\bp^*}(\R^{\le \alpha}) \cap O_{\varrho+\zeta} \setminus \overline{O_{\varrho-\zeta}}  \cap \Omega\subset \bigcup_{k \in I} B^C_{r_k}(p_k) \mm{ resp. } B^{H_i}_{r_k}(p_k).
\end{equation}
From this we can define $\xi$-\emph{ball extensions}: for an open set $V \subset \Omega \cap O_{\varrho+ 2 \cdot \zeta} \setminus \overline{O_{\varrho-2 \cdot \zeta}}$ we get another open set
$V_\xi \subset \Omega \cap O_{\varrho+ 2 \cdot \zeta} \setminus \overline{O_{\varrho-2 \cdot \zeta}}$ with $V \setminus \bigcup_{k \in I} B^C_{r_k}(p_k)=V_\xi \setminus  \bigcup_{k \in I} B^C_{r_k}(p_k)$ so that
\begin{equation}\label{vbx}
\textstyle \p V_\xi \cap  \bigcup_{k \in I} \overline{B^C_{r_k}(p_k)} \subset \bigcup_{k \in I}\p B^C_{r_k}(p_k) \mm{ and }  \mu_{\sima}^{\Lambda, n-1}(\p V_\xi \cap \bigcup_{k \in I} \p B_{r_k}(p_k))  < \xi
\end{equation}
We call $V_\xi$ an \textbf{$\boldsymbol\xi$-ball extension} of $V \setminus \bigcup_{k \in I} \overline{B^C_{r_k}(p_k)} $. Since $\mathcal{T}_{H_i}^{\beta} \cap \Omega$ is area minimizing we get
\begin{equation}\label{bbb}
\mu^{\Lambda, n-1}_{\sima,{H_i}}[\beta,\G^{H_i}_\beta]\big(\mathcal{T}_{H_i}^{\beta} \cap \Omega\big) \le \xi + (1+\ve_i)^{n-1} \cdot  \mu^{\Lambda, n-1}_{\sima,C}[\beta,\G^C_\beta](\mathcal{T}_C^\beta \cap \Omega),
\end{equation}
where the right hand side estimates the area of $\D(\mathcal{T}_C^{\beta} \cap \delta^{-1}_{\bp^*}(\R^{\ge \alpha}) \cap \Omega)$ plus a $\xi$-ball extension and $\G^{H_i}_\beta \in {\cal{C}}[{H_i}, \Lambda]$ and $\G^C_\beta \in {\cal{C}}[C, \Lambda]$ denote placeholder functions with $\G^{H_i}_\beta \circ \D=\G^C_\beta$ on $\delta^{-1}_{\bp^*}(\R^{\ge \alpha}) \cap O_{\varrho+ 2 \cdot \zeta} \setminus \overline{O_{\varrho-2 \cdot \zeta}} \cap \Omega$. For $i \ra \infty$ and $\xi \ra 0$ this shows
\begin{equation}\label{fl}
\limsup_{i \ra \infty} \mu^{\Lambda, n-1}_{\sima,{H_i}}[\beta,\G^{H_i}_\beta]\big(\mathcal{T}_{H_i}^{\beta} \cap \Omega\big)  \le  \mu^{\Lambda, n-1}_{\sima,C}[\beta,\G^C_\beta](\mathcal{T}_C^{\beta}  \cap \Omega).
\end{equation}

\textbf{Step 2 (Compact Convergence in $C \setminus \Sigma_C$)}\, Now we exchange the r\^{o}les of $H_i$ and $C$ and consider $\xi$-ball extensions of $\D^{-1}$-map images of the $\mathcal{T}_{H_i}^{\beta} \cap \delta^{-1}_{\bp^*}(\R^{\ge \alpha}) \cap \Omega$ in $C$ and write the open sets they bound as $W_i$. From (\ref{fl}) we get a flat norm subconvergence to some area minimizer, we write it again as $\mathcal{T}_C^{\beta}  \cap \Omega$ bounding $U_C$. We claim this implies compact Hausdorff convergence on $C \setminus \Sigma_C$. Assume that there is a ball $B_r(p) \in C \setminus \Sigma_C$ so that $B_r(p) \cap \mathcal{T}_C^{\beta} \v$ and $p \in W_i \Delta  U_C$ for all $i$. Then we consider the $\D$-map image where $\mathcal{T}_{H_i}^{\beta}$ is an area minimizer to argue from Cor.\:\ref{grr} (\ref{est}) that $\mu^{\Lambda}_{\sima,C}(W_i \Delta  U_C)>c$, for some $c>0$ independent of $i$, contradicting the flat norm convergence.\\

\textbf{Step 3 (Growth Estimates in  $H_i$)} We return to $H_i$ to derive (i) from  a contradiction. From step 2 we know after applying the $\D$-map that
\begin{equation}\label{par}
\mathcal{T}_{H_i}^{\beta} \cap \delta^{-1}_{\bp^*}(\R^{\ge \alpha}) \cap \Omega \ra \D(\mathcal{T}_C^{\beta} \cap \delta^{-1}_{\bp^*}(\R^{\ge \alpha}) \cap \Omega) \mm{ for } i \ra \infty
\end{equation}
in Hausdorff-norm. We consider $\xi$-ball extensions of $\mathcal{T}_{H_i}^{\beta} \cap \delta^{-1}_{\bp^*}(\R^{\ge \alpha}) \cap \Omega$ and assume that for any $i$ there is some $p_i \in \mathcal{T}_{H_i}^{\beta} \cap \delta^{-1}_{\bp^*}(\R^{\le \alpha}) \cap \Omega \cap \p (O^{H_i}_{\varrho+\zeta^*}\setminus O^{H_i}_{\varrho-\zeta^*})$, for $\zeta^*=101/100 \cdot \zeta$. Then we get from Corollary \ref{grr} (\ref{est}): $\mu^{\Lambda, n-1}_{\sima}[\beta,\G_\beta](\mathcal{T}_{H_i}^{\beta} \cap B_{1/100 \cdot \zeta}(p_i))>c$ for some $c>0$ independent of $i$.
This contradicts the Hausdorff closeness of $\mathcal{T}_{H_i}^{\beta} \cap \delta^{-1}_{\bp^*}(\R^{\ge \alpha}) \cap \Omega$ to $\D(\mathcal{T}_C^{\beta} \cap \delta^{-1}_{\bp^*}(\R^{\ge \alpha}) \cap \Omega)$ for large $i$ from (\ref{par}),  for sufficiently small $\xi >0$. \qed

To formulate the generalization of Cor.\:\ref{perg} we think of  $H$ as being  \emph{locally} closely approximated by $\R^{n-k} \times C^k$ in a common ambient space $\R^{n+1}$. We consider $\mathbb{O}^{\eta_q, \ell_q}_{q,\xi} \subset \R^{n-k} \times C^k$ defined from (\ref{oo}). $\mathbb{O}^{\eta_q, \ell_q}_{q,\xi}$ is union of trimmed cylinders we can write as intersections of trimmed cylinders $Z^{n+1}_f \subset \R^{n+1}$ surrounding the faces $f$ of  $Q_1^{n-k} \times \{0\} \subset \R^{n-k} \times C^k \subset \R^{n+1}$:
\begin{equation}\label{zy}
T[\eta_q ,f]\big(O_\xi \cap \Omega_q\big)=:Z^{n+1}_f \cap (\R^{n-k} \times C^k) \subset \R^{n-k} \times C^k
\end{equation}
We use the intersection of $Z^{n+1}_f $ with $H$ to define the \textbf{core} of the \textbf{periodic bump} in $H$ we
get for  $\bigcup_{v \in S} (Q^{n-k}_1+ v) \subset Q^{n-k}_{L,...,L} \subset  \R^{n-k}$ in $\R^{n-k} \times C^k$.
\begin{equation}\label{zy2}
\textstyle O^{\eta_q,\ell_q,H}_{f,q,\xi} \oplus v:=(Z^{n+1}_f +v)\cap  H \subset H,\,\, \mathbb{O}^H_{n-k}(\bigcup_{v \in S} (Q^{n-k}_1+ v)):=\bigcup_{v \in S}\bigcup_{f\in F_q} O^{\eta_q,\ell_q,H}_{f,q,r[n]}\oplus v
\end{equation}
where we assume an accuracy  $\ve \in [0,\ve_1]$, as in \ref{tran} and \ref{luem}. We transfer the original construction on cones for each (q)-face viewed
as a subset of $\R^q \times (\R^{n-k-q} \times C^k)$ with the assumptions and the notations of Cor.\:\ref{luem},  where we redefine $\zeta_\bullet:=1021/1000 \cdot \zeta$.

\begin{corollary}\emph{\textbf{(Periodic Shields on $H$)}}\label{omv} There is an $\ve_2(n)>0$ so that for an approximation with $\ve \in [0,\ve_2]$ as in (\ref{tran}),  $\G \in {\cal{C}}(H^n, \Lambda[n])$  with \emph{$\su |\G -1| \cap \su |\Phi^H_\square -1|  \v$} and  the $\ell_q, \eta_q,\beta_q$ from \ref{oma}  we have for any (q)-face $f\in F_q$ of  $Q_1^{n-k} \times \{0\}$ and $v \in S$ some
\begin{equation}\label{upen0}
\mathcal{T}_{f+v}  = \p U_{f+v}  \in \mathcal{P}_H[q,\ell_q,\G_{\beta_q}] \mm{ for an open } O^{\eta_q,\ell_q,H}_{f,q,r[n]} \oplus v \subset  U_{f+v}^H  \subset O^{\eta_q,\ell_q,H}_{f,q,1} \oplus v,
\end{equation}
where $\mathcal{T}_{f+v} $ is an area minimizer of \ref{luem} for the $\D$-map image of a bump element assigned to the face $f+v$ for some $\G_{\beta_q}\in {\cal{C}}(C^n, \Lambda[n])$ with  $\G_{\beta_q}:=\G \circ \D \circ (T[\eta_q ,f]+v)$ on $\delta^{-1}_{\bp^*}(\R^{\ge \beta_\bullet}) \cap O_{\varrho+ 2 \cdot \zeta} \setminus \overline{O_{\varrho-2 \cdot \zeta}} \cap \Omega_q$. From these $U^H_f$ we get a \textbf{periodic shield}:
\begin{equation}\label{open0}
\textstyle \mathbb{U}:=\bigcup^{n-k}_{q=0} \mathbb{U}^{\beta_q,\eta_q,\ell_q}_{q}:= \bigcup^{n-k}_{q=0} \bigcup_{f\in F_q} U^H_f,
\end{equation}
with $\mathbb{U} \supset \mathbb{O}^H_{n-k}(\bigcup_{v \in S} (Q^{n-k}_1+ v))  \supset  \Sigma_H \cap \mathbb{O}^H_{n-k}(\bigcup_{v \in S} (Q^{n-k}_1+ v))$ so that
 \begin{equation}\label{ibm0}
 \p \mathbb{U} \mm{  is a local \textbf{inner area minimizer} relative to } \mathbb{U} \mm{ and }(\G \cdot \Phi^H_\square \cdot \Psi_\circ )^{4/n-2} \cdot g_{H}.
 \end{equation}
For any $a\in \Z^{\ge 0}$ there is some $\ve^*_2(n,a) \in (0,\ve_2(n)]$ so that after rescaling the periodic bump on $\R^{n-k} \times C^k$,  described in \ref{perg} and \ref{perg2}, we get (up to uniform factors for the radii, when compared to the cone case,  converging to $1$, for $\ve \ra 0$, we omit for the sake of readability):
 \begin{enumerate}
\item $\D \big(U_{2^{-a}}(\bigcup_{v \in S} (Q^{n-k}_1+ v)) \setminus U_{2^{-a} \cdot r^\circ[n]}(\R^{n-k} \times \{0\})\big) \supset \su \bs_H(\bigcup_{v \in S} (Q^{n-k}_1+ v))$,
    \item $\{ z \in C \setminus \Sigma_C \,|\, 2 ^{-a} \ge \delta_{\bp}(z) \ge 2 ^{-a}\cdot r^\circ[n] \}  \supset \su \bs[a](\bigcup_{v \in S} (Q^{n-k}_1+ v))$,
\item $\mathbb{U}[a]\big(\bigcup_{v \in S} (Q^{n-k}_1+ v)\big) \supset \mathbb{O}^H_{n-k}[a]\big(\bigcup_{v \in S} (Q^{n-k}_1+ v)\big)  \supset  \Sigma_H \cap \mathbb{O}^H_{n-k}[a](\bigcup_{v \in S} (Q^{n-k}_1+ v))$,
\end{enumerate}
where $\mathbb{U}[a]$ and  $\mathbb{O}^H_{n-k}[a]$ are the counterparts of the rescaled shields and shielded cores in \ref{perg} for the scaled faces  (q)-face $f\in F_q$ of  $Q_1^{n-k} \times \{0\}$
and their $\D$-map transfer to $H$.
\end{corollary}

\textbf{Proof}\,  We define $\ve_2(n)$ to be the minimum of  the $\ve_1[n,\ell_q]$, which we multiply by $\eta_q$ to compensate the scaling of ($q$)-bump elements in Cor.\:\ref{luem}.  With this choice, the same bump element-wise argument as in Prop.\:\ref{oma}, based on the tightness result \ref{luem}, applies.\qeda

\begin{remark}\label{alr} \textbf{(Almost Tangent Cones)}\, In the results \ref{tran} - \ref{omv} on the transfer and evaluation of periodic bumps we only mentioned the
typical case where $(\R^{n-k} \times C^k,0)$ is a tangent cone in a singular point $p \in H$. We can slightly relax the coupling to $H$ since we only used the presence of a quantitatively controlled approximation by $\TP_{\R^{n-k} \times C^k}(0,\omega,R,r) \subset (\R^{n-k} \times C^k,0)$. The results remain valid provides the $\D$-map over $\TP_{\R^{n-k} \times C^k}(0,\omega,R,r)$ satisfies the conditions of (\ref{tran}) for some $\ve \in [0,\ve_2]$. In this case we call $(\R^{n-k} \times C^k,0)$ an $\boldsymbol{\ve}$\textbf{-almost tangent cone}. \qeda
\end{remark}

\subsubsection{Global Bumps and Alignments} \label{ifp}

When assembling a global bump we use self-similar ball covers to suitably place almost periodic bumps along the singular set $\Sigma_H \subset H$.
An essential observation, used to transfer the bumps, is that singular directions in any two tangent cones associated to neighboring balls align in a simple way. We start with a model situation that will be the limit case we use to understand fine tangent cone approximations for self-similar covers of high accuracy.

\begin{lemma}\emph{\textbf{(Linear Auto-Alignments)}}\label{comp1}\, For  $H \in {\cal{H}}^{\R}_n$ and $p_1,..,p_m \in \Sigma_H$, let $C_{p_1},..C_{p_m}$ be tangent cones with basepoints $p_1,..,p_m$ so that  $H=C_{p_1}=....C_{p_m} \subset \R^{n+1}$.  When the affine span of $\{p_1,..,p_m\}$ is  $q$-dimensional we have (up to a common rotation and translation)
\begin{equation}\label{oj1}
H=\R^q \times C^{n-q}\mm{ and } p_1,..,p_m \in \R^q \times \{0\}, \mm{ for some cone } C^{n-q} \in {\cal{SC}}_{n-q},
\end{equation}
in particular, we note that $q \le n-7$.
\end{lemma}

\textbf{Proof}\, We may assume that $p_1=0$. $C_{p_1}$ is scaling invariant and, thus, from scaling $C_{p_1}$ by $\lambda>0$ around $0$  we infer that for $C_{p_2}$ with basepoint $p_2 \neq 0$, using $H=C_{p_1}$ in the first equality and $H=C_{p_2}$ in the second and last one:
\[H= \lambda \cdot H= \lambda \cdot C_{p_2} = C_{p_2} - (1- \lambda) \cdot p_2 = H - (1- \lambda) \cdot p_2.\]
From this we see that $H$ is translation invariant in $p_2$-direction. After suitable rotation and translation this means $H^n = \R \times C^{n-1}$, with $p_1,p_2 \in \R \times \{0\}$, for some cone $C^{n-1} \in {\cal{SC}}_{n-1}$. For $p_3,..,p_m$ we continue inductively and note that the step is non-trivial only when the newly added point does not belong to the span of its predecessors.
\qed

We can rewrite (\ref{oj1}) as $B_1(p_k) \subset O_1^{\R^q \times C^{n-q}} \subset \R^q \times C^{n-q}$, for any $k=1,...m$. In the following we show that a completely similar inclusion locally holds for any layer ${\cal{A}}[i]$  in a self-similar cover of some sufficiently high accuracy.

\begin{proposition}\emph{\textbf{(Auto-Alignments)}}\label{comp2} For $H \in {\cal{G}}^c_n$, $R>10^2$ and $\eta \in (0,\ve_2)$ there is an $\ve_3(n,\eta,R) \in (0,\ve_2)$, so that for any self-similar cover of accuracy $\ve \in (0,\ve_3]$, $s \in (0,s_0)$ and $p  \in A^\bullet$ there is an $\eta$-almost tangent cone $C_p=\R^m \times C^{n-m}$, $m \in \Z^{\ge 0}$ in $p$ with
\begin{equation}\label{oj2}
 B_1(z) \subset O^{s^{-i} \cdot H}_2\mm{ for any } B_1(z) \subset B_R(p) \subset  s^{-i} \cdot H \mm{ with } z  \in A[i]  \mm{ when } p  \in A^\bullet[i].
\end{equation}
In general, and even in a fixed point $p$, $m$ is not uniquely determined.
\end{proposition}

\textbf{Proof}\, Otherwise we had  some $\eta \in (0,\ve_2)$ and a sequence $\ve_k \ra 0$, for  $k \ra \infty$, so that there are self-similar covers of accuracy $\ve_k$ and points $p_k \in A^\bullet[i_k]$, so that
\begin{equation}\label{oj3}
\textstyle \bigcup_{z  \in A[i_k], B_1(z) \subset B_{R}(p_k)} B_{1}(z) \nsubseteq  O^{s_{i_k}^{-1} \cdot H}_{2}.
\end{equation}
for any $\eta$-almost tangent cone $C$ in $p_k$ recalling from (\ref{oar}) that $O_\rho^H$ is defined relative to $C$.
\begin{itemize}[leftmargin=*]
\item  The number of balls $B_1(q_k)$ with $B_1(q_k) \subset B_{R}(p_k) \subset s^{-i_k} \cdot H$ and $q_k \in A[i_k]$ is uniformly upper estimated depending only on $n$ and $R$ from the Ahlfors $n$-regularity of  $(H,d_H,\mu_H)$, giving lower and upper volume bounds for the $B_{s^{i_k}/2}(q_k), B_{R \cdot s^{i_k}}(p_k) \subset H$ as in Prop.\ref{brop} since $B_{s^{i_k}/2}(q_k) \cap B_{s^{i_k}/2}(q^*_k)\v$, for $q_k \neq q^*_k \in  A[i_k]$.
\item Since $H$ is compact, we observe that $i_k \ra \infty$, for $k \ra \infty$. After selecting a subsequence the pairs
    \[\big((1/s^{i_k}\cdot H,p[\ve_k]),(1/s^{i_k}\cdot  (\overline{B_{(R-1) \cdot s^{i_k}}(p_k)} \cap A[i_k])\big)\] of pointed spaces and of finite subsets, of upper bounded cardinality, converge to some limit $\big((H_\infty,p_\infty),\overline{B_{R-1}(p_\infty)} \cap A_\infty\big)$, where $H_\infty \in {\cal{H}}^{\R}_n$ and $\overline{B_{R-1}(p_\infty)} \cap  A_\infty \subset \Sigma_{H_\infty}$, from Allard theory. In turn, for $k \ra \infty$, the refining  $\ve_k$-accuracy and the compactness of the cone space ${\cal{SC}}_{n}$ show that  $H_\infty$ actually is a cone $C_\infty$ and any sequence of $\ve_k$-accurately approximating tangent cones $C_{p_k}$, with basepoint $p_k$, subconverges to $C_\infty$. From this $C_\infty$ is an $\eta_k$-almost tangent cone of $(1/s^{i_k}\cdot H,p_k)$, for some $\eta_k \ra 0$, when $k \ra \infty$.
\item  For further subsequences there are $\ve_k$-accurately approximating tangent cones in the points $1/s^{i_k}\cdot (\overline{B_{(R-1) \cdot s^{i_k}}(p_k)} \cap A[i_k])$ converging to tangent cones of $C_\infty$ in corresponding points of $A_\infty$. Each of the limit cones coincides with $C_\infty$ since $\ve_k \ra 0$. That is, we reach a configuration as in \ref{comp1} and thus we have $C_\infty=\R^m \times C^{n-m}$, for some $m \in \Z^{\ge 0}$ and $\overline{B_{R-1}(p_\infty)} \cap A_\infty \subset O^{\R^m \times C^{n-m}}_1$. \\ For large enough $k$, $C_\infty$ is an $\eta_k$-almost tangent cone with $\eta_k \le  \min \{ \eta, \ve_2\}$ and the convergence of  points in $A[i_k]$ to corresponding points in $A_\infty$  eventually contradicts  (\ref{oj3})  for $(1/s^{i_k}\cdot H,p_k)$. This proves the claim for $H \in {\cal{G}}^c_n$, that is, we get an   $\ve_3(H,\eta,R)$ as asserted.
\item In the  argument we do not need that the hypothetical  sequence of points $p_k$ belongs to a fixed $H \in {\cal{G}}^c_n$, that is, after an appropriate adjustment of the $i_k$, we even get an $\ve_3(n,\eta,R)$ independent of $H$.\qeda
\end{itemize}

\begin{theorem} \emph{\textbf{(Global Bumps)}}\label{gbb} For any $H \in {\cal{G}}_n^c$ and any neighborhood $W$ of $\Sigma_H$, there is a finite family of almost periodic bumps, disjointly supported in $W$ so that their sum $\bx$ is a \textbf{global bump}: the associated metric  $\big(\Phi^H_{\bx} \cdot \Psi_H\big)^{4/(n-2)} \cdot g_H$  has $\boldsymbol{\scal>0}$, the union $\mathbb{U}(\bx)$ of the associated almost periodic shields covers $\Sigma_H$ and  $\p \mathbb{U}(\bx)$ is a local \textbf{inner area minimizer}.
\end{theorem}

 \textbf{Proof}\, We apply the auto-alignment to each layer ${\cal{A}}^\bullet[i]$ of a self-similar cover separately and place periodic bumps disjointly shielding all balls in ${\cal{A}}[i]$.
 The balls and, thus, the assigned bumps in one sublayer ${\cal{A}}^\bullet[i,j]$, i.e. for a fixed $j$, do not intersect. For intersecting balls in different of the at most $c(n)$ sublayers we can use the rescaled bumps from  Cor.\ref{omv}(i) and (ii) to ensure that the support of the bumps does not intersect, cf.~Fig.\:\ref{fig:ali}. This argument equally applies for any $i$ and we only need to choose $s>0$ small enough, again independent of $i$, to ensure that the balls in ${\cal{A}}^\bullet[k]$, for $k>i$, do not intersect the support of the periodic bumps for ${\cal{A}}^\bullet[l]$, for $l \le i$.

{\begin{figure}[htbp]
\centering
\includegraphics[width=0.7\textwidth]{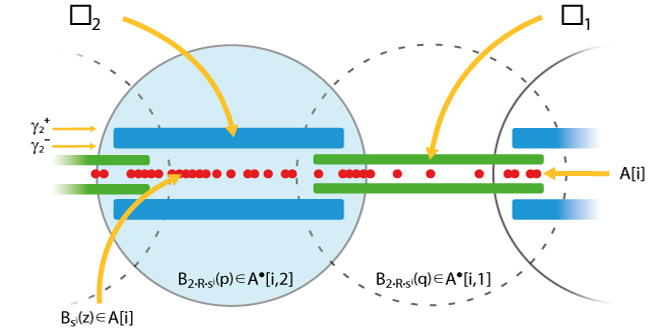}
\caption{\footnotesize Disjoint placement of periodic bumps $\bs_1$ and $\bs_2$ for intersecting balls in different sublayers of ${\cal{A}}^\bullet[i]$: ${\cal{A}}^\bullet[i,1]$ \emph{= dashed lines} and ${\cal{A}}^\bullet[i,2]$ \emph{= straight lines}. \normalsize}
\label{fig:ali}
\end{figure}}

\textbf{Step 1 (Prepared Cones)}\, From  Cor.\ref{perg}, Remark \ref{perg2}, Example \ref{perg3} and suitable rescaling there is a large $R_n\ge 10^2$ so that for $B^n_{R_n}(0) \subset \R^{n-k} \times C^k=C^n$,  for any $k \ge 7$, $C^n \in{\cal{SC}}_n$ we have $c(n)$  periodic bumps $\bs_b$ supported in $B^n_{3 \cdot R_n}(0)$, for $b=1,..,c$, with
\begin{equation}\label{bed}
\su \bs_b \subset  \{z \in B^n_{3 \cdot R_n}(0) \,|\, \gamma^+_b \ge \delta_{\bp}(z) \ge \gamma^-_b \}
\end{equation}
for $\gamma_b^\pm=2^{m_b^\pm}$, for some $m_b^\pm \in \Z^{\ge 1}$, depending only on $n$,  with $\gamma^+_b> \gamma^-_b\ge 4$ and $\gamma^-_{b+1} \ge  2 \cdot \gamma^+_b$,   and so that
each of the $\bs_b$ provides a shield $\mathbb{U}^C_b$ for $B^{n-k}_2(0)$, that is,
\begin{equation}\label{bed2}
\{z \in B^n_{3 \cdot R_n}(0) \,|\, \gamma^+_b \ge \delta_{\bp}(z)\} \supset \mathbb{U}^C_b \supset B^{n-k}_2(0).
\end{equation}
For the transfer to  $H \in {\cal{G}}_n^c$ we choose an $\ve_4 \in (0,\ve_3]$ so that $R_n \in [10^2,\ve^{-1}/10^2]$, for any $\ve \in (0,\ve_4]$, and we consider a self-similar cover  ${\cal{A}}^\bullet$ of accuracy  $\ve$.\\

\textbf{Step 2 (Transfer to $H$)}\, For any of the sublayers  ${\cal{A}}^\bullet[i,j]$ of ${\cal{A}}^\bullet[i]$, we consider the balls $\overline{B_{R_n \cdot s^i}(q)} \in {\cal{A}}^\bullet[i,j]$ and, after scaling by $s^{-i}$, we transfer $\bs_j$
via $\D$-map from some $\eta$-almost tangent cone as provided by Prop.\ref{comp2} to  $\overline{B_{3 \cdot R_n \cdot s^i}(q)}$. From Prop.\ref{omv} and  Prop.\ref{comp2} we can choose the accuracy $\ve  \in (0,\ve_3)$ so that for $\eta  \in (0,\ve_2)$ the transfer of the periodic bumps to $\overline{B_{3 \cdot R_n \cdot s^i}(q)}$ yields a periodic shield $\mathbb{U}^H_j \subset \overline{B_{3 \cdot R_n \cdot s^i}(q)}$ that contains all $\overline{B_{s^i}(z)} \subset \overline{B_{R_n \cdot s^i}(q)}$ with $z  \in A[i]$. We repeat this for any $j$ and any ball in ${\cal{A}}^\bullet[i,j]$ and, hence, in ${\cal{A}}^\bullet[i]$ and notice that the periodic shields are disjointly supported and contain all balls in ${\cal{A}}[i]$.\\

\textbf{Step 3 (Global Shield)}\, We repeat this bump placement for all layers ${\cal{A}}^\bullet[k]$ of ${\cal{A}}^\bullet$  with the same estimates, up to appropriate scaling, for the bumps and shields.
From estimate (\ref{bus}), which carries over to $H$ and $p \in \Sigma_H$, we observe that for $s>0$ small enough we get disjointly supported  periodic bumps.
We recall from Prop.\ref{brop}(ii) that $\bigcup_{i \le k} {\cal{A}}[i]$ covers $\bigcup_{i \le k} \Sigma_i$. From Cor.\ref{omv} each of the, at first infinitely many, periodic bumps we associated to the balls in ${\cal{A}}^\bullet$ shields a core independent of the placeholder function $\G \in{\cal{C}}[H,\Lambda]$ we use for that particular bump. These cores form an open cover of the compact set $\Sigma_H$. We take a finite subcover.\\
We define the global bump metric from adding the periodic bumps associated to the cores in this finite subcover to $(H,d_{\sima},\mu_{\sima})$. Now, we consider each of these periodic bumps $\bs$ separately, and get a shield $\U_{\bs}$ relative to the global bump metric. The union $\mathbb{U}(\bx)$ of all these $\U_{\bs}$ is bounded by locally inner minimizer $\p \mathbb{U}(\bx)$, from the same argument as in the proof of Prop.\ref{oma}. For suitably small $\ve>0$ we can also make sure that $\mathbb{U}(\bx) \subset W$. \qed

Finally, we show how \ref{gbb} implies our main splitting result.

\begin{theorem} \emph{\textbf{(Splitting with Boundary)}}\label{spth} Let $H^{n} \subset M^{n+1}$ be an almost minimizer $H \in {\cal{G}}^c$ with singular set $\Sigma_H$,  so that the conformal Laplacian $\bp^{-2} \cdot L_H$ has a positive principal eigenvalue  $\lambda^{\bp}_H>0$.  Then there are arbitrarily small neighborhoods $U$ of $\Sigma$ so that $H \setminus U$ is \textbf{conformal} to a $\boldsymbol{\scal>0}$\textbf{-manifold} $X_U$  with  \textbf{minimal}  boundary $\p X_U$.
\end{theorem}

This is Theorem 3 from the introduction. It contains Theorem 1 as a special case from  \cite[Theorem 2(i)]{L4}. For Theorem 2, where $H$ has a boundary $\p H \cap \Sigma \v$, we note that the assembly of a global bump and the proof of its shielding effect only use properties close to $\Sigma$. This leaves the argument unchanged, cf.~\cite[Rm.\:1.13]{L1} and \cite[Rm.\:3.10]{L4}, and, hence, the following argument equally applies to Theorem 2.\\

\textbf{Proof}\, For $\ve>0$ we consider the $10 \cdot \ve$-neighborhood $U_{10 \cdot \ve}$ of $\Sigma_H$ in $(H,d_{\sima},\mu_{\sima})$. Using smooth approximations of the distance function we may assume that $\p U_{10 \cdot \ve}$ is smooth. Now we make an auxiliary deformation that, however, becomes invisible towards the end of the argument. We deform $H$ in an $\ve$-distance tube of $\p U_{10 \cdot \ve}$ so that $\p U_{10 \cdot \ve}$ becomes positively mean curved relative to $U_{10 \cdot \ve}$. Then  $\p U_{10 \cdot \ve}$ is locally inner minimizing relative to $H \setminus U_{10 \cdot \ve}$. We claim that for any $\eta>0$ there is a neighborhood $V_\eta \subset U_{10 \cdot \ve}$ of $\Sigma_H$ with
\begin{equation}\label{2d}
\dist(\p V_\eta,\p U_{10 \cdot \ve})> 5 \ve \mm{  and  }\mu^{\Lambda, n-1}_{\sima}(\p V_\ve) \le  \eta.
\end{equation}
To check (\ref{2d}) we use Cor.\:\ref{ctrl2}: for any $\eta > 0$ and $r \in (0,1)$ there is some ball cover $B_{r_i}(p_i)$, $p_i \in \Sigma$, $i \in I$, of $\Sigma$ with $r_i \le r$ so that $\sum_{i \in I}\mu_{\sima}^{\Lambda, n-1}(\p B_{r_i}(p_i), d_{\sima}) \le b^\circ \cdot \sum_{i \in I} r_i^{n-1} < \eta$.
For sufficiently small $r \in (0, \ve)$, $V_\ve:=\bigcup_{i \in I}B_{r_i}(p_i)$ satisfies (\ref{2d}).\\
From \ref{gbb} we can choose a global bump $\bx$ supported in $V_\eta$. Since $\p \mathbb{U}(\bx)$ is locally inner minimizing  relative to $\mathbb{U}(\bx)$,
we get an open Caccioppoli set $W_\eta$ with  $\mathbb{U}(\bx) \subset W_\eta \subset U_{10 \cdot \ve}$ so that $\p W_\eta$ is area minimizing under all such $W$ and thus
\begin{equation}\label{w}
\mu^{\Lambda, n-1}_{\sima}(\p W_\eta) \le \mu^{n-1}_{\sima}(\p V_\eta) \le  \eta.
\end{equation}
In turn, we get from Cor.\:\ref{grr} (\ref{est}) some $a>0$, so that for any ball $B_\ve(p)$ for some $p \in  U_{10 \cdot \ve}$ with $\dist(p,\p U_{10 \cdot \ve}) =  3 \cdot \ve$ and any area minimizer $L^{n-1} \subset (H,d_{\sima},\mu_{\sima})$ passing through $p$ we have $\mu_{\sima}^{n-1}(L^{n-1} \cap B_\ve(p))=\mu_{\sima}^{\Lambda,n-1}(L^{n-1} \cap B_\ve(p))>a$ and, hence,
 for $\eta>0$ small enough, $\dist(\p W_\eta,\p U_{10 \cdot \ve}) \ge 2 \cdot \ve$.  That is, $W_\eta$ is disjoint from the support of the auxiliary deformation of $(H,d_{\sima},\mu_{\sima})$ in the $\ve$-neighborhood of $\p U_{10 \cdot \ve}$. Summarizing, we find for any $\ve>0$ some $\eta>0$ and some global bump $\bx$, with support in $V_\eta$,  so that  $\big(\Phi^H_{\bx} \cdot \Psi_H\big)^{4/(n-2)} \cdot g_H$ and $U:=W_\eta$ have the asserted properties.\qed

\footnotesize
\renewcommand{\refname}{\fontsize{14}{0}\selectfont 
\end{document}